\documentclass[11pt]{article}
\usepackage{amsfonts}
\usepackage{mathrsfs}
\usepackage{latexsym,amssymb,graphics}
\usepackage{amsmath,amssymb,epsfig,subfigure,url,color}
\usepackage{graphicx}
\usepackage{subfigure}
\usepackage{exscale}
\usepackage{relsize}
\usepackage{amsmath}
\usepackage{txfonts}
\usepackage{url}
\usepackage[flushleft]{caption2}
\usepackage{multirow} \usepackage{makecell}
\usepackage{multicol}
\usepackage[numbers,sort&compress]{natbib}
\date{}

\begin{document}
\title{New explicit thresholding/shrinkage formulas for one class of regularization problems with overlapping group sparsity and their applications\thanks{The work of Gang Liu, Ting-Zhu Huang and Jun Liu is supported by 973 Program (2013CB329404), NSFC (61370147), Sichuan Province Sci. \& Tech. Research Project (2012GZX0080). The work of Xiao-Guang Lv is supported by Postdoctoral Research Funds (2013M540454, 1301064B).}}
\author{Gang Liu\thanks{School of Mathematical Sciences, University of Electronic Science and Technology of China, Chengdu, Sichuan, 611731, P. R. China ({\em wd5577@163.com}).},
Ting-Zhu Huang\thanks{School of Mathematical Sciences, University of Electronic Science and Technology of China, Chengdu, Sichuan, 611731, P. R. China ({\em tingzhuhuang@126.com}).},
Xiao-Guang Lv\thanks{School of Mathematical Sciences, Nanjing Normal University, Nanjing, Jiangsu, 210097, P. R. China ({\em xiaoguanglv@126.com}).},
Jun Liu\thanks{School of Mathematical Sciences, University of Electronic Science and Technology of China, Chengdu, Sichuan, 611731, P. R. China ({\em junliucd@163.com}).}}

\maketitle

\begin{abstract} The least-square regression problems or inverse problems have been widely studied in many fields such as compressive sensing, signal processing, and image processing. To solve this kind of ill-posed problems, a regularization term (i.e., regularizer) should be introduced, under the assumption that the solutions have some specific properties, such as sparsity and group sparsity. Widely used regularizers include the $\ell_1$ norm, total variation (TV) semi-norm, and so on.
 Recently, a new regularization term with overlapping group sparsity has been considered. Majorization minimization iteration method or variable duplication methods are often applied to solve them. However, there have been no direct methods for solve the relevant problems because of the difficulty of overlapping. In this paper, we proposed new explicit shrinkage formulas for one class of these relevant problems, whose regularization terms have translation invariant overlapping groups. Moreover, we apply our results in TV deblurring and denoising with overlapping group sparsity. We use alternating direction method of multipliers to iterate solve it. Numerical results also verify the validity and effectiveness of our new explicit shrinkage formulas.
\end{abstract}

\section{Introduction}

The least-square regression problems or inverse problems have been widely studied in many fields such as compressive sensing, signal processing, image processing, statistics and machine learning. Regularization terms with sparse representations (for instance the $\ell_1$ norm regularizer) have been developed into an important tool
in these applications recently \cite{DYZ2013,HZ2010,CW2013,SPH2009}.
These methods are based on the assumption that signals or images have a sparse representation, that is, only containing a few nonzero entries.
To further improve the solutions, more recent studies suggested to go beyond sparsity and took into account additional information about the
underlying structure of the solutions \cite{DYZ2013,JOV2009,HZ2010}. Particularly, a wide class of solutions which with specific ※group sparsity§ structure are considered. In this case, a group sparse vector can be divided into groups of components satisfying a) only a few of groups contain nonzero values and b) these groups are not needed to be sparse. This property sometimes calls ``joint sparsity'' that a set of sparse vectors with the union of their supports being sparse \cite{CW2013}.
Many literature had consider this new sparse problems \cite{DYZ2013,JOV2009,HZ2010,SPH2009,CW2013,MW2009}. Putting
such group vectors into a matrix as row vectors of the marix, this matrix will only have
few nonzero rows and these rows may be not sparse.
These problems are typically obtained
by replacing the problem
\begin{equation}\label{norm11}
  \min_{\mathbf{z}} \ \ \ \|\mathbf{z}\|_{1} + \frac{\beta}{2}\|\mathbf{z} -\mathbf{x}\|_2^2,\ \
\end{equation}
with
\begin{equation}\label{norm21nop}
  \min_{\mathbf{z}} \ \ \ \sum_{i=1}^r \|\mathbf{z}[i]\|_{2}+ \frac{\beta}{2}\|\mathbf{z} -\mathbf{x}\|_2^2,
\end{equation}
where $\mathbf{x}\in \mathbb{R}^n$ is a given vector, $\mathbf{z}\in \mathbb{R}^n$, $\|\mathbf{z}\|_{p}=(\sum_{i=1}^n |z_i|^p)^{\frac{1}{p}}$ with $(p=1,2)$ represents the $\ell_p$ norm of vector $z$, $|z_i|$ is the absolute value of $z_i$,
and $\mathbf{z}[i]$ is the $i$th group of $z$ with  $\mathbf{z}[i]\cap \mathbf{z}[j] = \oslash$ and  $\bigcup_{i=1}^r \mathbf{z}[i]= \mathbf{z}$. The first term of former equations is called the regularization term, the second term is called the fidelity term, and $\beta>0$ is the regularization parameters.

Group sparsity solutions have better representation and have been widely studied both for convex and nonconvex cases \cite{DYZ2013,MW2009,EV2011,GSC2012,LWZ2011,CW2013}. More recently, overlapping group sparsity (OGS) had been considered \cite{DYZ2013,CW2013,LHS2013,LHLL2013,SC2013,FB2011,JAB2011,K2009,PF2011,BJMO2011,B2011}. These methods are based on the assumption that signals or image have a special sparse representation with OGS. The task is to solve the following problem
\begin{equation}\label{norm21ogsp}
  \min_{\mathbf{z}} \ \ \ \|\mathbf{z}\|_{2,1}+ \frac{\beta}{2}\|\mathbf{z} -\mathbf{x}\|_2^2,
\end{equation}
where $\|\mathbf{z}\|_{2,1} = \sum_{i=1}^n \|(z_i)_g\|_2$ is the generalized $\ell_{2,1}$-norm. Here, each $(z_i)_g$ is a group vector containing $s$ (called group size) elements that surrounding the $i$th entry of $z$. For example, $(z_i)_g = (z_{i-1},z_{i},z_{i+1})$ with $s=3$. In this case, $(z_i)_g$, $(z_{i+1})_g$and $(z_{i+2})_g$ contain the $(i+1)$th entry of $z$ ($z_{i+1}$), which means overlapping different from the form of group sparsity (\ref{norm21nop}). Particularly, if $s=1$, the generalized $\ell_{2,1}$-norm degenerates into the original $\ell_{1}$-norm, and the relevant regularization problem (\ref{norm21ogsp}) degenerates to (\ref{norm11}).

To be more general, we consider the weighted generalized $\ell_{2,1}$-norm $\|z\|_{w,2,1} = \sum_{i=1}^n \|w_g\circ(z_i)_g\|_2$ (we only consider that each group has the same weight, which means translation invariant) instead the former generalized $\ell_{2,1}$-norm, the task can be extended to
\begin{equation}\label{norm21wogsp}
  \min_{\mathbf{z}} \ \ \ \|\mathbf{z}\|_{w,2,1}+ \frac{\beta}{2}\|\mathbf{z} -\mathbf{x}\|_2^2,
\end{equation}
where $w_g$ is a nonnegative real vector with the same size as $(z_i)_g$ and ``$\circ$'' is the point-wise product or hadamard product. For instance, $w_g\circ(z_i)_g= ((w_g)_1 z_i,(w_g)_2 z_{i+1},(w_g)_3 z_{i+2})$ with $s=3$ as the former example. Particularly, the weighted generalized $\ell_{2,1}$-norm degenerates into the generalized $\ell_{2,1}$-norm if each entry of $w_g$ equals to 1.
The problems (\ref{norm21ogsp}) and (\ref{norm21wogsp}) had been considered in \cite{FB2011,JAB2011,K2009,PF2011,BJMO2011,B2011,DYZ2013}. They solve the relevant problems by using variable  duplication methods (variable  splitting,  latent/auxilliary  variables,  etc.). Particularly, Deng et. al in \cite{DYZ2013} introduced a diagonal matrix $G$ for this variable duplication methods. This matrix $G$ was not easy to find and would break the structure of the coefficient matrix, which made the difficulty of solving solutions under high dimensional vector cases. Moreover, it is difficult to extend this method to the matrix case.

Considering the matrix case of the problem (\ref{norm21wogsp}), we can get
\begin{equation}\label{norm21wogspm}
  \min_{A} \ \ \ \|A\|_{W,2,1}+ \frac{\beta}{2}\|A -X\|_F^2,
\end{equation}
where $X,A\in \mathbb{R}^{m\times n}$, $\|A\|_{W,2,1} = \sum_{i=1}^m\sum_{j=1}^n \|W_g\circ(A_{i,j})_g\|_F$. Here, each $(A_{i,j})_g$ is a group matrix containing $K_1\times K_2$ (called group size) elements that surrounding the $(i,j)$th entry of $A$. For example,
\begin{equation*}\label{sampKK}
W_g\circ(A_{i,j})_g =\left[\begin{array}{cccc}  (W_g)_{1,1}A_{i-l_1,j-l_2}  &(W_g)_{1,2}A_{i-l_1,j-l_2+1}  &\cdots &(W_g)_{1,K_2}A_{i-l_1,j+r_2}\\
                                         (W_g)_{2,1}A_{i-l_1+1,j-l_2}  &(W_g)_{2,2}A_{i-l_1+1,j-l_2+1}  &\cdots &(W_g)_{2,K_2}A_{i-l_1+1,j+r_2}\\
                                          \vdots    &\vdots       &\ddots &\vdots\\
                                         (W_g)_{K_1,1}A_{i+r_1,j-l_2}&(W_g)_{K_1,2}A_{i+r_1,j-l_1+1}&\cdots &(W_g)_{K_1,K_2}A_{i+r_1,j+r_2}\\ \end{array}\right]\in \mathbb{R}^{K_1\times K_2},
\end{equation*}
where $l_1=[\frac{K_1-1}{2}]$, $l_2=[\frac{K_2-1}{2}]$, $r_1=[\frac{K_1}{2}]$, $r_2=[\frac{K_2}{2}]$ (with $l_1+r_1+1=K_1$ and $l_2+r_2+1=K_2$) and [$x$] denotes the largest integer less than or equal to $x$. Particularly, if $K_2=1$, this problem degenerate to the former vector case (\ref{norm21wogsp}). If $K_1=K_2=1$, this problem degenerate to the original $\ell_1$ regularization problem (\ref{norm11}) for the matrix case. If $(W_g)_{i,j}\equiv 1$ for $i=1, \cdots, K_1$,  $j=1, \cdots, K_2$, this problem had been considered in Chen et. al \cite{CS2014}. However, they used an iterative algorithm based on the principle of majorization minimization (MM) to solve this problem.

In this paper, we propose new explicit shrinkage formulas for all the former problems (\ref{norm21ogsp}), (\ref{norm21wogsp}) and (\ref{norm21wogspm}), which can get accurate solutions without iteration in \cite{CS2014}, variable  duplication (variable  splitting,  latent/auxilliary  variables,  etc.) in \cite{FB2011,JAB2011,K2009,PF2011,BJMO2011,B2011}, or finding matrix $G$ in \cite{DYZ2013}. Numerical results will also verify the validity and effectiveness of our new explicit shrinkage formulas.
Moreover, this new method can be used as a subproblem in many other OGS problems such as compressive sensing with $\ell_1$ regularization and image restoration with total variation (TV) regularization. According to the framework of ADMM, the relevant convergence theory results of these OGS problems can be easy to be obtained because of the accurate solution of their subproblems with the application of our new explicit shrinkage formulas. For example, we will apply our results in image restoration using TV with OGS in our work, and we will obtain the convergence theorems. Numerical results will also verify the validity and effectiveness of our new methods.

The outline of the rest of this paper is as follows. In Section 2 we detailed deduce our explicit shrinkage formulas for OGS problems (\ref{norm21ogsp}), (\ref{norm21wogsp}) and (\ref{norm21wogspm}). In Section 3, we propose some extension for these shrinkage formulas. In Section 4, we apply our results in image deblurring and denoising problems with OGS TV. Numerical results are given in Section 5. Finally, we conclude this paper in Section 6.

\section{OGS-shrinkage}
\subsection{Original shrinkage}
For the original sparse represent solutions we often want to solve the following problems
\begin{equation}\label{norm1p}
  \min_{\mathbf{z}} \ \ \ \|\mathbf{z}\|_{p} + \frac{\beta}{2}\|\mathbf{z} -\mathbf{x}\|_2^2, \ \ p=1,2.
\end{equation}
\\
\textbf{Definition 1}.  Define shrinkage mappings $Sh_1$ and $Sh_2$ from $\mathbb{R}^N\times R^+$ to $\mathbb{R}^N$ by
\begin{equation}\label{shrinkage1}
  {Sh_1(\mathbf{x},\beta)}_i = {\rm sgn}(x_i) \max\{ |x_i| - \frac{1}{\beta},0\},
\end{equation}
\begin{equation}\label{shrinkage2}
  {Sh_2(\mathbf{x},\beta)} = \frac{\mathbf{x}}{\|\mathbf{x}\|_2} \max\{ \|\mathbf{x}\|_2 - \frac{1}{\beta},0\},
\end{equation}
where both expressions are taken to be zero when the second factor is zero, and ``sgn'' represents the signum function indicating the sign of a number, that is, sgn($x$)=0 if $x=0$, sgn($x$)=-1 if $x<0$ and sgn($x$)=1 if $x>0$.

The shrinkage (\ref{shrinkage1}) is known as soft thresholding and occurs in many algorithms related to sparsity since it is
the proximal mapping for the $\ell_1$ norm. Then the minimizer of (\ref{norm1p}) with $p=1$ is the following equation (\ref{shrinkage1g}).
\begin{equation}\label{shrinkage1g}
  \arg\min\limits_{\mathbf{z}} \ \ \ \|\mathbf{z}\|_{1} + \frac{\beta}{2}\|\mathbf{z} -\mathbf{x}\|_2^2 = {Sh_1(\mathbf{x},\beta)}.
\end{equation}
Thanks to the additivity and separability of both the $\ell_1$ norm and the square of the $\ell_2$ norm, the shrinkage (\ref{shrinkage1}) can be deduced easily by the following formula:
\begin{equation}\label{norm11p}
  \min_{\mathbf{z}} \ \ \ \|\mathbf{z}\|_{1} + \frac{\beta}{2}\|\mathbf{z} -\mathbf{x}\|_2^2 = \sum_{i=1}^n \min_{z_i}|z_i| + \frac{\beta}{2}|z_i-x_i|^2.
\end{equation}

The minimizer of (\ref{norm1p}) with $p=2$ is the following equation (\ref{shrinkage2g}).
\begin{equation}\label{shrinkage2g}
  \arg\min\limits_{\mathbf{z}} \ \ \ \|\mathbf{z}\|_{2} + \frac{\beta}{2}\|\mathbf{z} -\mathbf{x}\|_2^2 = {Sh_2(\mathbf{x},\beta)}.
\end{equation}
This formula is deduced by the Euler equation of (\ref{norm1p}) with $p=2$. Clearly,
\begin{equation}\label{shrinkage2gEuler}
 \beta\left( \mathbf{z} -\mathbf{x} \right)+ \frac{\mathbf{z}}{\|\mathbf{z}\|_{2}} \ni \mathbf{0},
\end{equation}
\begin{equation}\label{shrinkage2gEuler2}
\left(1 + \frac{1}{\beta}\frac{1}{\|\mathbf{z}\|_{2}} \right)\mathbf{z}-\mathbf{x} \ni \mathbf{0}.
\end{equation}
We can easily get that the necessary condition is that the vector $\mathbf{z}$ is parallel to the vector $\mathbf{x}$. That is $\frac{\mathbf{z}}{\|\mathbf{z}\|_{2}}=\frac{\mathbf{x}}{\|\mathbf{x}\|_{2}}$. Substitute into (\ref{shrinkage2gEuler2}), and the formula (\ref{shrinkage2}) is obtained.
More details please refer to \cite{WY2008,YYZW2009}.

Our new explicit OGS shrinkage formulas are based on these observations, especially the properties of additivity, separability and parallelity. \\
\textbf{Remark 1}. The problem (\ref{norm21nop}) is also easy to be solved by a simple shrinkage formula, which is not used in this work. More details refer to \cite{DYZ2013,WCT2012,SRSE2011}.
\subsection{The OGS shrinkage formulas}
Now we focus on the problem (\ref{norm21ogsp}) firstly. The difficulty of this problem is overlapping. Therefore, we must take some special techniques to avoid overlapping. That is the point of our new explicit OGS shrinkage formulas.

It is obvious that the first term of the problem (\ref{norm21ogsp}) is additive and separable. So if we find some relative rules such that the second term of the problem (\ref{norm21ogsp}) has the same properties with the same variable as the first term, the solution of (\ref{norm21ogsp}) can be easily found similar as (\ref{norm11p}).

Assuming period boundary conditions is used here, we observe that each entry $z_i$ of vector $\mathbf{z}$ would appear exactly $s$ times in the first term. Therefore, to hold on the uniformity of vectors $\mathbf{z}$ and $\mathbf{x}$, we need multiply the second term by $s$. To maintain the invariability of the problem (\ref{norm21ogsp}), and after some manipulations, we have
\begin{equation}\label{norm21ogspdo}
\begin{array}{r*{20}{l}}
  f_m(\mathbf{z})=&\min\limits_{\mathbf{z}} &\|\mathbf{z}\|_{2,1}+ \frac{\beta}{2}\|\mathbf{z} -\mathbf{x}\|_2^2\\
  =&\min\limits_{\mathbf{z}} &\sum_{i=1}^n \|(z_i)_g\|_2 + \frac{\beta}{2s}s\|\mathbf{z} -\mathbf{x}\|_2^2\\
  =&\min\limits_{\mathbf{z}} &\sum_{i=1}^n \|(z_i)_g\|_2 + \frac{\beta}{2s}\sum_{i=1}^n\|(z_i)_g -(x_i)_g\|_2^2\\
    =&\min\limits_{\mathbf{z}} &\sum_{i=1}^n \left(\|(z_i)_g\|_2 + \frac{\beta}{2s}\|(z_i)_g -(x_i)_g\|_2^2\right),\\
\end{array}
\end{equation}
where $(x_i)_g$ is same as $(z_i)_g$ defined before.

\begin{figure}
  \centering
\subfigure[Vector]{\includegraphics[width=0.4\textwidth,clip]{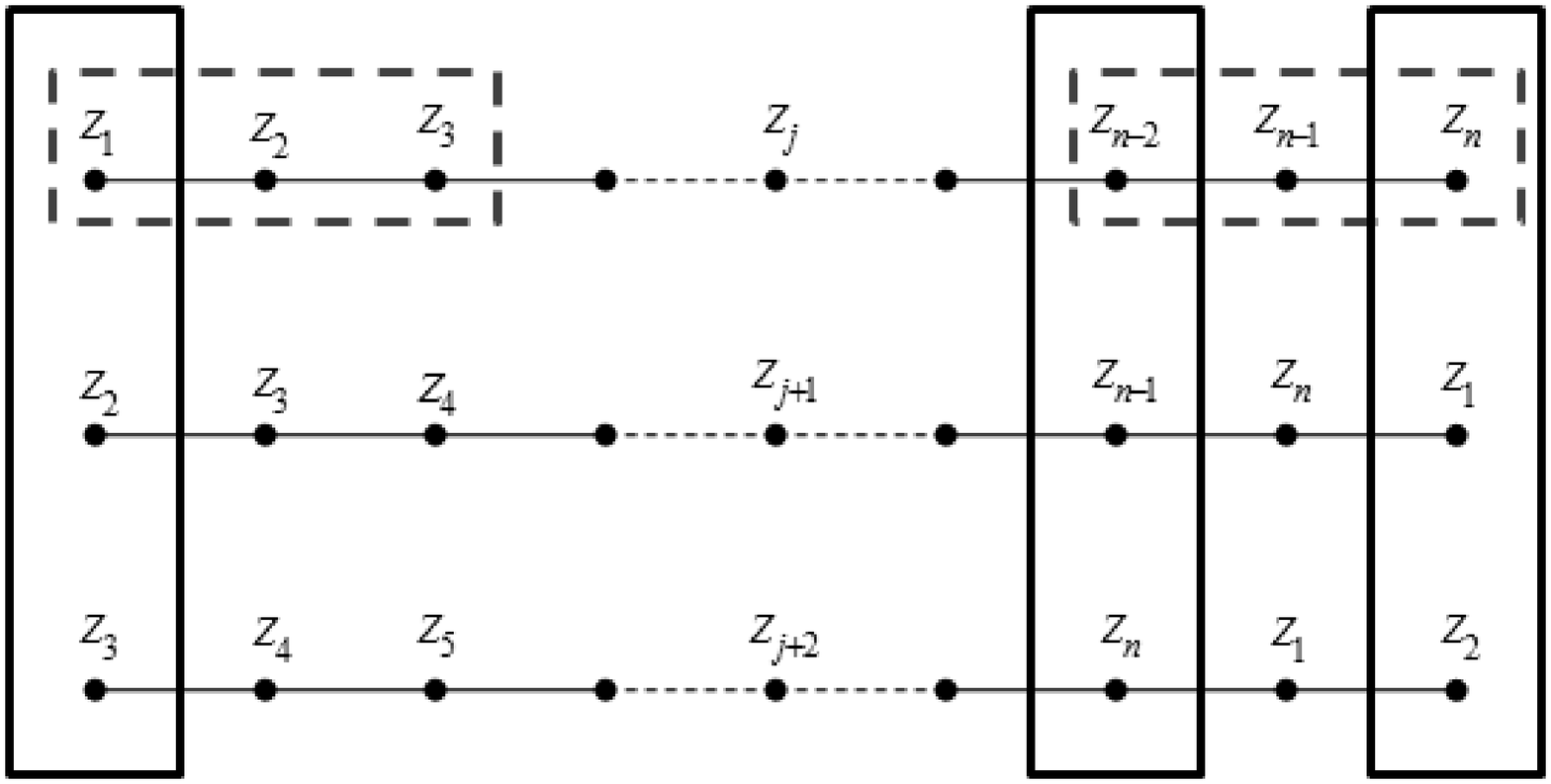}}
  \subfigure[weighted]{\includegraphics[width=0.4\textwidth,clip]{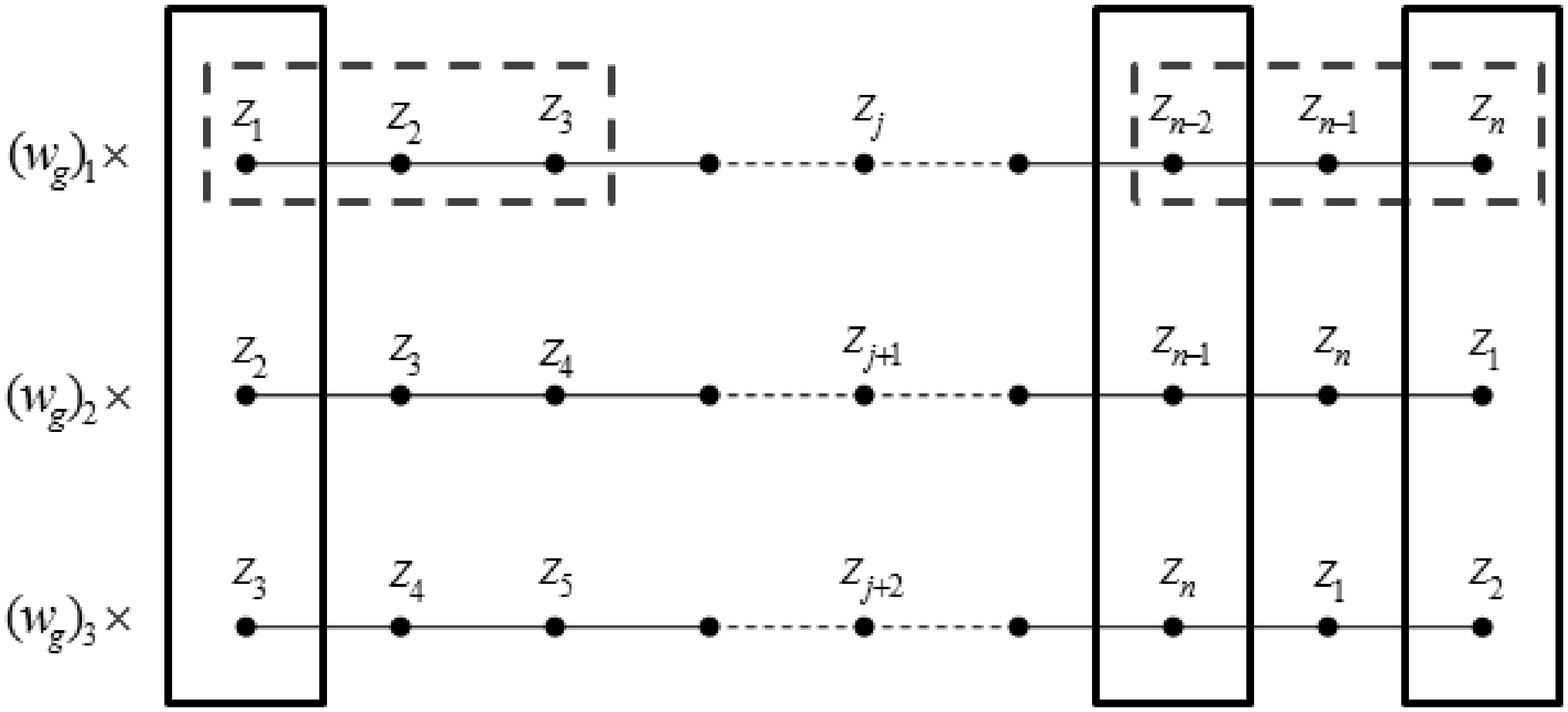}}
  \caption{Vector case.}
  \label{vector}
 \end{figure}

For example, we set $s=3$ and define $(z_i)_g=(z_i,z_{i+1},z_{i+2})$. The generalized $\ell_{2,1}$-norm $\|\mathbf{z}\|_{2,1}$ can be treated as the generalized $\ell_1$ norm of generalized points, whose entry $(z_i)_g$ is also a vector, and the absolute value of each entry is treated as the $\ell_2$ norm of $(z_i)_g$. See Figure~\ref{vector}(a) intuitively, where the top line is the vector $\mathbf{z}$, the rectangles with dashed line are original $(z_i)_g$, and the rectangles with solid line are the
generalized points. Because of the period boundary conditions, we know that each line of Figure~\ref{vector}(a) is translated equal. We treat the vector $\mathbf{x}$ same as the vector $\mathbf{z}$. Putting these generalized points (rectangles with solid line in the figure) as the columns of a matrix, then $s\|\mathbf{z} -\mathbf{x}\|_2^2$ can be regarded as matrix Frobenius norm $\|\left((z_i)_g\right)-\left((x_i)_g\right)\|_F$ in Figure~\ref{vector}(a) with every line being the row of the matrix. This is why the second equality in (\ref{norm21ogspdo}) holds. Therefore, generally, for each $i$ of the last line of (\ref{norm21ogspdo}), from the equation (\ref{shrinkage2}) and (\ref{shrinkage2g}), we can obtain
\begin{equation*}
  \arg\min\limits_{(z_i)_g} \|(z_i)_g\|_2 + \frac{\beta}{2s}\|(z_i)_g -(x_i)_g\|_2^2 = {Sh_2((x_i)_g,\frac{\beta}{s})},
\end{equation*}then
\begin{equation}\label{subshrinkage2g}
  {(z_i)_g} = \max \left\{ \|{(x_i)_g}\|_2 - \frac{s}{\beta}, 0\right\} \frac{{(x_i)_g}}{\|{(x_i)_g}\|_2}, \ \ \  {(z_i)_g} = \left([{(z_i)_g}]_1,[{(z_i)_g}]_2,...,[{(z_i)_g}]_s\right).
\end{equation}

Similarly as Figure~\ref{vector}(a), for each $i$, the $i$th entry $x_i$ (or $z_i$) of the vector $\mathbf{x}$ (or $\mathbf{z}$) may appear $s$ times, so we need compute each $z_i$ for $s$ times in $s$ different groups.

However, the results from (\ref{subshrinkage2g}) are wrong, because the results $z_i$ in $s$ different groups are different from (\ref{subshrinkage2g}). That means the results are not able to be satisfied simultaneously in this way. Moreover, for each $i$ of the last line in (\ref{norm21ogspdo}), the result (\ref{subshrinkage2g}) is given by that the vector $(z_i)_g$ is parallel to the vector $(x_i)_g$. Notice this point and ignore that $(z_i)_g=\mathbf{0}$ or $(x_i)_g=\mathbf{0}$, particularly for $s=4$ and $(z_i)_g=(z_{i-1},z_{i},z_{i+1},z_{i+2})$, the vector $\mathbf{z}$ can be split as follows,
\begin{equation}\label{subshrinkagenew1}
  \begin{array}{rr*{10}l}
  \mathbf{z}=&&(z_1,&z_2,&z_3,&z_4,&z_5,&\cdots,&z_{n-2},&z_{n-1},&z_{n})\\
=&+\frac{1}{4}&(z_1,&z_2,&z_3,&0,&0,&\cdots,&0,&0,&z_{n}).\\
&\frac{1}{4}&(z_1,&z_2,&z_3,&z_4,&0,&\cdots,&0,&0,&0)\\
&+\frac{1}{4}&(0,&z_2,&z_3,&z_4,&z_5,&\cdots,&0,&0,&0)\\
&+&&\cdots&&&&&&&\\
&+\frac{1}{4}&(z_1,&0,&0,&0,&0,&\cdots,&z_{n-2},&z_{n-1},&z_{n})\\
&+\frac{1}{4}&(z_1,&z_2,&0,&0,&0,&\cdots,&0,&z_{n-1},&z_{n})\\
  \end{array}
\end{equation}
Let $(z_i)'_g = (0,\cdots,0,z_{i-1},z_i,z_{i+1},z_{i+2},0,\cdots,0)$ be the expansion of $(z_i)_g$, with
$(z_1)'_g = (z_1,z_{2},z_{3},0,$ $\cdots,0,z_{n})$, $(z_{n-1})'_g = (z_1,0,0,\cdots,0,z_{n-2},z_{n-1},z_{n})$, $(z_{n})'_g = (z_1,z_2,0,\cdots,0,z_{n-1},z_{n})$. Let $(x_i)'_g = (0,\cdots,0,x_{i-1},x_i,x_{i+1},x_{i+2},0,\cdots,0)$ be the expansion of $(x_i)_g$ similarly as $(z_i)'_g$. Then, we have $\mathbf{z}=\frac{1}{4}\sum_{i=1}^{n}(z_i)'_g$,
and $\mathbf{x}=\frac{1}{4}\sum_{i=1}^{n}(x_i)'_g$. Moreover, we can easily obtain that $\|(z_i)'_g\|_2=\|(z_i)_g\|_2$ and $\|(x_i)'_g\|_2=\|(x_i)_g\|_2$ for every $i$.

On one hand, the Euler equation of $f_m(\mathbf{z})$ (with $s=4$) is given by

\begin{equation}\label{subshrinkagenew2222}
  \beta\left( \mathbf{z}- \mathbf{x}\right) + \frac{(z_1)'_g}{\|(z_1)'_g\|_2} + \cdots + \frac{(z'_n)_g}{\|(z'_n)_g\|_2}\\
  \ni \mathbf{0},
\end{equation}
\begin{equation}\label{subshrinkagenewnew}
\frac{\beta}{4}\sum_{i=1}^{n}\left((z_i)'_g - (x_i)'_g \right) + \frac{(z_1)'_g}{\|(z_1)'_g\|_2} +\cdots + \frac{(z_i)'_g}{\|(z_i)'_g\|_2} + \cdots + \frac{(z_n)'_g}{\|(z_n)'_g\|_2}\ni \mathbf{0}.
\end{equation}
From the deduction of the 2-dimensional shrinkage formula (\ref{shrinkage2}) in Section 2.1, we know that the necessary condition of minimizing the $i$th term of the last line in (\ref{norm21ogspdo}) is that
$(z_i)_g$ is parallel to $(x_i)_g$. That is, $(z_i)'_g$ is parallel to $(x_i)'_g$ for every $i$. For example,
\begin{equation}\label{subshrinkagenew2}
 (z_2)'_g = (z_1,z_2,z_3,z_4,0,\cdots,0,0) // (x_2)'_g = (x_1,x_2,x_3,x_4,0,\cdots,0,0).
\end{equation}
Then we obtain  $\frac{(z_i)'_g}{\|(z_i)'_g\|_2} =\frac{(x_i)'_g}{\|(x_i)'_g\|_2}$.
Therefore, (\ref{subshrinkagenewnew}) changes to
\begin{equation}\label{subshrinkagenew2221}
\frac{\beta}{4}\sum_{i=1}^{n}\left((z_i)'_g - (x_i)'_g \right) + \frac{(x_1)'_g}{\|(x_1)'_g\|_2} +\cdots + \frac{(x_i)'_g}{\|(x_i)'_g\|_2}+ \cdots + \frac{(x_n)'_g}{\|(x_n)'_g\|_2}\ni \mathbf{0},
\end{equation}
\begin{equation}\label{subshrinkagenew22211}
 \beta\left( \mathbf{z}- \mathbf{x}\right)+ \frac{(x_1)'_g}{\|(x_1)'_g\|_2} +\cdots + \frac{(x_i)'_g}{\|(x_i)'_g\|_2}+ \cdots + \frac{(x_n)'_g}{\|(x_n)'_g\|_2}\ni \mathbf{0},
\end{equation}

\begin{equation}\label{subshrinkagenew22211new}
\mathbf{z} \ni \mathbf{x}-\frac{ 1}{\beta} \left(\frac{(x_1)'_g}{\|(x_1)'_g\|_2} +\cdots + \frac{(x_i)'_g}{\|(x_i)'_g\|_2}+ \cdots + \frac{(x_n)'_g}{\|(x_n)'_g\|_2}\right),
\end{equation}
for each component, we obtained
\begin{equation}\label{subshrinkagenew22211newi}
z_i \ni x_i -\frac{ 1}{\beta} \left(\frac{x_i}{\|(x_{i-2})'_g\|_2} +\frac{x_i}{\|(x_{i-1})'_g\|_2} + \frac{x_i}{\|(x_i)'_g\|_2}+ \frac{x_i}{\|(x_{i+1})'_g\|_2}\right).
\end{equation}

Therefore, when $(z_i)_g\neq\mathbf{0}$ and $(x_i)_g\neq\mathbf{0}$,
we find a minimizer of (\ref{norm21ogsp}) on the direction that all the vectors $(z_i)_g$ are parallel
 to the vectors $(x_i)_g$.
 In addition, when $ \frac{\beta}{4}\left((z_i)'_g - (x_i)'_g \right) + \frac{(z_i)'_g}{\|(z_i)'_g\|_2} = 0$, (\ref{subshrinkagenewnew}) holds,
  then $\left(\frac{\beta}{4}+ \frac{1}{\|(z_i)'_g\|_2}\right)(z_i)'_g =(x_i)'_g$, therefore,
  $\frac{(z_i)'_g}{\|(z_i)'_g\|_2} =\frac{(x_i)'_g}{\|(x_i)'_g\|_2}$ holds.
 Moreover, because of the strict
 convexity of $f_m(\mathbf{z})$, we know that the minimizer is unique. This minimizer $\mathbf{z}$ is accurate.

On the other hand, when $(z_i)_g=\mathbf{0}$ or $(x_i)_g=\mathbf{0}$, our method may not obtain the accurate minimizer.
When $(x_i)_g=\mathbf{0}$, we know that the minimizer of the subproblem $\min\limits_{(z_i)_g}\|(z_i)_g\|_2 + \frac{\beta}{2s}\|(z_i)_g -(x_i)_g\|_2^2$ is exactly that $(z_i)_g=\mathbf{0}$. When $(x_i)_g\neq\mathbf{0}$ and the minimizer of the subproblem $\min\limits_{(z_i)_g}\|(z_i)_g\|_2 + \frac{\beta}{2s}\|(z_i)_g -(x_i)_g\|_2^2$ is that $(z_i)_g=\mathbf{0}$ (this is because that the parameter $\beta/s$ is two small), our method is not able to obtain the accurate minimizer. For example, that $(z_i)_g=\mathbf{0}$ while $(z_{i+1})_g\neq\mathbf{0}$ makes the element $z_i$ in $\mathbf{z}$ take different values in different subproblems. However, we can obtain an approximate minimizer in this case, which is that the element $z_i$ is a simple summation of corresponding subproblems containing $z_i$. We will show that in experiments of Section 5 the approximate minimizer is also good. Moreover, when we take this problem as a subproblem of the image processing problem, we can set the parameter $\beta$ to be large enough to make sure that the minimizer of the subproblem is accurate. Therefore, the convergence theorem results can be obtained by this accuracy, which will be applied in Section 4.

In addition, form (\ref{subshrinkagenew22211newi}), we can know that the element $z_i$ of the minimizer can be treated as in $s$ subproblems independently and then combine them. After some manipulations, in conclusion, we can get the following two general formulas.

1). \begin{equation}\label{shrinkageogs1}
 \arg\min\limits_{\mathbf{z}} \ \ \ \|\mathbf{z}\|_{2,1} + \frac{\beta}{2}\|\mathbf{z} -\mathbf{x}\|_2^2 = {Sh_{OGS}(\mathbf{x},\beta)},
\end{equation}
 with
\begin{equation}\label{shrinkageogs21}
  {Sh_{OGS}(\mathbf{x},\beta)}_i =z_i =\max \left\{ 1 - \frac{1}{\beta}{F(x_i)}, 0\right\} {x_i}.
\end{equation}
Here, for instance, when group size $s=4$, $(z_i)_g=(z_{i-1},z_{i},z_{i+1},z_{i+2})$ in $\|\mathbf{z}\|_{2,1}$ and $(x_i)_g$ is defined similarly as $(z_i)_g$, we have
$F(x_i) = \left(\frac{1}{\|(x_{i-2})_g\|_2} +\frac{1}{\|(x_{i-1})_g\|_2}+\frac{1}{\|(x_i)_g\|_2}+\frac{1}{\|(x_{i+1})_g\|_2}\right)$. The $\|(x_{j})_g\|_2$ is contained in $F(x_i)$ if and only if $\|(x_{j})_g\|_2$ has the component $x_i$, and we follow the convention $(1/0)= 1$ in $F(x_i)$ because $\|(x_i)_g\|_2=0$ implies $x_i=0$ and the value of $F(x_i)$ is insignificant in (\ref{shrinkageogs21}).

2). \begin{equation}\label{shrinkageogs}
 \arg\min\limits_{\mathbf{z}} \ \ \ \|\mathbf{z}\|_{2,1} + \frac{\beta}{2}\|\mathbf{z} -\mathbf{x}\|_2^2 = {Sh_{OGS}(\mathbf{x},\beta)},
\end{equation}
 with
\begin{equation}\label{shrinkageogs2}
  {Sh_{OGS}(\mathbf{x},\beta)}_i =z_i ={G(x_i)}\cdot{x_i}.
\end{equation}
Here, symbols are the same as 1), and
$G(x_i) = \max\left(\frac{1}{s} - \frac{1}{\beta\|(x_{i-2})_g\|_2},0\right)$ $ +\max\left(\frac{1}{s}- \frac{1}{\beta\|(x_{i-1})_g\|_2},0\right)$ \\$+\max\left(\frac{1}{s} - \frac{1}{\beta\|(x_i)_g\|_2},0\right)+\max\left(\frac{1}{s} - \frac{1}{\beta\|(x_{i+1})_g\|_2},0\right)$.

When $\beta$ is sufficiently large or sufficiently small, the former two formulas are the same and both are accurate. For the other values of $\beta$ , from the experiments, we find that 2) is better approximate than 1), so we choose the formula 2). Then, we obtain the following algorithm for finding the minimizer of (\ref{norm21ogsp}).\\
\begin{tabular}{l}
\hline
\hline
{\textsc{\textbf{Algorithm 1}}} \textup{Direct shrinkage algorithm for the minimization problem (\ref{norm21ogsp})}\\
\hline
 \textit{\textbf{Input}}:  \\ \ \ Given vector $\mathbf{x}=(x_1,x_2,\cdots,x_n)$, group size $s$, parameter $\beta$.
\\
 \textit{\textbf{Compute}}:\\ Definition of $(x_i)_g$, for example $(x_i)_g = \left(x_{i-s_l},\cdots,x_i,\cdots,x_{i+s_r}\right)$, \\
 with $s_l =[\frac{s-1}{2}]$ and $s_r =[\frac{s}{2}]$ ($s_l+s_r+1=s$). \\
 $\mathbf{w}$ = ones(1,s) = [1,1,$\cdots$,1].\\
If $s$ is even, then $\mathbf{w} = [\mathbf{w},0]$ and $s=s+1$. \\
Let $\mathbf{w}_r$ = fliplr($\mathbf{w}$) be fliping $\mathbf{w}$ over 180 degrees.\\
 Compute  $X_{n}=\left(\|(x_1)_g\|_2,\cdots,\|(x_i)_g\|_2,\cdots,\|(x_n)_g\|_2\right)$, by convolution of $\mathbf{w}$ and $\mathbf{x}$.\\
Compute $X'_{n}=\max\left( \frac{1}{s} - 1./(X_{n}\cdot\beta),0\right)$ pointwise. \\
Compute $G(x_i)$ by correlation of $\mathbf{w}$ and $X'_{n}$, or by convolution of $\mathbf{w}_r$ and $X'_{n}$.\\
\\
\hline
\end{tabular}\\

We can see that Algorithm 1 only need 2 times convolution computations with time complexity $n*s$, which is just the same time complexity
as one step iteration in the MM method in \cite{CS2014}. Therefore, our method is much more efficient than MM method or other variable  duplication methods. In Section 5, we will give the numerical experiments for comparison between our method and the MM method. Moreover, if $s=1$, our Algorithm 1 degenerates to the classic soft thresholding as our formula (\ref{shrinkageogs2}) degenerates to  (\ref{shrinkage1}). Moreover, when $\beta$ is sufficiently large or sufficiently small, our algorithm is accurate while MM method is also approximate.
\\\textbf{Remark 2} Our new explicit algorithm can be treated as an average estimation algorithm for solve all the overlapping group subproblem independently.
In Section 5, our numerical experiments show that Our new explicit algorithm
is accurate when $\beta$ is sufficiently large or sufficiently small, and is approximate to the other methods for instance the MM method for other $\beta$.

For the problem (\ref{norm21wogsp}), similar to (\ref{norm21ogspdo}), we can get
\begin{equation}\label{norm21wogspdo}
\begin{array}{r*{20}{l}}
  f_w(\mathbf{z})=&\min\limits_{\mathbf{z}} &\|\mathbf{z}\|_{w,2,1}+ \frac{\beta}{2}\|\mathbf{z} -\mathbf{x}\|_2^2\\
  =&\min\limits_{\mathbf{z}} &\sum_{i=1}^n \|w_g\circ(z_i)_g\|_2 + \frac{\beta}{2\sum_{k=1}^s(w_g)_k^2}\sum_{k=1}^s(w_g)_k^2\|\mathbf{z} -\mathbf{x}\|_2^2\\
  =&\min\limits_{\mathbf{z}} &\sum_{i=1}^n \|w_g\circ(z_i)_g\|_2 + \frac{\beta}{2\|w_g\|_2^2}\sum_{i=1}^n\|w_g\circ\left((z_i)_g -(x_i)_g\right)\|_2^2\\
  =&\min\limits_{\mathbf{z}} &\sum_{i=1}^n \|w_g\circ(z_i)_g\|_2 + \frac{\beta}{2\|w_g\|_2^2}\sum_{i=1}^n\|w_g\circ(z_i)_g -w_g\circ(x_i)_g\|_2^2\\
  =&\min\limits_{\mathbf{z}} &\sum_{i=1}^n \left(\|w_g\circ(z_i)_g\|_2 + \frac{\beta}{2\|w_g\|_2^2}\|w_g\circ(z_i)_g -w_g\circ(x_i)_g\|_2^2\right).\\
\end{array}
\end{equation}
See Figure~\ref{vector}(b) intuitively. All the symbols are the same as before. Similarly as before, we know that the necessary condition of minimizing the $i$th term of the last line in (\ref{norm21wogspdo}) is that
($w_g\circ(z_i)_g$) is parallel to ($w_g\circ(x_i)_g$). That is, $\frac{w_g\circ(z_i)_g}{\|w_g\circ(z_i)_g\|_2} =\frac{w_g\circ(x_i)_g}{\|w_g\circ(x_i)_g\|_2}$ for every $i$.
On the other hand, Let $W$=diag($w_g$) be a diagonal matrix with diagonal being the vector $w_g$, then $w_g\circ(x_i)_g=W(x_i)_g=W^T(x_i)_g$. If the vector $\mathbf{x'}$ (the same size as $(x_i)_g$) is parallel to the vector $\mathbf{z'}$, we have $\mathbf{x'}=\alpha\mathbf{z'}$.
Then, $W\mathbf{x'}=W\alpha\mathbf{z'}=\alpha W\mathbf{z'}$. We obtain that the vector $W\mathbf{x'}$ is also parallel to the vector $W\mathbf{z'}$.
Therefore, $\frac{W(z_i)_g}{\|W(z_i)_g\|_2} =\frac{W(x_i)_g}{\|W(x_i)_g\|_2}$, and each of them is a unit vector. Then, $\frac{W^T W(z_i)_g}{\|W(z_i)_g\|_2} =\frac{W^T W(x_i)_g}{\|W(x_i)_g\|_2}$.
That is, $$\frac{w_g\circ w_g\circ(z_i)_g}{\|w_g\circ(z_i)_g\|_2} =\frac{w_g\circ w_g\circ(x_i)_g}{\|w_g\circ(x_i)_g\|_2}.$$

Particularly, we first consider that $s=4$, $(z_i)_g=(z_{i-1},z_{i},z_{i+1},z_{i+2})$ and $w_g=(w_1,w_2,w_3,w_4)$. We mark $(w_g\circ(z_i)_g)'$ be the expansion of ($w_g\circ(z_i)_g$) similarly as $(z_i)_g'$, and we can get $\mathbf{x}=\frac{1}{\sum_{i=1}^4 w_i}\sum_{i=1}^{n} (w_g\circ(z_i)_g)'$. Then,
 the Euler equation of $f_w(\mathbf{z})$ is given by

\begin{equation}\label{subshrinkagenew2222w}
  \beta\left( \mathbf{z}- \mathbf{x}\right) + \frac{(w_g\circ w_g\circ(z_1)_g)'}{\|(w_g\circ(z_1)_g)'\|_2} + \cdots + \frac{(w_g\circ w_g\circ(z_n)_g)'}{\|(w_g\circ(z_n)_g)'\|_2}
  \ni \mathbf{0},
\end{equation}
\begin{equation}\label{subshrinkagenew2221w}
\frac{\beta}{\sum_{i=1}^4 w_i^2}\sum_{i=1}^{n}\left(w_g\circ(z_i)_g - w_g\circ(x_i)_g \right) + \frac{(w_g\circ w_g\circ(z_1)_g)'}{\|(w_g\circ(z_1)_g)'\|_2} + \cdots + \frac{(w_g\circ w_g\circ(z_n)_g)'}{\|(w_g\circ(z_n)_g)'\|_2}\ni \mathbf{0},
\end{equation}
\begin{equation}\label{subshrinkagenew22211w}
 \beta\left( \mathbf{z}- \mathbf{x}\right)+ \frac{(w_g\circ w_g\circ(x_1)_g)'}{\|(w_g\circ(x_1)_g)'\|_2} + \cdots + \frac{(w_g\circ w_g\circ(x_n)_g)'}{\|(w_g\circ(x_n)_g)'\|_2}\ni \mathbf{0},
\end{equation}

\begin{equation}\label{subshrinkagenew22211neww}
\mathbf{z} \ni \mathbf{x}-\frac{ 1}{\beta} \left(\frac{(w_g\circ w_g\circ(x_1)_g)'}{\|(w_g\circ(x_1)_g)'\|_2} + \cdots + \frac{(w_g\circ w_g\circ(x_n)_g)'}{\|(w_g\circ(x_n)_g)'\|_2}\right),
\end{equation}
for each component, we obtained
\begin{equation}\label{subshrinkagenew22211newiw}
z_i \ni x_i -\frac{ 1}{\beta} \left(\frac{w_4^2 x_i}{\|w_g\circ(x_{i-2})_g\|_2} +\frac{w_3^2 x_i}{\|w_g\circ(x_{i-1})_g\|_2} + \frac{w_2^2 x_i}{\|w_g\circ(x_i)_g\|_2}+ \frac{w_1^2 x_i}{\|w_g\circ(x_{i+1})_g\|_2}\right).
\end{equation}
Another expression is as follows.
\begin{equation}\label{subshrinkagenew22211newiw1}
z_i \ni \left(\frac{w_4^2 x_i}{\|w_g\|_2^2} -\frac{ 1}{\beta}\frac{w_4^2 x_i}{\|w_g\circ(x_{i-2})_g\|_2}\right) +\left(\frac{w_3^2 x_i}{\|w_g\|_2^2} -\frac{ 1}{\beta}\frac{w_3^2 x_i}{\|w_g\circ(x_{i-1})_g\|_2}\right) + \cdots+ \left(\frac{w_1^2 x_i}{\|w_g\|_2^2} -\frac{ 1}{\beta}\frac{w_1^2 x_i}{\|w_g\circ(x_{i+1})_g\|_2}\right).
\end{equation}

Similarly as Algorithm 1, we obtain the following algorithm for finding the minimizer of (\ref{norm21ogsp}).\\
\begin{tabular}{l}
\hline
\hline
{\textsc{\textbf{Algorithm 2}}} \textup{Direct shrinkage algorithm for the minimization problem (\ref{norm21ogsp})}\\
\hline
 \textit{\textbf{Input}}:  \\ \ \ Given vector $\mathbf{x}$, group size $s$, parameter $\beta$, weight vector $w_g=(w_1,w_2,\cdots,w_s)$.
\\
 \textit{\textbf{Compute}}:\\ Definition of $(x_i)_g$, for example $(x_i)_g = \left(x_{i-s_l},\cdots,x_i,\cdots,x_{i+s_r}\right)$, \\
 with $s_l =[\frac{s-1}{2}]$ and $s_r =[\frac{s}{2}]$. \\
If $s$ is even, then $w_g = [0, w_g]$ and $s=s+1$. \\
Let $w_g^r$ = fliplr($w_g$) be fliping $w_g$ over 180 degrees.\\
 Compute  $X_{n}=\left(\|(x_1)_g\|_2,\cdots,\|(x_i)_g\|_2,\cdots,\|(x_n)_g\|_2\right)$, by correlation of $w_g.^2$ (pointwise square)\\
 \quad \quad and $\mathbf{x}$, or by convolution of $w_g^r.^2$ and $\mathbf{x}$.\\
Compute $X'_{n}=\max\left( \frac{1}{\|w_g\|_2^2} - 1./(X_{n}\cdot\beta),0\right)$ pointwise. \\
Compute $F(x_i)$ by convolution of $w_g.^2$ and $X'_{n}$, or by correlation of $w_g^r.^2$ and $X'_{n}$.\\
\\
\hline
\end{tabular}\\

We can also see that Algorithm 2 only need 2 times convolution computations with time complexity $n*s$, which is just the same time complexity
as one step iteration in the MM method in \cite{CS2014}. Therefore, our method is much more efficient than MM method.

Here, thanks to the properties of inequalities, without loss of generality, let $\mathbf{x}\in[0,1]^n$, then we obtain that if
 \(\beta\leq \frac{\|w_g\|_2}{\sqrt{s}}\leq \frac{\|w_g\|_2}{\|x_g\|_2}=\frac{\|w_g\|_2^2}{\|w_g\|\|x_g\|_2}\leq \frac{\|w_g\|_2^2}{\|w_g\cdot x_g\|_2}\),
 then \(\beta\) is sufficiently small. However, we do not directly know what \(\beta\) is sufficiently large. In Section 5, after more than one thousand tests, we find that
 \(\beta\geq 30\cdot\frac{\|w_g\|_2}{\sqrt{s}}\) is sufficiently large generally.

For the problem (\ref{norm21wogspm}), similar to (\ref{norm21wogspdo}), we can obtain
\begin{equation}\label{norm21wogspdom}
\begin{array}{r*{20}{l}}
  f_W (A)=&\min\limits_{A}\|A\|_{W,2,1}+ \frac{\beta}{2}\|A -X\|_F^2\\
  =&\min\limits_{A}\sum_{i=1}^m\sum_{j=1}^n \|W_g\circ(A_{i,j})_g\|_F + \frac{\beta}{2\sum_{k_1=1}^{K_1}\sum_{k_2=1}^{K_2}(W_g)_{k_1,k_2}^2}\sum_{k_1=1}^{K_1}\sum_{k_2=1}^{K_2}(W_g)_{k_1,k_2}^2\|A -X\|_F^2\\
  =&\min\limits_{A}\sum_{i=1}^m\sum_{j=1}^n \|W_g\circ(A_{i,j})_g\|_F+ \frac{\beta}{2\|W_g\|_F^2}\sum_{i=1}^m\sum_{j=1}^n\|W_g\circ\left((A_{i,j})_g -(X_{i,j})_g\right)\|_F^2\\
  =&\min\limits_{A}\sum_{i=1}^m\sum_{j=1}^n \|W_g\circ(A_{i,j})_g\|_F+ \frac{\beta}{2\|W_g\|_F^2}\sum_{i=1}^m\sum_{j=1}^n\|W_g\circ(A_{i,j})_g -W_g\circ(X_{i,j})_g\|_F^2\\
  =&\min\limits_{A} \sum_{i=1}^m\sum_{j=1}^n \left(\|W_g\circ(A_{i,j})_g\|_F + \frac{\beta}{2\|W_g\|_F^2}\|W_g\circ(A_{i,j})_g -W_g\circ(X_{i,j})_g\|_F^2\right).\\
\end{array}
\end{equation}

 For example, we set $K_1=2,K_2=2$, and define $(A_{i,j})_g=(A_{i,j},A_{i,j+1};A_{i+1,j},A_{i+1,j+1})$ and $W_g=(W_{1,1},W_{1,2};W_{2,1},W_{2,2})$. Similar as the vector case, each $(A_{i,j})_g$ is a matrix with $\|(A_{i,j})_g\|_F=\sqrt{\sum_{k_1=1}^{K_1}\sum_{k_2=1}^{K_2}((A_{i,j})_g)_{k_1,k_2}^2}$. Notice that the Frobenius norm of a matrix is equal to the $\ell_2$ norm of a vector reshaped by the matrix. Then,
 the Euler equation of $f_w(\mathbf{z})$ is given by

\begin{equation}\label{subshrinkagenew2222wm}
  \beta\left( A- F\right) + \frac{(W_g\circ W_g\circ(A_{1,1})_g)'}{\|(W_g\circ(A_{1,1})_g)'\|_2} + \cdots + \frac{(W_g\circ W_g\circ(A_{n,n})_g)'}{\|(W_g\circ(A_{n,n})_g)'\|_2}
  \ni \mathbf{0},
\end{equation}
where $((A_{i,j})_g)'$ is defined similarly as $((z_i)_g)'$, which is an expansion of $(A_{i,j})_g$. These symbols remain consistent as default through this paper.

\begin{equation}\label{subshrinkagenew2221wm}
\beta\left( A- X\right) + \frac{(W_g\circ W_g\circ(X_{1,1})_g)'}{\|(W_g\circ(X_{1,1})_g)'\|_2} + \cdots + \frac{(W_g\circ W_g\circ(X_{n,n})_g)'}{\|(W_g\circ(X_{n,n})_g)'\|_2}\ni \mathbf{0},
\end{equation}

\begin{equation}\label{subshrinkagenew22211newwm}
A \ni X-\frac{ 1}{\beta} \left(\frac{(W_g\circ W_g\circ(X_{1,1})_g)'}{\|(W_g\circ(X_{1,1})_g)'\|_2} + \cdots + \frac{(W_g\circ W_g\circ(X_{n,n})_g)'}{\|(W_g\circ(X_{n,n})_g)'\|_2}\right),
\end{equation}
for each component, we obtained
\begin{equation}\label{subshrinkagenew22211newiwm}
A_{i,j} \ni X_{i,j} -\frac{ 1}{\beta} \left(\frac{W_{2,2}^2 X_{i,j}}{\|W_g\circ(X_{i-1,j-1})_g\|_2} +\frac{W_{2,1}^2 X_{i,j}}{\|W_g\circ(X_{i-1,j})_g\|_2} + \frac{W_{1,2}^2 X_{i,j}}{\|W_g\circ(X_{i,j-1})_g\|_2}+ \frac{W_{1,1}^2 X_{i,j}}{\|W_g\circ(X_{i,j})_g\|_2}\right).
\end{equation}

Therefore, we can obtain a similar algorithm on the former formula (\ref{subshrinkagenew22211newiwm}) for finding the minimizer of (\ref{norm21ogsp}) as Algorithm 2.

\section{Several extensions}
\subsection{Other boundary conditions}
In Section 2, we gave the explicit shrinkage formulas for one class of OGS problems (\ref{norm21ogsp}), (\ref{norm21wogsp}) and (\ref{norm21wogspm}). In order to achieve a simply deduction, we assume that period boundary conditions are used. One may confuse that whether period boundary conditions are always good for regularization problems such as signal processing or image processing, since natural signals or images are often asymmetric. However, in these problems, assuming a kind of boundary condition is necessary for simplifying the problem and making the computation possible \cite{HN2006}. There are kinds of boundary conditions, such as zero boundary conditions and reflective boundary conditions. Period boundary conditions are often used in optimization because it can be computed fast as before or other fast computation for example, computation of matrix of block circulant with circulant blocks by fast Fourier transforms \cite{HN2006,WY2008,YYZW2009}.
 \begin{figure}
  \centering
\subfigure[zero boundary conditions]{\includegraphics[width=0.48\textwidth,clip]{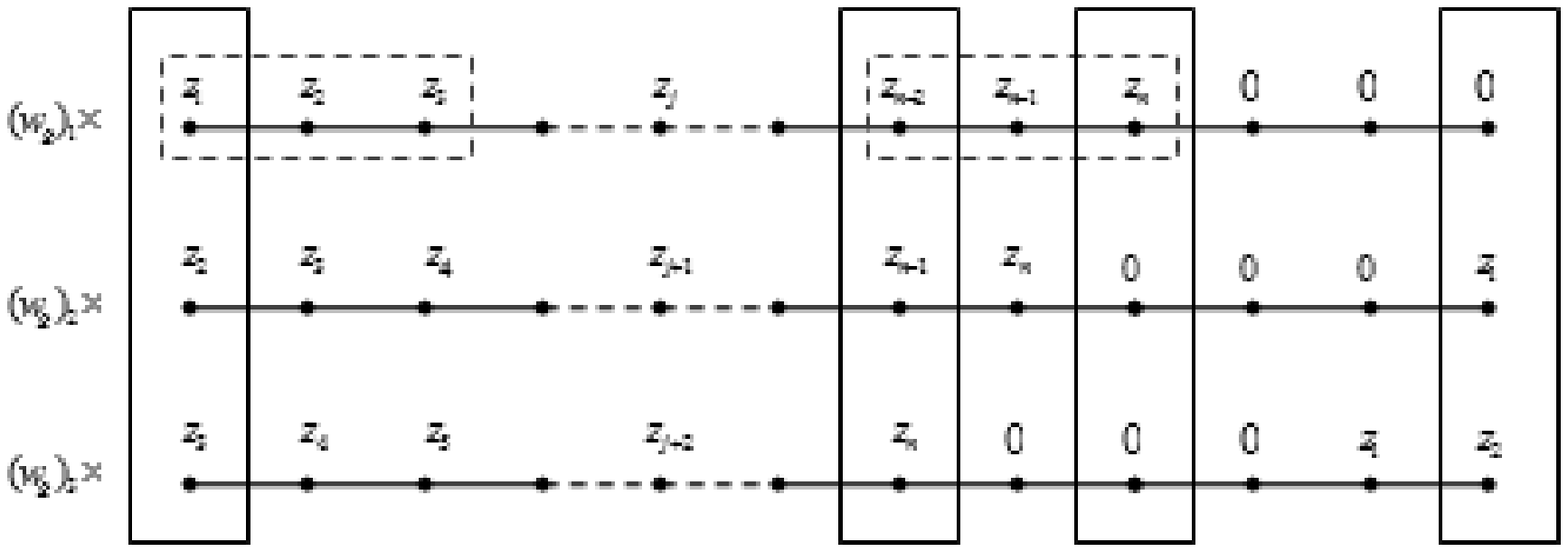}}
  \subfigure[reflective boundary conditions]{\includegraphics[width=0.48\textwidth,clip]{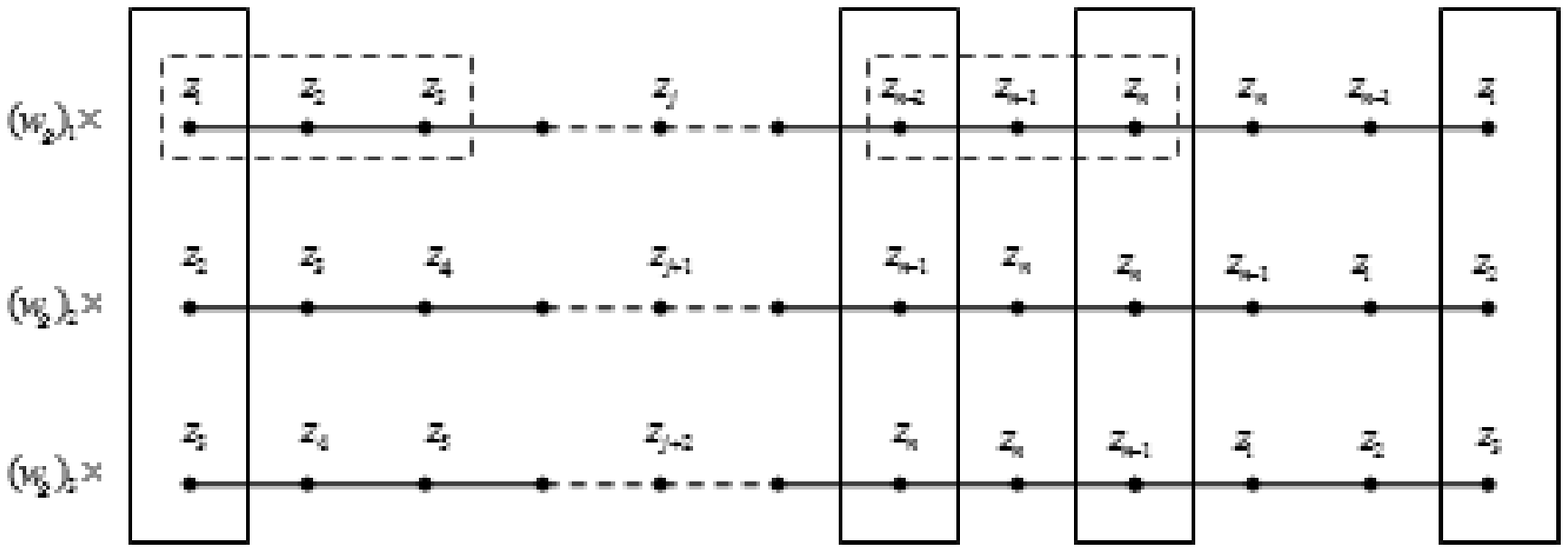}}
  \caption{Vector case with other boundary conditions.}
  \label{vectoro}
 \end{figure}

In this section, we consider other boundary conditions such as zero boundary conditions and reflective boundary conditions. For simplification, we only consider the vector case, while the results can be easily expanded to the matrix case similar as Section 2. When zero boundary conditions are used, we can expand the original vectors (or signal) $\mathbf{z}$ ($=(z_i)_{i=1}^{n}$) and $\mathbf{x}$ by two $s$-length vectors on the both hands of the original vectors, which are $\tilde{\mathbf{z}}$ ($=\left[\tilde{z}_{-s},\cdots,\tilde{z}_{-1},(\tilde{z}_i)_{i=1}^{n},\tilde{z}_{n+1},\cdots,\tilde{z}_{n+s}\right]$ $=\left[\mathbf{0}_s,(z_i)_{i=1}^{n},\mathbf{0}_s\right]$) and $\tilde{\mathbf{x}}$ respectively. Then, the results and algorithms in Section 2 are similar as assuming period boundary conditions on $\tilde{\mathbf{z}}$ and $\tilde{\mathbf{x}}$. See Figure~\ref{vectoro}(a) intuitively. Moreover, according to the definition of weighted generalized $\ell_{2,1}$-norm, although zero boundary conditions seems better than than period boundary conditions, our numerical results will show that the results form different boundary conditions are almost the same in practice (see Section 5.1). Therefore, we will choose the zero boundary conditions to solve the problems (\ref{norm21ogsp}), (\ref{norm21wogsp}) and (\ref{norm21wogspm}) in the following sections.

When reflective boundary conditions are assumed, the results are also the same. We only need to extend the original vector $\mathbf{x}$ to $\hat{\mathbf{x}}$ with reflective boundary conditions. See Figure~\ref{vectoro}(b) intuitively.

\subsection{Nonpositive weights and different weights in groups}
In this section, we show that the weights vector $w_g$ and matrix $W_g$ in the former sections can contain arbitrary $s$ entries with arbitrary real numbers.
On one hand, the zero value can be the arbitrary entries of the weights vector or matrix. For example, the original $\ell_1$ regularization problem (\ref{norm11}) can be seen as a special form of weighted generalized $\ell_{2,1}$-norm regularization problems (\ref{norm21wogsp}) with $s=1$, or $s=3$ and $w_g=(1,0,0)$ for the example of Figure~\ref{vector}(b).

On the other hand, the norm is with the property of positive homogeneity, $\|(-1)\cdot w_g\|_{any} = \|w_g\|_{any}$, where ``$\|\cdot\|_{any}$'' can be arbitrary norm such as $\|\cdot\|_{1}$ and $\|\cdot\|_{2}$. Therefore, for the regularization problems (\ref{norm21wogsp}) (or (\ref{norm21wogspm})), the weights vector $w_g$ (or matrix $W_g$) is the same as $|w_g|$ (or $|W_g|$), where $|\cdot|$ the point-wise absolute value. Therefore, in general, our results are true, whatever the number of entries included in the weights vector $w_g$ and matrix $W_g$, and whatever the real number value of these entries.

However, if the weight $w_g$ is dependent on the group index $i$, that is for example, $w_{i,g}\circ (z_i)_g = \left((w_{i,g})_1 z_i, (w_{i,g})_2 z_{i+1},(w_{i,g})_3 z_{i+2}\right)$, we cannot solve the relevant problems easily since our method fails. As we mentioned before, we only focus on that the weight $w_g$ is independent on the group index $i$, that is $w_{i,g}=w_g$ for all $i$, which means that it is with translation invariant overlapping groups.

\section{Applications in TV regularization problems with OGS}

The TV regularizer was firstly introduced by Rudin et.al\cite{ROF1992}(ROF). It is used in many fields, for instance, denoising and deblurring problems. Several fast algorithms have
 been proposed, such as Chambelle \cite{Chamb2004} and fast TV deconvolution (FTVd) \cite{WY2008,YYZW2009}. Its corresponding minimization task is:
\begin{equation}\label{DROF}
  \min_f \|f \|_\text{TV} + \frac{\mu}{2} \|Hf-g\|_2^2 , \ \
\end{equation}
where $ \|f \|_\text{TV}:= \sum\limits_{1\leqslant i,j \leqslant n} \|(\nabla f)_{i,j}\|_2 = \sum\limits_{1\leqslant i,j \leqslant n}$ $ \sqrt{{|(\nabla_x f)}_{i,j}|^2 +
 |(\nabla_y f)_{i,j}|^2 }$ (called isotropic TV, ITV), or $ \|f \|_{TV}:= \sum\limits_{1\leqslant i,j \leqslant n}$    $\|(\nabla f)_{i,j}\|_1$ $= \sum\limits_{1\leqslant i,j \leqslant n} |(\nabla_x f)_{i,j}| + |(\nabla_y f)_{i,j}|$ (called anisotropic TV, ATV), $H$ denotes the blur matrix, and $g$ denotes the given observed image with blur and noise.
Operator $\nabla: \mathbb{R}^{n^2}\rightarrow\mathbb{R}^{2\times{n^2}}$ denotes the discrete gradient
 operator (under periodic boundary conditions) which is defined by
 $(\nabla f)_{i,j}=((\nabla_x f)_{i,j},(\nabla_y f)_{i,j} )$,
with
\begin{equation*}
(\nabla_x f)_{i,j} =\left\{\begin{array}{lll}  f_{i+1,j}-f_{i,j}&{\rm if}&i<n,\\
                                            f_{1,j}-f_{n,j}&{\rm if}&i=n,\end{array}\right.
(\nabla_y f)_{i,j} =\left\{\begin{array}{lll}  f_{i,j+1}-f_{i,j}&{\rm if}&j<n,\\
                                            f_{i,1}-f_{i,n}&{\rm if}&j=n,\end{array}\right.
\end{equation*}
for $i,j=1,2,\cdots,n$, where $f_{i,j}$ refers to the $((j-1)n+i)$th entry
of the vector $f$ (it is the $(i,j)$th pixel location of the
$n\times n$ image, and this notation remains valid throughout the paper
unless otherwise specified). Notice that $H$ is a matrix of block circulant with circulant blocks (BCCB) structure when periodic boundary conditions are applied or other structures when other boundary conditions are applied \cite{HN2006}.

Recently, Selesnick et. al \cite{SC2013} proposed an OGS TV regularizer to one-dimensional signal denoising. They applied the MM method to solve their model. Their numerical experiments showed that their method can overcome staircase effects effectively and get better results. However, their method has the disadvantages of the low speed of computation and the difficulty to be extended to the two-dimensional image case because they did not choose a variable substitution method. More recently, Liu et. al \cite{LHS2013} proposed an OGS TV regularizer for two-dimensional image denoising and deblurring under Gaussian noise, and Liu et. al \cite{LHLL2013} proposed an OGS TV regularizer for image deblurring under impulse noise. Both of them used a variable substitution method and the ADMM framelet with an inner MM iteration for solving the subproblems similar as (\ref{norm21ogsp}). Therefore, they did not get convergence theorem results because of the inner iterations. However, from our results in Section 2, we can solve the subproblems exactly, and we can get convergence theorem results under the ADMM framelet. Moreover, when the MM method is used in the inner iterations in \cite{LHS2013,LHLL2013}, they can only solve the ATV case but not the ATV case with OGS while our methods can solve both ATV and ITV cases.

Firstly, we defined the ATV case with OGS under Gaussian noise and impulse noise respectively which is similar like \cite{LHS2013,LHLL2013}. For the Gaussian noise case,
\begin{equation}\label{L2OGSTV}
  \min_f \|(\nabla_x f)\|_{W,2,1} +\|(\nabla_y f)\|_{W,2,1} + \frac{\mu}{2} \|Hf-g\|_2^2.
\end{equation}
For the impulse case,
\begin{equation}\label{L1OGSTV}
  \min_f \|(\nabla_x f)\|_{W,2,1} +\|(\nabla_y f)\|_{W,2,1} + \mu \|Hf-g\|_1.
\end{equation}
We call the former model as TV OGS $L_2$ model and the latter model as TV OGS $L_1$ model respectively.

Then, we defined the ITV case with OGS. For Gaussian noise and impulse noise, we only change the former two terms of (\ref{L2OGSTV}) and (\ref{L1OGSTV}) respectively by $\|A\|_{W,2,1}$.
Here, $A$ is a high-dimensional matrix with each entry $A_{i,j} = ((\nabla_x f)_{i,j};(\nabla_y f)_{i,j})$. \\
\textbf{Remark 3}. Here and the following sections, $A$ can be treated as $((\nabla_x f);(\nabla_y f))$ for simplicity in vector and matrix computation, although the computation of $\|A\|_{W,2,1}$ is should be tread as pointwise with $A_{i,j} = ((\nabla_x f)_{i,j};(\nabla_y f)_{i,j})$ since the computation are almost all pointwise.

Moreover, we also consider constrained model as listing in \cite{LHS2013,LHLL2013,CTY2013}. For any
true digital image, its pixel value can attain only a finite number
of values. Hence, it is natural to require all pixel values of
the restored image to lie in a certain interval $[a, b]$,
see \cite{CTY2013} for more details. In general, with the easy computation and the certified results in \cite{CTY2013}, we only consider all the images located on the standard range $[0, 1]$. Therefore,we define a projection operator $\mathcal{P}_{\Omega}$ on the set $\Omega =\left\{f\in \mathbb{R}^{n\times n}|0\leqslant f\leqslant 1\right\}$,
\begin{equation}\label{Cpro}
  \mathcal{P}_{\Omega}(f)_{i,j}=
  \left\{\begin{array}{ll}
  0,&f_{i,j}<0,\\
  f_{i,j},&f_{i,j}\in[0,1],\\
  1,&f_{i,j}>1.
  \end{array}\right.
\end{equation}
\subsection{Constrained TV OGS $L_2$ model}
For the constrained model (the ATV case) (called CATVOGSL2), we have
\begin{equation}\label{L2COGSTV}
  \min_{u\in\Omega,v_x,v_y,f} \left\{\|(v_x)\|_{W,2,1} +\|(v_y)\|_{W,2,1} + \frac{\mu}{2} \|Hf-g\|_2^2:\ v_x=\nabla_x f, v_y=\nabla_y f,u=f\right\}.
\end{equation}
The augmented Lagrangian
function of (\ref{L2COGSTV}) is
\begin{equation}\label{L2COGSTVAL}
  \begin{array}{rl}
   \mathcal{L}(v_x,v_y,u,f;\lambda_1,\lambda_2,\lambda_3)=&\|v_x\|_{W,2,1} - \lambda_1^T (v_x-\nabla_x f)
    +\frac{\beta_1}{2}\|v_x-\nabla_x f\|_2^2\\
    &+\|v_y\|_{W,2,1} - \lambda_2^T (v_y-\nabla_y f)
    +\frac{\beta_1}{2}\|v_y-\nabla_y f\|_2^2\\
    &-\lambda_3^T (u -f) + \frac{\beta_2}{2}\|u-f\|_2^2 + \frac{\mu}{2} \|Hf-g\|_2^2,
       \end{array}
\end{equation}
where $\beta_1, \beta_2>0$ are penalty parameters and $\lambda_1, \lambda_2, \lambda_3 \in\mathbb{R}^{n^2}$ are the Lagrange multipliers.
The solutions are according to the scheme of ADMM in Gabay \cite{GM1976}, and we refer to some applications in image
processing which can be solved by ADMM, e.g., \cite{EB1992,RG1984,HY1998,Esser2009,NWY2010,ZBBO2010,
YZY2010,ZBO2010}. For a given $(v_x^k,v_y^k,u^k,f^k$; $\lambda_1^k,\lambda_2^k,\lambda_3^k)$, the next iteration
$(v_x^{k+1},v_y^{k+1},u^{k+1}$, $f^{k+1}$; $\lambda_1^{k+1},\lambda_2^{k+1},\lambda_3^{k+1})$ is generated as follows:

1. Fix $f = f^k$, $\lambda_1=\lambda_1^k, \lambda_2=\lambda_2^k, \lambda_3=\lambda_3^k$,
and minimize (\ref{L2COGSTVAL}) with respect to $v_x$, $v_y$ and $u$.
Respect to $v_x$ and $v_y$,
\begin{equation}\label{COTVSUB1}\begin{array}{rl}
  v_x^{k+1}&=\arg\min \|v_x\|_{W,2,1} - {\lambda_1^k}^T (v_x-\nabla_x f^k)
    +\frac{\beta_1}{2}\|v_x-\nabla_x f^k\|_2^2\\
    &=\arg\min \|v_x\|_{W,2,1} +\frac{\beta_1}{2}\|v_x-\nabla_x f^k - \frac{\lambda_1^k}{\beta_1}\|_2^2,\\
  \end{array}
\end{equation}
\begin{equation}\label{COTVSUB11}\begin{array}{rl}
  v_y^{k+1}&=\arg\min \|v_y\|_{W,2,1} - {\lambda_2^k}^T (v_y-\nabla_y f^k)
    +\frac{\beta_1}{2}\|v_y-\nabla_y f^k\|_2^2\\
    &=\arg\min \|v_y\|_{W,2,1} +\frac{\beta_1}{2}\|v_y-\nabla_y f^k - \frac{\lambda_2^k}{\beta_1}\|_2^2.\\
  \end{array}
\end{equation}
It is obvious that problems (\ref{COTVSUB1}) and (\ref{COTVSUB11}) match
the framework of the problem (\ref{norm21wogspm}), thus the minimizers of (\ref{COTVSUB1}) and (\ref{COTVSUB11})
can be obtained by using the formulas in Section 2.2.

Respect to $u$,
\begin{equation*}\begin{array}{rl}
  u^{k+1}&=\arg\min -{\lambda_3^k}^T (u -f^k) + \frac{\beta_2}{2}\|u-f^k\|_2^2\\
    &=\arg\min \frac{\beta_2}{2}\|u-f^k - \frac{\lambda_3^k}{\beta_2}\|_2^2.\\
  \end{array}
\end{equation*}
The minimizer is given explicitly by
\begin{equation}\label{COTVSUB3}
  u^{k+1}=\mathcal{P}_{\Omega}\left[f^k + \frac{\lambda_3^k}{\beta_2}\right].
\end{equation}

2.  Compute $f^{k+1}$ by solving the normal equation
\begin{equation}\label{COTVSUB4}
 \begin{array}{l}
 (\beta_1 (\nabla_x^* \nabla_x+\nabla_y^* \nabla_y)+\mu H^* H + \beta_2 I)f^{k+1} \\
=\nabla_x^*(\beta_1 v_x^{k+1}-\lambda^k_1)+{\nabla_y}^*(\beta_1 v_y^{k+1}-\lambda_2^k)+\mu H^* g + \beta_2 (u^{k+1} - \frac{\lambda_3^k}{\beta_2}),\\
\end{array}
\end{equation}
where ``$*$'' denotes the conjugate transpose, see \cite{CX2010} for more details.
Since all the parameters are positive, the coefficient matrix in
(\ref{COTVSUB4}) are always invertible and symmetric positive-definite. In addition, note that $H$, $\nabla_x$, $\nabla_y$ have BCCB
structure under periodic boundary conditions. We
know that the computations with BCCB matrix can be very
efficient by using fast Fourier transforms.

3.  Update the multipliers via
\begin{equation}\label{COTVSUB5}
  \left\{
  \begin{array}{*{20}{l}}
  \lambda_1^{k+1}&=&\lambda_1 ^k - \gamma\beta_1(v_x^{k+1}-\nabla_x f^{k+1}),\\
  \lambda_2^{k+1}&=&\lambda_2 ^k - \gamma\beta_1(v_x^{k+1}-\nabla_x f^{k+1}),\\
  \lambda_3^{k+1}&=&\lambda_3 ^k - \gamma\beta_2(u^{k+1} -f^{k+1}).\\
  \end{array} \right.
\end{equation}

Based on the discussions above, we present the ADMM algorithm for solving the convex CATVOGSL2 model (\ref{L2COGSTV}), which is
shown as Algorithm 3.\\\\
\begin{tabular}{l}
\hline
\hline
{\textsc{\textbf{Algorithm 3}}} \textup{CATVOGSL2 for the minimization problem (\ref{L2COGSTV})}\\
\hline
 \textit{\textbf{initialization}}:  \\ \ \ Starting point $f^0=g$, $k=0$, $\beta_1$, $\beta_2$, $\gamma$, $\mu$,
 group size $K_1\times K_2$,\\
 \ \ weighted matrix $W_g$, $\lambda_i^0=0$, $i=1,2,3$.
\\
 \textit{\textbf{iteration}}:\\ 
$\begin{array}{l}
 1.\ Compute\ v_x^{k+1}\ and\ v_y^{k+1}\  according\ to\ (\ref{COTVSUB1})\ and\ (\ref{COTVSUB11}), \\
\quad and\  compute\ u^{k+1}\ according\ to\ (\ref{COTVSUB3}).\\
\end{array}$\\ 
$\begin{array}{l}
2.\ Compute\ f^{k+1}\ by\ solving\ (\ref{COTVSUB4}).\\
3.\ update\ \lambda_i^0=0, i=1,2,3\ according\ to\ (\ref{COTVSUB5}).\quad\quad\quad\quad\quad\ ~\\
4.\ k=k+1.\\
\end{array}$\\
\textit{\textbf{until a stopping criterion is satisfied.}}\\
\hline
\end{tabular}\\

Algorithm 3 is an application of ADMM for the case with two blocks of variables $(v_x, v_y, u)$ and $f$.  Thus, its convergence is guaranteed by the theory of ADMM \cite{GM1976,EB1992, HY1998}, and we
summarize it in the following theorem.\\
\textbf{Theorem 1}.
1). When $K_1=K_2=1$, our method CATVOGSL2 degenerates to the constrained ATV $L_2$ model in Chan. \cite{CTY2013}. Therefore, we have for $\beta_1,\beta_2>0$ and $\gamma\in(0,\frac{1+\sqrt{5}}{2})$, the sequence $(v_x^k,v_y^k,u^k,f^k$; $\lambda_1^k,\lambda_2^k,\lambda_3^k)$
generated by Algorithm 1 from any initial point $(f^0$; $\lambda_1^0,\lambda_2^0,\lambda_3^0)$ converges to $(v_x^\diamond,v_y^\diamond,u^\diamond,f^\diamond$; $\lambda_1^\diamond,\lambda_2^\diamond,\lambda_3^\diamond)$, where $(v_x^\diamond,v_y^\diamond,u^\diamond,f^\diamond)$ is a solution of (\ref{L2COGSTV}).\\
2). When $K_1\neq 1$ or $K_2\neq 1$, form Section 2.2, when $\beta_1$ is sufficiently large, the computation in every step of Algorithm 3 is accurate. Therefore, we have for $\beta_2>0$, $\beta_1>L>0$ (L is a sufficiently large real number) and $\gamma\in(0,\frac{1+\sqrt{5}}{2})$, the sequence $(v_x^k,v_y^k,u^k,f^k$; $\lambda_1^k,\lambda_2^k,\lambda_3^k)$
generated by Algorithm 1 from any initial point $(f^0$; $\lambda_1^0,\lambda_2^0,\lambda_3^0)$ converges to $(v_x^\diamond,v_y^\diamond,u^\diamond,f^\diamond$; $\lambda_1^\diamond,\lambda_2^\diamond,\lambda_3^\diamond)$, where $(v_x^\diamond,v_y^\diamond,u^\diamond,f^\diamond)$ is a solution of (\ref{L2COGSTV}).\\

\textbf{Remark 4}. Here, for the image $f$, we use period boundary conditions because of the fast computation. However, for $v_x$, $v_y$, we use zero boundary conditions. Because $v_x$ and $v_y$ substitute the gradient of the image, zero boundary conditions seems better for the definition of the generalized norm $\ell_{2,1}$ on $v_x$ and $v_y$ as mentioned in Section 3.1. Therefore, the boundary conditions of the image and its gradient are different and independent. These remain valid throughout the paper unless otherwise specified.

For the ITV case, the constrained model (CITVOGSL2) is
\begin{equation}\label{L2COGSITV}
  \min_{u\in\Omega,A,f} \left\{\|A\|_{W,2,1} + \frac{\mu}{2} \|Hf-g\|_2^2:\ A= (\nabla_x f;\nabla_y f),u=f\right\}.
\end{equation}
We can also get convergence theorem similar as Theorem 1, whose detail will be presented in the next section for the constrained TV OGS $L_1$ model. We call this relevant algorithm CITVOGSL2.

\subsection{Constrained TV OGS $L_1$ model}

For the constrained model (the ITV case) (called CITVOGSL1), we have
\begin{equation}\label{L1COGSTV}
  \min_{u\in\Omega,A,f} \left\{\|A\|_{W,2,1} + {\mu} \|z\|_1:\ z=Hf-g, A=(\nabla_x f;\nabla_y f),u=f \right\}.
\end{equation}
The augmented Lagrangian
function of (\ref{L1COGSTV}) is
\begin{equation}\label{L1COGSTVAL}
  \begin{array}{rl}
   \mathcal{L}(v_x,v_y,z,w,f;\lambda_1,\lambda_2,\lambda_3,\lambda_4)=&\|A\|_{W,2,1} - (\lambda_1^T,\lambda_2^T) (A-\left[\begin{matrix}\nabla_x f\\ \nabla_y f\\ \end{matrix}\right])
    +\frac{\beta_1}{2}\left\|A-\left[\begin{matrix}\nabla_x f\\ \nabla_y f\\ \end{matrix}\right]\right\|_2^2\\
    &+\mu \|z\|_1 - \lambda_3^T \left(z-(Hf - g)\right) + \frac{\beta_2}{2}\|z-(Hf-g)\|_2^2\\
    &-\lambda_4^T (w -f) + \frac{\beta_3}{2}\|w-f\|_2^2,
       \end{array}
\end{equation}
where $\beta_1, \beta_2, \beta_3>0$ are penalty parameters and $\lambda_1, \lambda_2, \lambda_3, \lambda_4 \in\mathbb{R}^{n^2}$ are the Lagrange multipliers.
For a given $(A^k,z^k,u^k,f^k$; $\lambda_1^k,\lambda_2^k,\lambda_3^k,\lambda_4^k)$, the next iteration
$(A^{k+1},z^{k+1},u^{k+1}$, $f^{k+1}$; $\lambda_1^{k+1},\lambda_2^{k+1},\lambda_3^{k+1},$ $\lambda_4^{k+1})$ is generated as follows:

1. Fix $f = f^k$, $\lambda_1=\lambda_1^k, \lambda_2=\lambda_2^k, \lambda_3=\lambda_3^k, \lambda_4=\lambda_4^k$,
and minimize (\ref{L1COGSTVAL}) with respect to $A$, $z$ and $u$. Respect to $A$,
\begin{equation}\label{COTVSUB21}\begin{array}{rl}
  A^{k+1}&=\arg\min \|A\|_{W,2,1} -  (\lambda_1^T,\lambda_2^T) (A-\left[\begin{matrix}\nabla_x f^k\\ \nabla_y f^k\\ \end{matrix}\right])
    +\frac{\beta_1}{2}
  \left\|A-\left[\begin{matrix}\nabla_x f^k\\ \nabla_y f^k\\ \end{matrix}\right]\right\|_2^2,\\
\left[\begin{matrix}A_1^{k+1}\\ A_2^{k+1}\\ \end{matrix}\right]&=\arg\min \|A\|_{W,2,1} +\frac{\beta_1}{2}\left\|\left[\begin{matrix}A_1\\ A_2\\ \end{matrix}\right]-\left[\begin{matrix}\nabla_x f^k\\ \nabla_y f^k\\ \end{matrix}\right]-
    \left[\begin{matrix}{\lambda_1^k}/{\beta_1}\\ {\lambda_2^k}/{\beta_1}\\ \end{matrix}\right]\right\|_2^2.\\
  \end{array}
\end{equation}
It is obvious that problems (\ref{COTVSUB21}) match
the framework of the problem (\ref{norm21wogspm}), thus the solution of (\ref{COTVSUB21})
can be obtained by using the formulas in Section 2.2.

 Respect to $z$,
\begin{equation*}\begin{array}{rl}
  z^{k+1}&=\arg\min\mu \|z\|_1 - {\lambda_3^k}^T \left(z-(Hf^k - g)\right) + \frac{\beta_2}{2}\|z-(Hf^k-g)\|_2^2\\
    &=\arg\min\mu \|z\|_1 + \frac{\beta_2}{2}\|z-(Hf^k-g) - \frac{\lambda_3^k}{\beta_2}\|_2^2.\\
  \end{array}
\end{equation*}

The minimization with respect to $z$ can be given by (\ref{shrinkage1}) and (\ref{shrinkage1g}) explicitly, that is,
\begin{equation}\label{COTVSUB211}
  z^{k+1}={\rm sgn}\left\{Hf^k - g +\frac{\lambda_3^k}{\beta_2}\right\}\circ \max\left\{|Hf^k - g +\frac{\lambda_3^k}{\beta_2}|-\frac{\mu}{\beta_2},0\right\}.
\end{equation}

Respect to $u$,
\begin{equation*}\begin{array}{rl}
  u^{k+1}&=\arg\min -{\lambda_4^k}^T (u -f^k) + \frac{\beta_3}{2}\|u-f^k\|_2^2\\
    &=\arg\min \frac{\beta_3}{2}\|u-f^k - \frac{\lambda_4^k}{\beta_3}\|_2^2.\\
\end{array}
\end{equation*}
The minimizer is given explicitly by
\begin{equation}\label{COTVSUB31}
  u^{k+1}=\mathcal{P}_{\Omega}\left[f^k + \frac{\lambda_4^k}{\beta_3}\right].
\end{equation}

2.  Compute $f^{k+1}$ by solving the following normal equation similarly as the last section.
\begin{equation}\label{COTVSUB41}
 \begin{array}{l}
(\beta_1 (\nabla_x^* \nabla_x+\nabla_y^* \nabla_y)+\beta_2 H^* H + \beta_3 I)f^{k+1} \\
=\nabla_x^*(\beta_1 A_1^{k+1}-\lambda^k_1)+{\nabla_y}^*(\beta_1 A_2^{k+1}-\lambda_2^k)+H^* (\beta_2 z^{k+1} -\lambda_3^k) +\beta_2 H^* g + \beta_3 (u^{k+1} - \frac{\lambda_4^k}{\beta_3}).\\
\end{array}
\end{equation}

3.  Update the multipliers via
\begin{equation}\label{COTVSUB51}
  \left\{
  \begin{array}{*{20}{l}}
\left[\begin{matrix}\lambda_1^{k+1}\\ \lambda_2^{k+1}\\ \end{matrix}\right]
  &=&\left[\begin{matrix}\lambda_1^{k}\\ \lambda_2^{k}\\ \end{matrix}\right]- \gamma\beta_1 \left(\left[\begin{matrix}A_1^{k+1}\\ A_2^{k+1}\\ \end{matrix}\right]-
  \left[\begin{matrix}\nabla_x f^{k+1}\\ \nabla_y f^{k+1}\\ \end{matrix}\right]\right),\\
  \lambda_3^{k+1}&=&\lambda_3 ^k - \gamma\beta_2(z^{k+1}-(Hf^{k+1}-g)),\\
  \lambda_4^{k+1}&=&\lambda_4 ^k - \gamma\beta_3(u^{k+1} -f^{k+1}).\\
  \end{array} \right.
\end{equation}

Based on the discussions above, we present the ADMM algorithm for solving the convex CITVOGSL1 model (\ref{L1COGSTV}), which is
shown as Algorithm 4.\\\\
\begin{tabular}{l}
\hline
\hline
{\textsc{\textbf{Algorithm 4}}} \textup{CITVOGSL1 for the minimization problem (\ref{L1COGSTV})}\\
\hline
 \textit{\textbf{initialization}}:  \\ \ \ Starting point $f^0=g$, $k=0$, $\beta_1$, $\beta_2$, $\beta_3$, $\gamma$, $\mu$,
 group size $K_1\times K_2$,\\
  \ \ weighted matrix $W_g$, $\lambda_i^0=0$, $i=1,2,3,4$.\\
 \textit{\textbf{iteration}}:\\
$\begin{array}{l}
 1.\ Compute\ A^{k+1}\  according\ to\ (\ref{COTVSUB21}), \\
 \quad and\  compute\ z^{k+1}\ according\ to\ (\ref{COTVSUB211}).\\
\quad and\  compute\ u^{k+1}\ according\ to\ (\ref{COTVSUB31}).\\
\end{array}$\\
$\begin{array}{l}
2.\ Compute\ f^{k+1}\ by\ solving\ (\ref{COTVSUB41}).\\
3.\ update\ \lambda_i^0=0, i=1,2,3,4\ according\ to\ (\ref{COTVSUB51}).\quad\quad\quad\quad\quad\quad~\\
4.\ k=k+1.\\
\end{array}$\\
\textit{\textbf{until a stopping criterion is satisfied.}}\\
\hline
\end{tabular}\\

Algorithm 2 is an application of ADMM for the case with two blocks of variables $(A, z, u)$ and $f$.  Thus, its convergence is guaranteed by the theory of ADMM \cite{GM1976,EB1992, HY1998}, and we
summarize it in the following theorem.\\
\textbf{Theorem 2}.
1). When $K_1=K_2=1$, our method CITVOGSL1 degenerates to the constrained ITV $L_1$ model similarly as Chan. \cite{CTY2013}. Therefore, we have for $\beta_1,\beta_2,\beta_3>0$ and $\gamma\in(0,\frac{1+\sqrt{5}}{2})$, the sequence $(A^k,z^k,u^k,f^k$; $\lambda_1^k,\lambda_2^k,\lambda_3^k,\lambda_4^k)$
generated by Algorithm 2 from any initial point $(f^0$; $\lambda_1^0,\lambda_2^0,\lambda_3^0,\lambda_4^0)$ converges to $(A^\diamond,z^\diamond,u^\diamond,f^\diamond$; $\lambda_1^\diamond,\lambda_2^\diamond,\lambda_3^\diamond,\lambda_4^\diamond)$, where $(A^\diamond,z^\diamond,u^\diamond,f^\diamond)$ is a solution of (\ref{L1COGSTV}).\\
2). When $K_1\neq 1$ or $K_2\neq 1$, form Section 2.2, when $\beta_1$ is sufficiently large, the the computation in every step of Algorithm 4 is accurate. Therefore, we have for $\beta_2,\beta_3>0$, $\beta_1>L>0$ (L is a sufficiently large real number) and $\gamma\in(0,\frac{1+\sqrt{5}}{2})$, the sequence $(A^k,z^k,u^k,f^k$; $\lambda_1^k,\lambda_2^k,\lambda_3^k,\lambda_4^k)$
generated by Algorithm 2 from any initial point $(f^0$; $\lambda_1^0,\lambda_2^0,\lambda_3^0,\lambda_4^0)$ converges to $(A^\diamond,z^\diamond,u^\diamond,f^\diamond$; $\lambda_1^\diamond,\lambda_2^\diamond,\lambda_3^\diamond,\lambda_4^\diamond)$, where $(A^\diamond,z^\diamond,u^\diamond,f^\diamond)$ is a solution of (\ref{L1COGSTV}).

For the ATV case, the constrained model (CATVOGSL1) is
\begin{equation}\label{L2COGSITV}
  \min_{u\in\Omega,v_x,v_y,f} \left\{\|(v_x)\|_{W,2,1} +\|(v_y)\|_{W,2,1} + {\mu} \|z\|_1:\ z=Hf-g, v_x=\nabla_x f, v_y=\nabla_y f,u=f\right\}
\end{equation}
we can also get convergence theorem similar as Theorem 5, whose detail has been presented in the last section for the constrained TV OGS $L_2$ model. We call this relevant algorithm CATVOGSL1.

%

\section{Numerical results}
In this section, we present several numerical results to illustrate the performance of the proposed method.
All experiments are carried out on a desktop computer using Matlab 2010a. Our computer is equipped with an Intel
Core i3-2130 CPU (3.4 GHz), 3.4 GB of RAM and a 32-bit Windows 7 operation system.
\subsection{Comparison with MM method on the problem (\ref{norm21wogspm})}
In this section, we solve an example of (\ref{norm21wogspm}) with $X=rand(100,100)$ and $X(45:55,45:55)=0$. We only compare the results of our explicit shrinkage formulas with the most recent MM iteration method proposed in \cite{CS2014} as a simple example. Particularly, we set weighted matrix $W_g\in\mathbb{R}^{3\times 3}((W_g)_{i,j}\equiv 1$), and list (\ref{norm21wogspm}) again
\begin{equation}\label{norm21w1ogspm}
  \min_{A} \ \ {\rm Fun}(A)=\ \|A\|_{2,1}+ \frac{\beta}{2}\|A -X\|_F^2.
\end{equation}

\begin{figure}
  \centering
  \subfigure{\includegraphics[width=0.18\textwidth,clip]{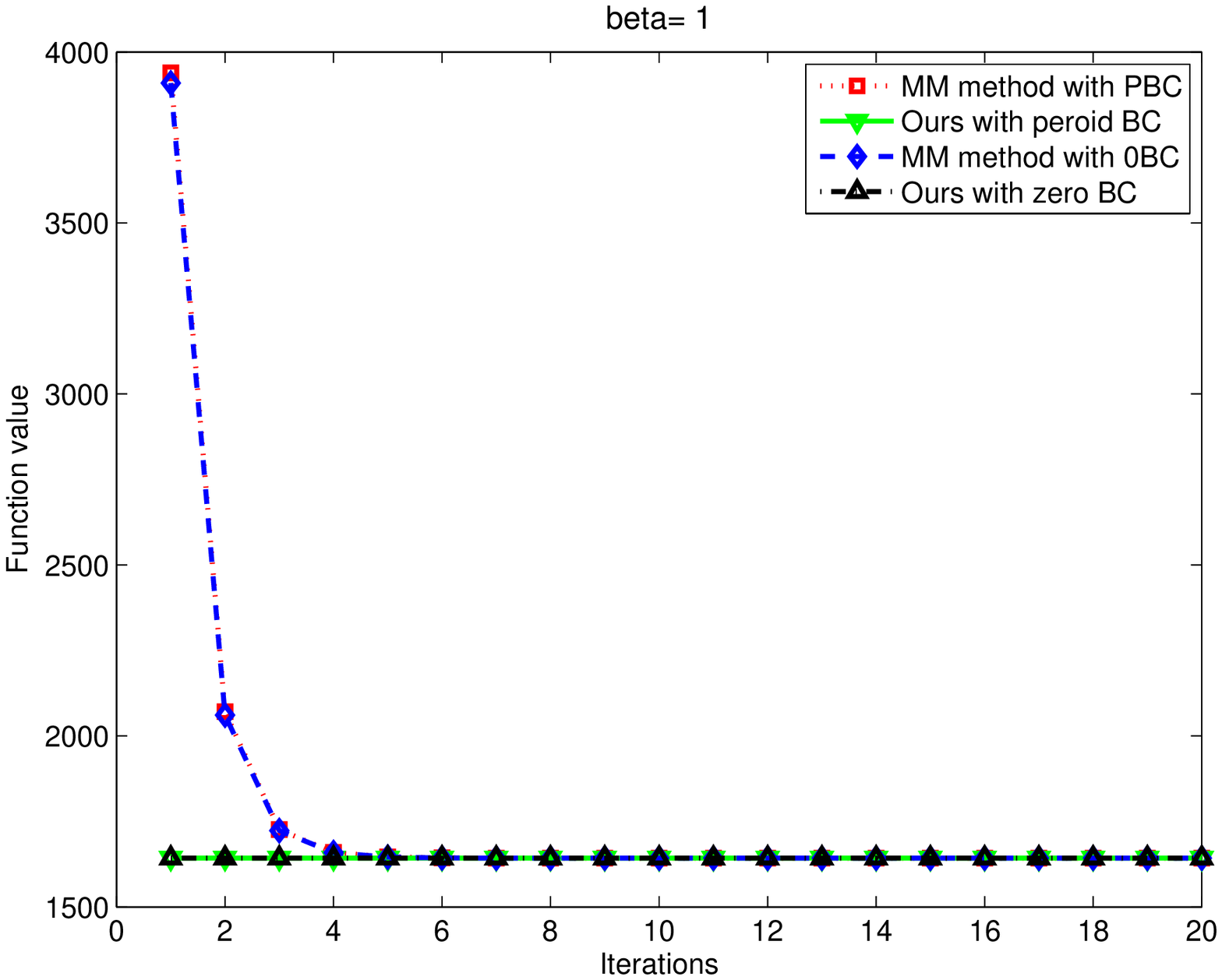}}
  \subfigure{\includegraphics[width=0.18\textwidth,clip]{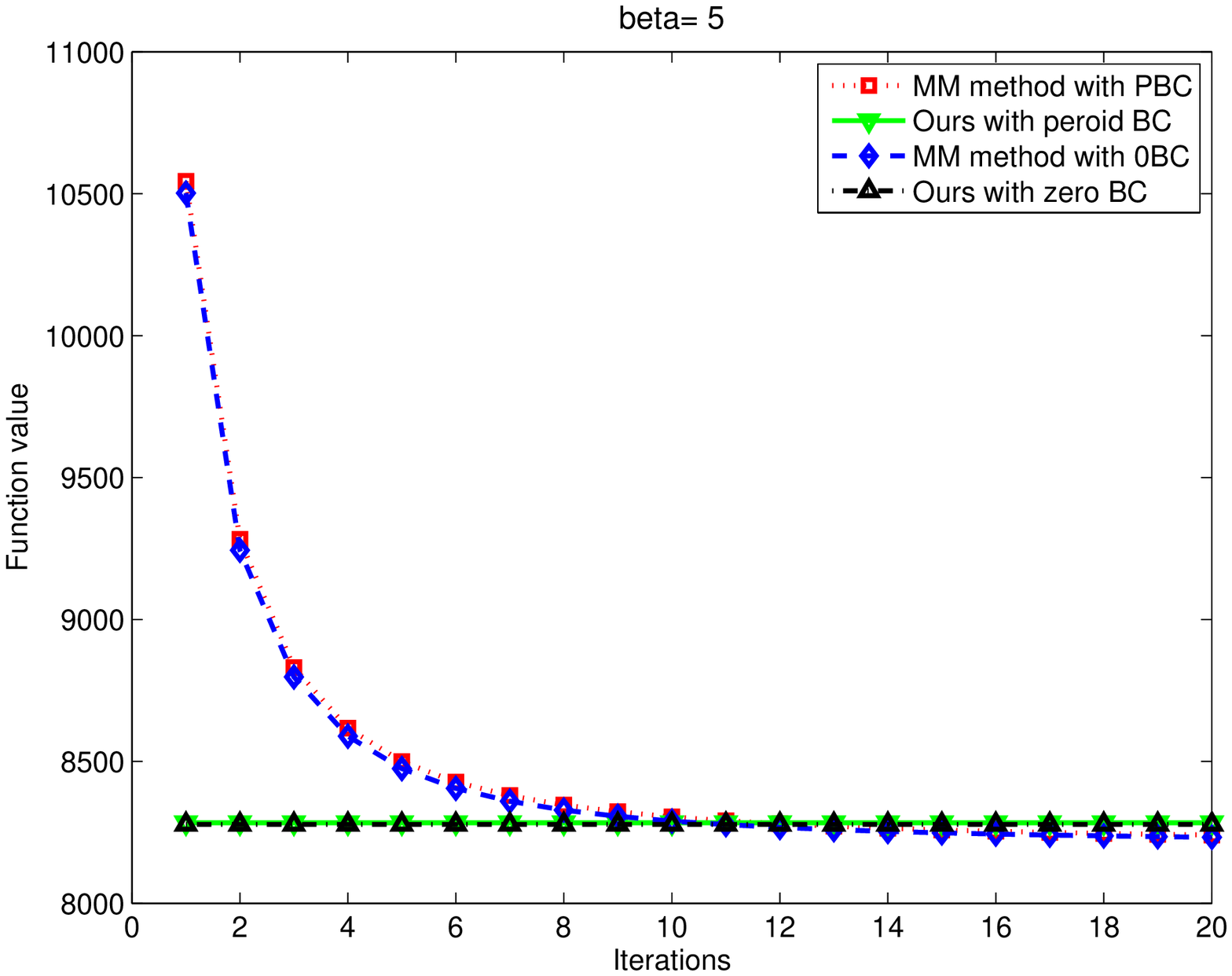}}
  \subfigure{\includegraphics[width=0.18\textwidth,clip]{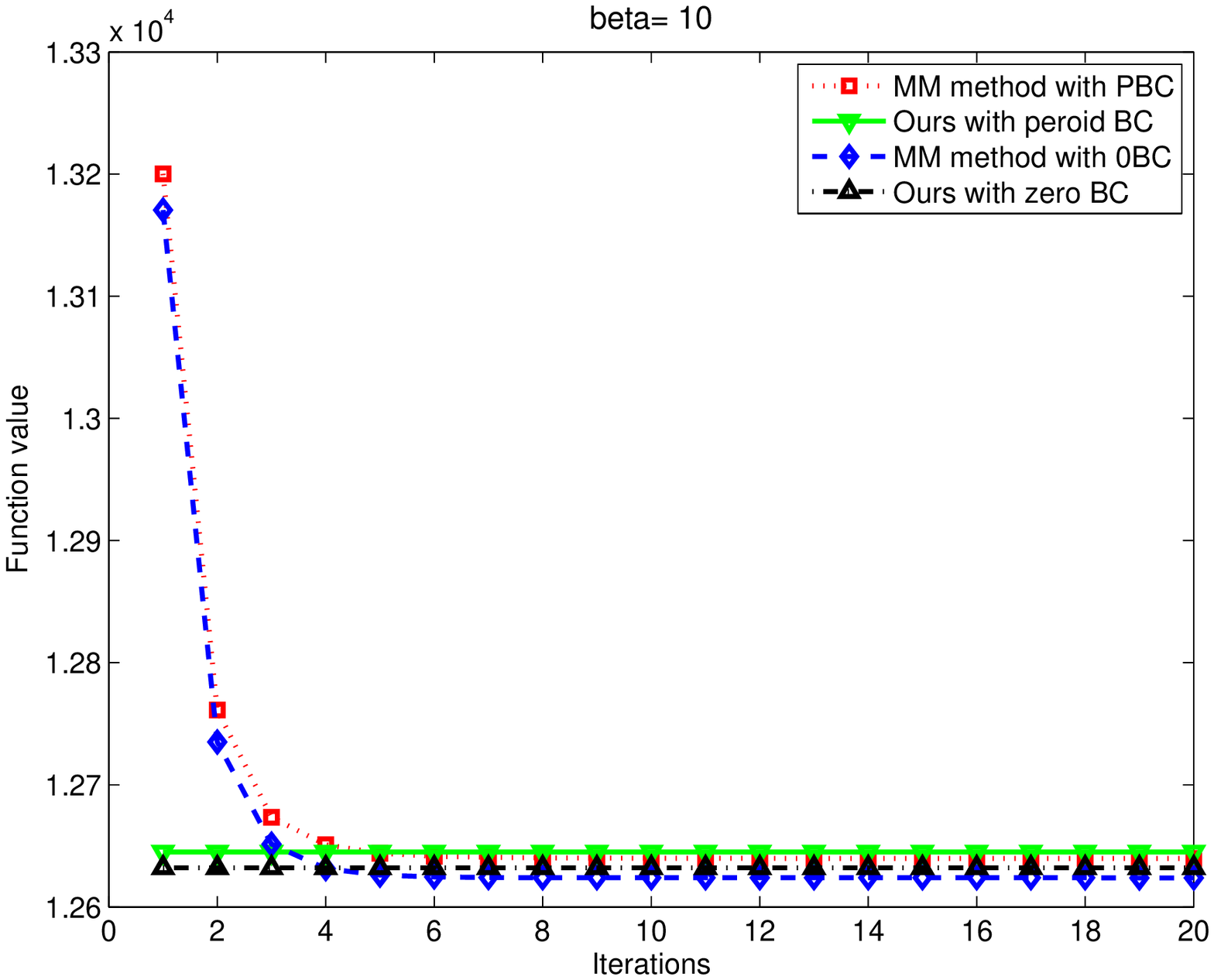}}
  \subfigure{\includegraphics[width=0.18\textwidth,clip]{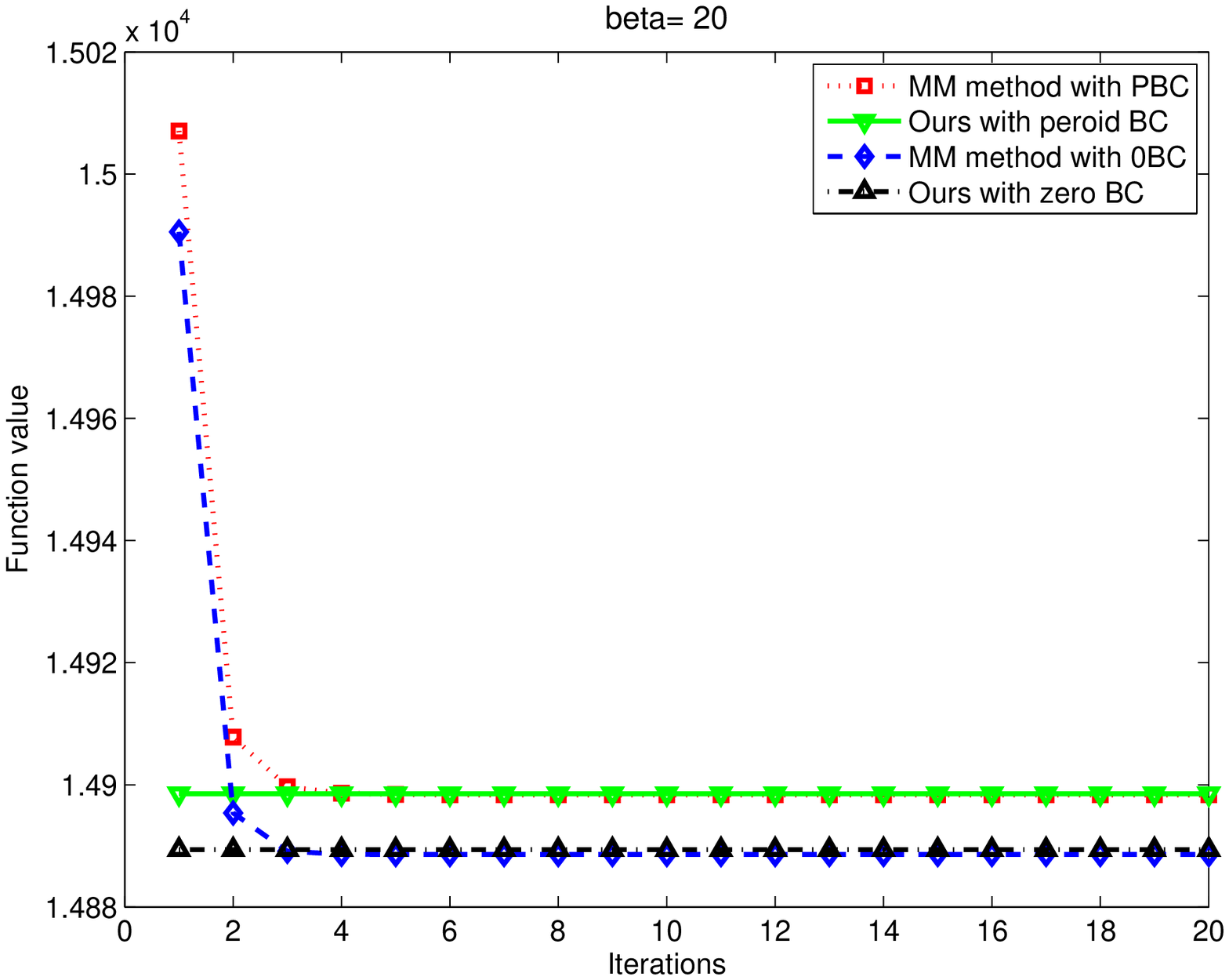}}
  \subfigure{\includegraphics[width=0.18\textwidth,clip]{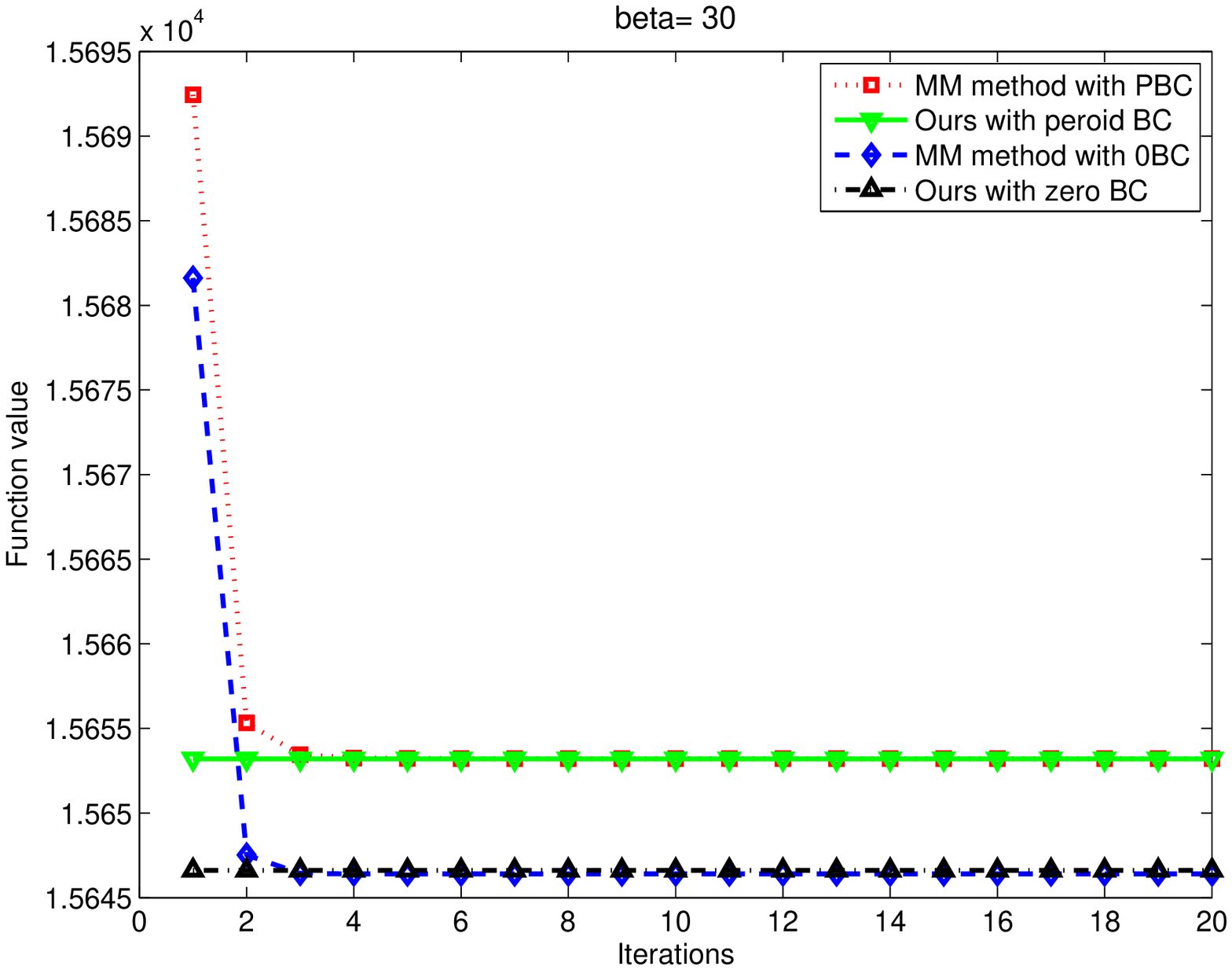}}
  \subfigure{\includegraphics[width=0.18\textwidth,clip]{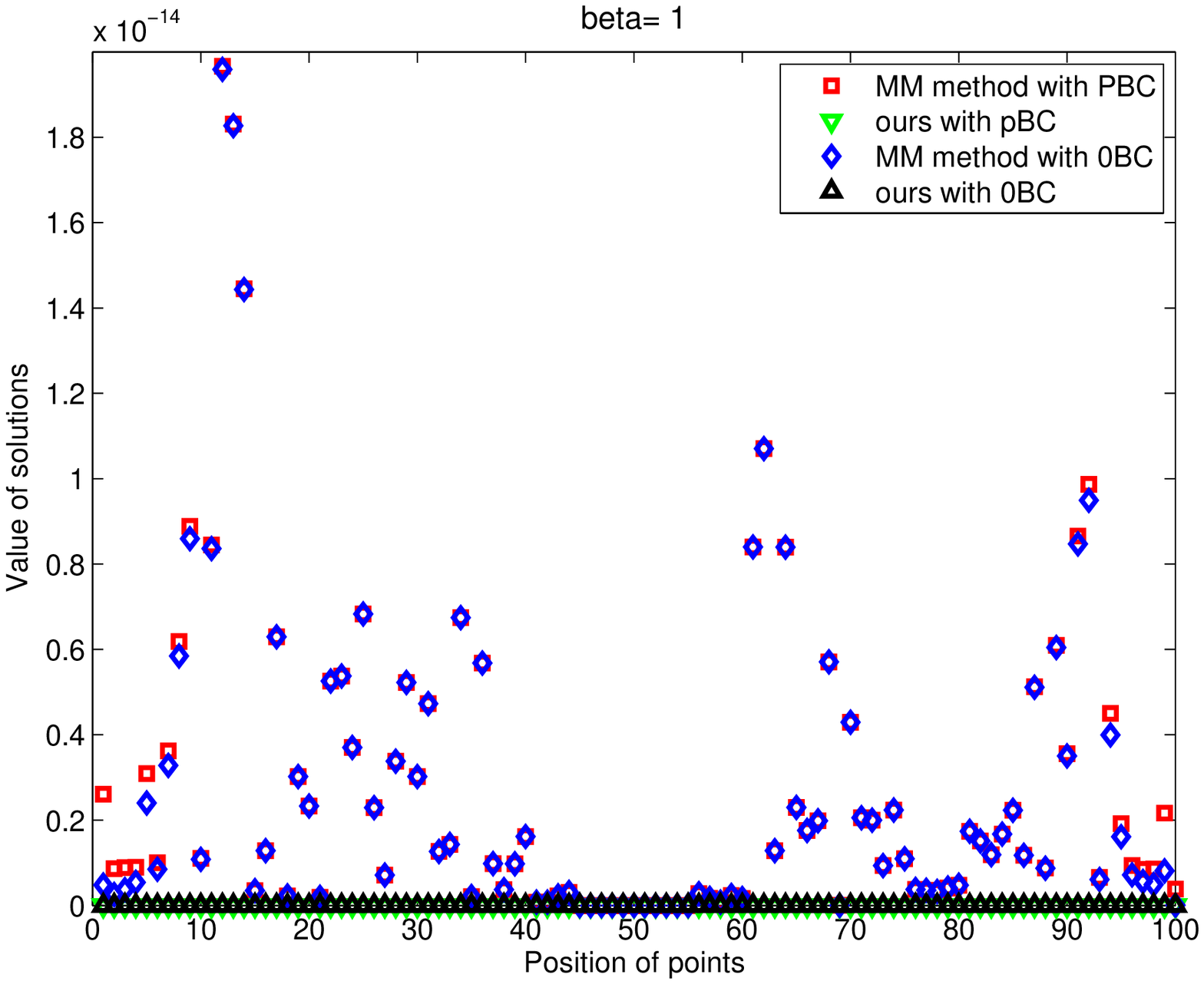}}
  \subfigure{\includegraphics[width=0.18\textwidth,clip]{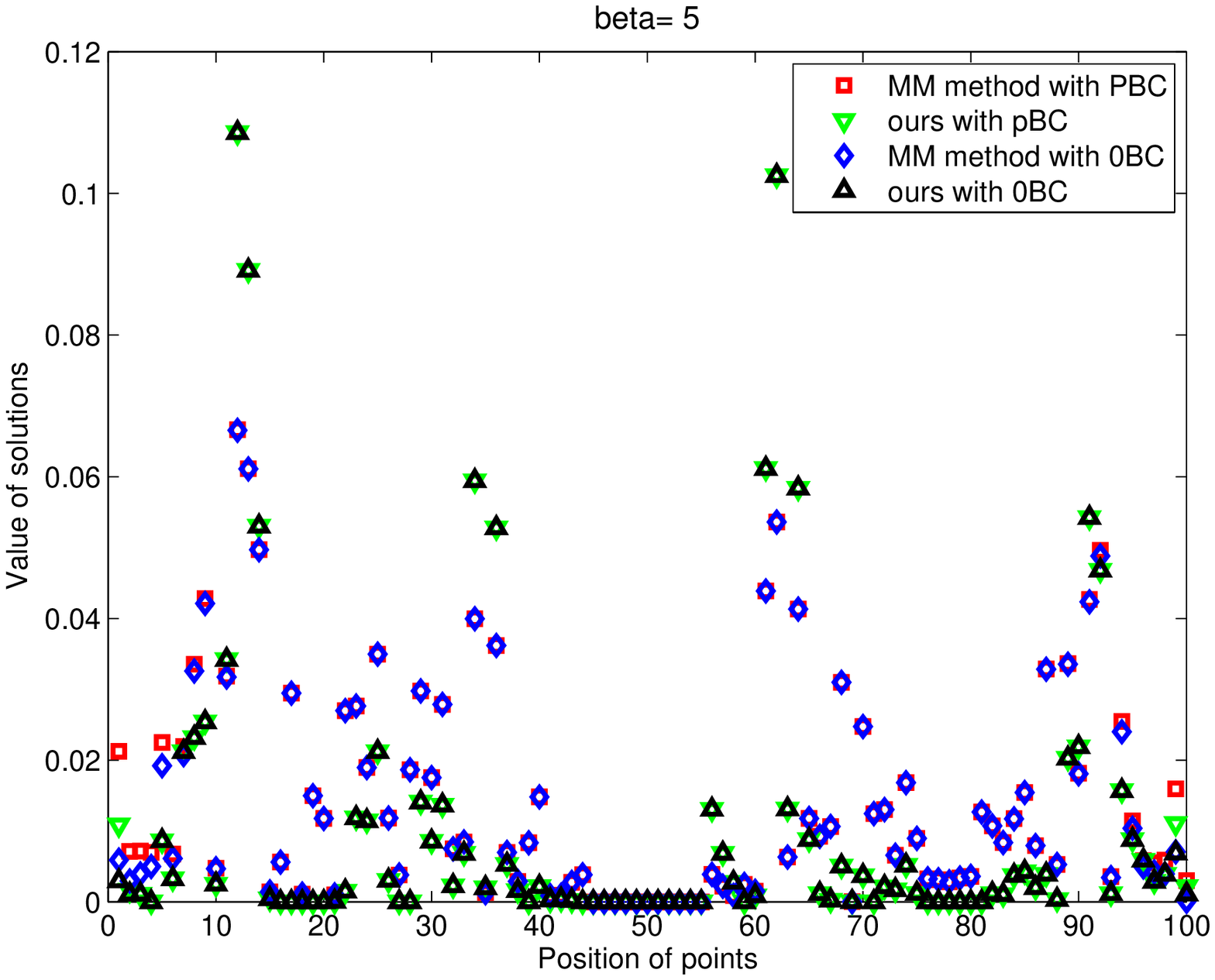}}
  \subfigure{\includegraphics[width=0.18\textwidth,clip]{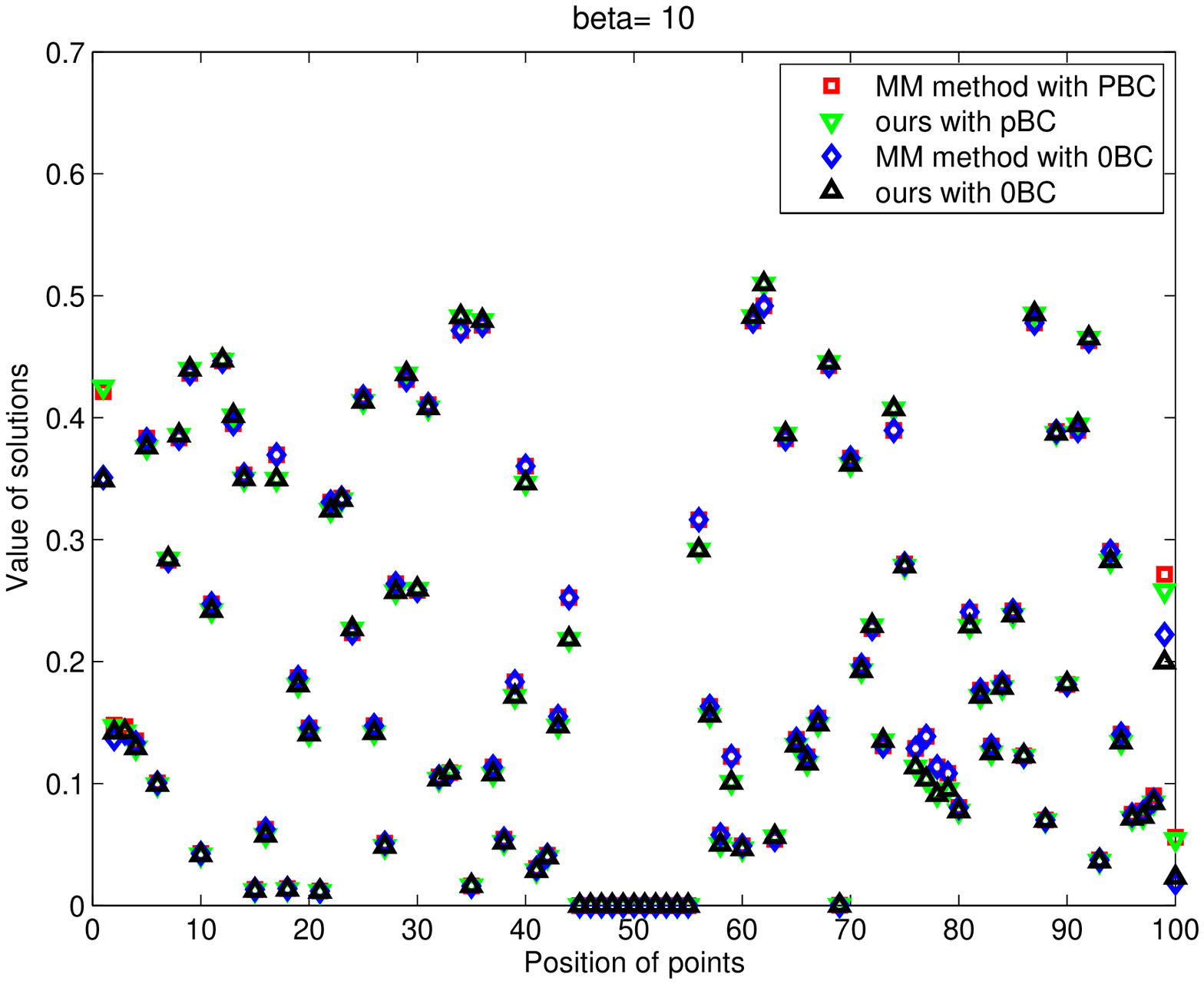}}
  \subfigure{\includegraphics[width=0.18\textwidth,clip]{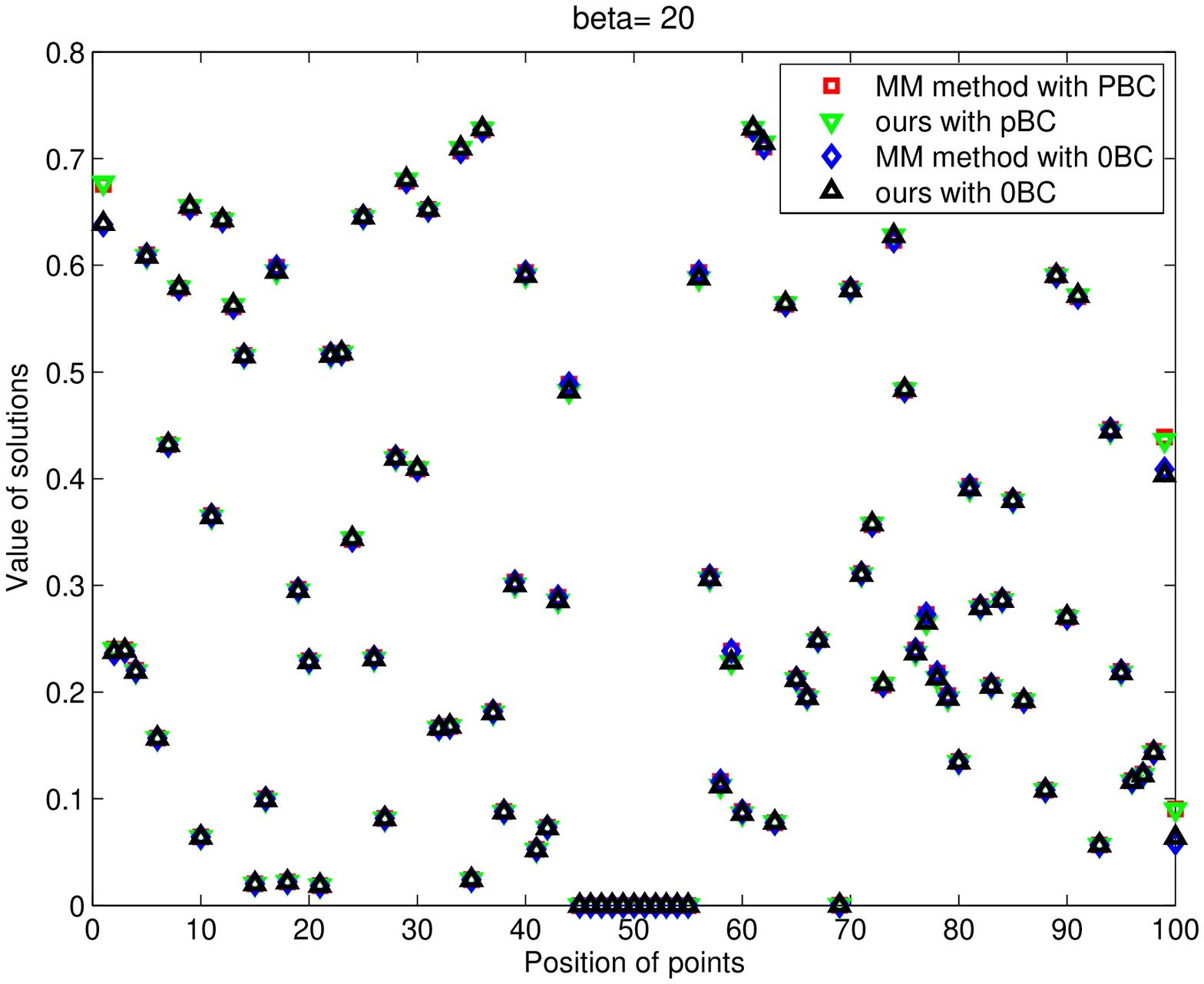}}
  \subfigure{\includegraphics[width=0.18\textwidth,clip]{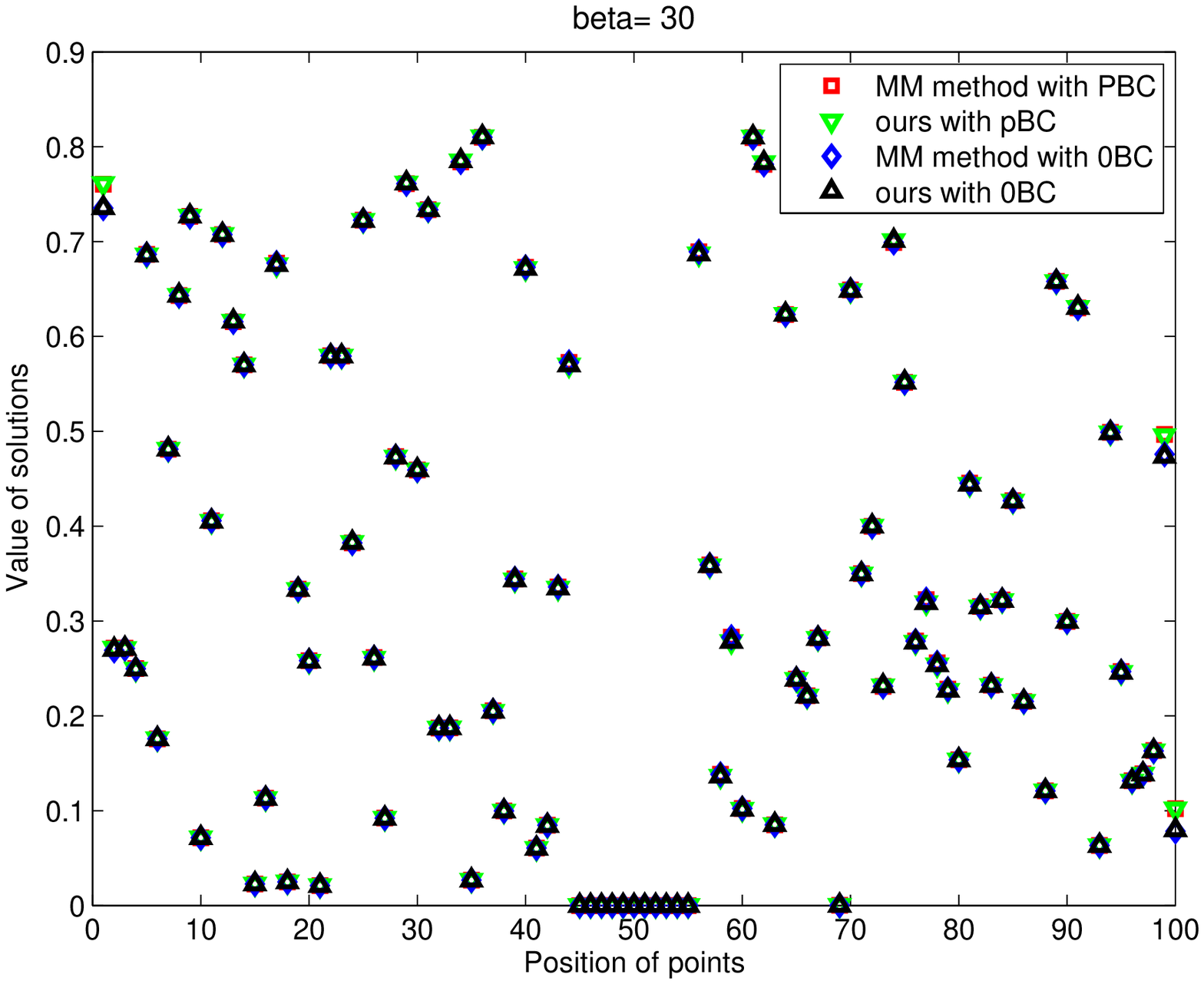}}
  \caption{\text{Comparison between our method and the MM method}.}
  \label{diffbeta}
 \end{figure}

We set different parameter $\beta$ and set the number of MM iteration by 20 steps. For comparison,
we expand our result by explicit shrinkage formulas to length 20 (the same as the MM iteration steps).
The values of the object function Fun($A$) against to iterations are illustrated in top line in Figure~\ref{diffbeta}.
The cross sectional elements of the minimizer $A$ are show in the bottom line in Figure~\ref{diffbeta}.
We choose both zero boundary conditions (0BC) and period boundary conditions (PBC) both for our method and
 the MM method. It is obvious that our method is faster than MM method because of the explicit shrinkage formulas.
 From the top line in Figure~\ref{diffbeta}, we observe that MM method is also fast because it only needs less than 20 steps
 (sometimes 5) to converge. The related error of the function value of two methods is much less than 1 percent
  for different parameters $\beta$. From the bottom line in Figure~\ref{diffbeta}, we can see that when
  $\beta\leq 1$ which is sufficiently small, and when $\beta\geq 30$ which is sufficiently large,
  our result is almost the same as the final results by the MM iteration method. This shows that our method is accurate and the error of the MM method is very small. When $1<\beta<30$,
  the minimizer computed by our method is approximate and the error is much large than
   the MM method although the error of the function value by two methods is very small.
   However, in this case, form the table and the figure, we can see that our method is also able to be seen as a approximate solution.
   Moreover, our method is much more efficient because our time complexity is just the same as one step iteration in the MM method.
In Table~\ref{gau75noi4NItda}, we show the numerical comparison of our method and the MM method on three parts, 
the related error of function value (ReE of  $f_A(X)$), related error of minimizer $X$ (ReE of $X $) and the mean absolute error of $X$ (MAE of $X$).
From the table, we can see that our formula can almost get same results as the MM method  when $\beta$ is sufficiently large, and approximate results similar as the MM method for other $\beta$. This is another proof for the feasibility of our formula.

   After more than one thousand similar tests on the example $X=rand(100,100)$,
   we find that when $\beta\geq 30$ which is sufficiently large, our method is always
   accurate and the minimizers by our method and the MM method (with 20 steps) are almost the same.
   Therefore, form Theorem 1 and Theorem 2, we choose L=30 (sufficiently large)  and choose a suit
    parameter $\beta_1$ in the next section to make sure the convergence. Moreover, we can see the efficiency of our method in the next section.

   Similarly as weighted matrix $W_g\in\mathbb{R}^{3\times 3}((W_g)_{i,j}\equiv 1$), we also tested other weighted matrix for more than one thousand times.
   For the examples when $X=rand(100,100)$ and $X(X>=0.5)=0$ (or $X(X<=0.5)=1$) (every element is in $[0,1]$), we find that, generally,
    \(\beta\leq \frac{\|w_g\|_2}{\sqrt{s}}\) is sufficiently small and \(\beta\geq 30\cdot\frac{\|W_g\|_2}{\sqrt{s}}\) is sufficiently large.
    This results illustrate that the former when $W_g\in\mathbb{R}^{3\times 3}((W_g)_{i,j}\equiv 1$) $\beta\geq 30$ which
     is sufficiently large again.
    However, in practice, if \(\beta\) is too large, then $A=X$,
     and the minimization problem (\ref{norm21w1ogspm}) is meaningless. In our experiments after more than one thousand tests, we find that, when every element of $X$ is in $[0,1]$,
     \(30\cdot\frac{\|W_g\|_2}{\sqrt{s}}\leq\beta\leq 300\cdot\frac{\|W_g\|_2}{\sqrt{s}}\) is sufficiently large
      but not too large to make the minimization problem meaningless. In this case, that means our formula is useful and effective in practice.
     On the other hand, when some elements of $A$ are not in $[0,1]$, we can first project or stretch $A$
     into the region $[0,1]$, and then choose the sufficiently small or large parameter $\beta$.

 \begin{table*}[hhbp]
\renewcommand{\captionlabeldelim}{.}
\setlength{\abovecaptionskip}{0pt}
\setlength{\belowcaptionskip}{10pt} \centering \caption{Comparison of our method and the MM method for two kind of boundary conditions (BC), zero boundary conditions (0BC) and period boundary conditions (PBC).} \centering
\begin{tabular}{|c|c|c|c|c|c|c|c|c|c|c|}
\hline 
  BC             &$\beta$=&1& 5&7&10&15&20&30&50   \\
\hline
\multirow{3}{*}{0BC}& ReE of $f_A(X)$ &5.9e-14&0.0055&0.0028&6.5e-4&1.4e-4&5.4e-5&1.4e-5&2.6e-6\\
 \cline{2-10}
                &ReE of $X $&---&---&0.0207&0.0017&1.9e-4&4.7e-5&7.5e-6&7.8e-7 \\
 \cline{2-10}
                                  &MAE of $X$ &3.4e-15&0.0094&0.0130&0.0067&0.0029&0.0016&6.9e-4&2.4e-4\\
\hline
\multirow{3}{*}{PBC}& ReE of $f_A(X)$ & 6.5e-14&0.0053&0.0025&4.1e-4&5.5e-5&1.1e-5&6.3e-7&1.3e-6 \\
\cline{2-10}
                   &ReE of $X$ &---&---&0.0193&0.0014&1.3e-4&2.9e-5&4.3e-6&4.5e-7 \\
 \cline{2-10}                &MAE of $X$ &3.8e-15&0.0097&0.0130&0.0063&0.0026&0.0014&5.9e-4&2.1e-4\\
\hline
\end{tabular}
\label{gau75noi4NItda}\end{table*}


\subsection{Comparison with TV methods and TV with OGS with inner iterations MM methods for image deblurring and denoising}
\begin{figure}
  \centering
 \subfigure[Man]{\includegraphics[width=0.24\textwidth,clip]{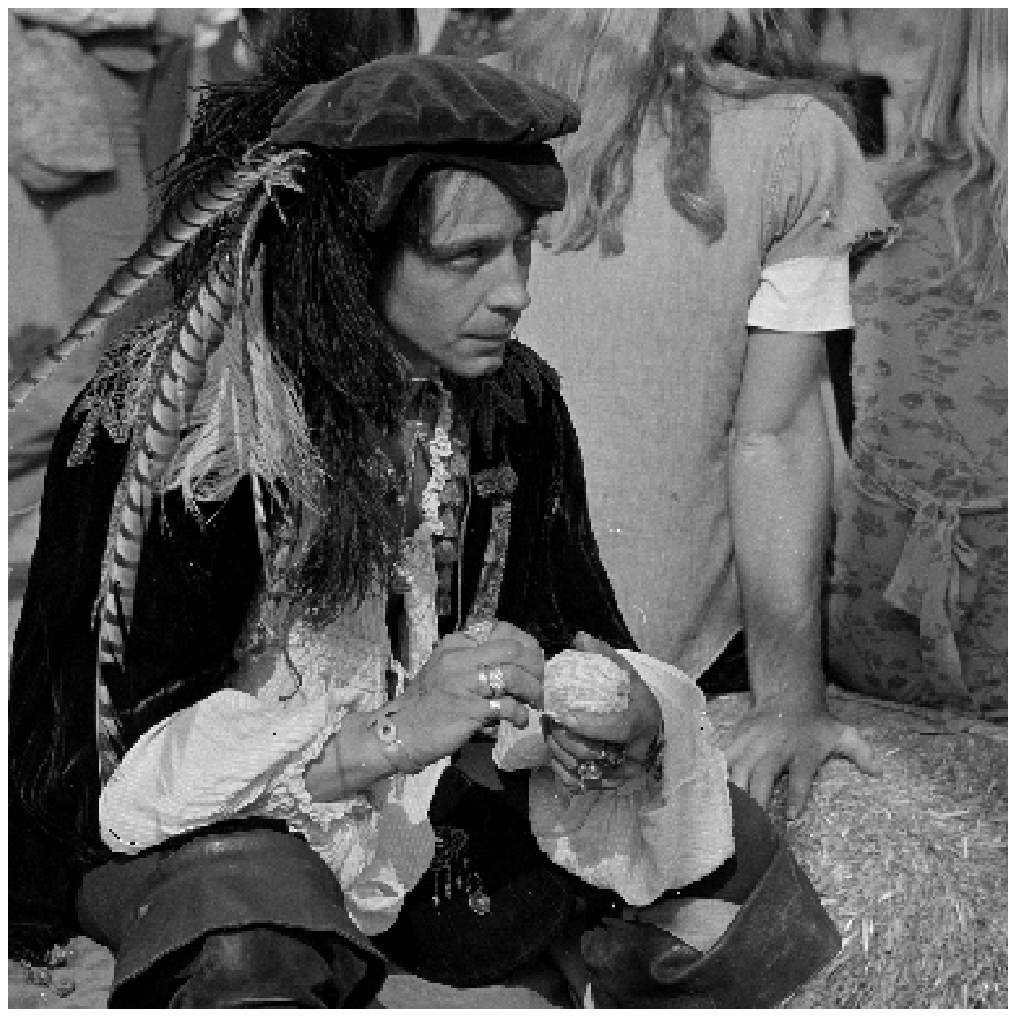}}
   \subfigure[Car]{\includegraphics[width=0.24\textwidth,clip]{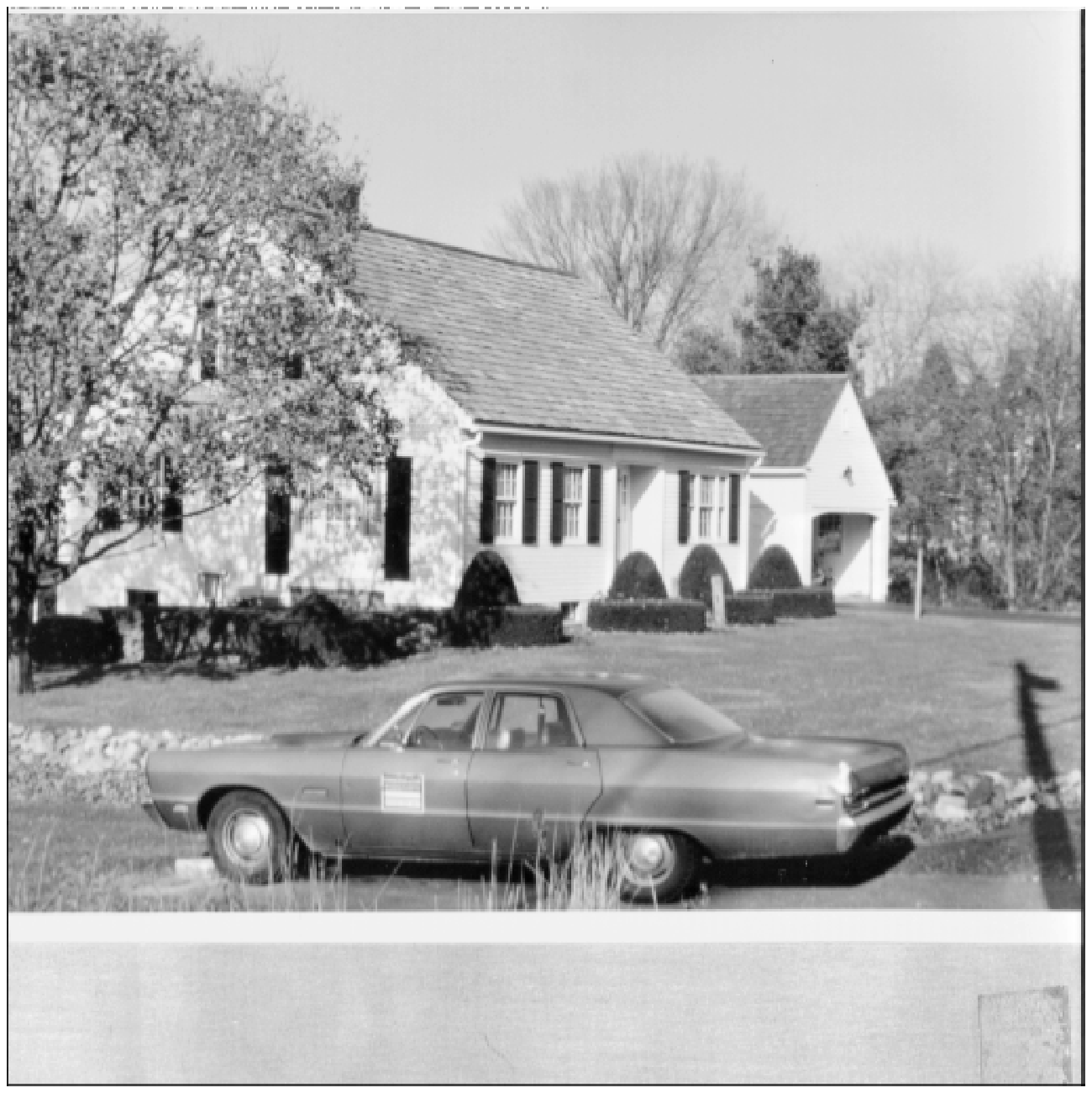}}
 \subfigure[Parlor]{\includegraphics[width=0.24\textwidth,clip]{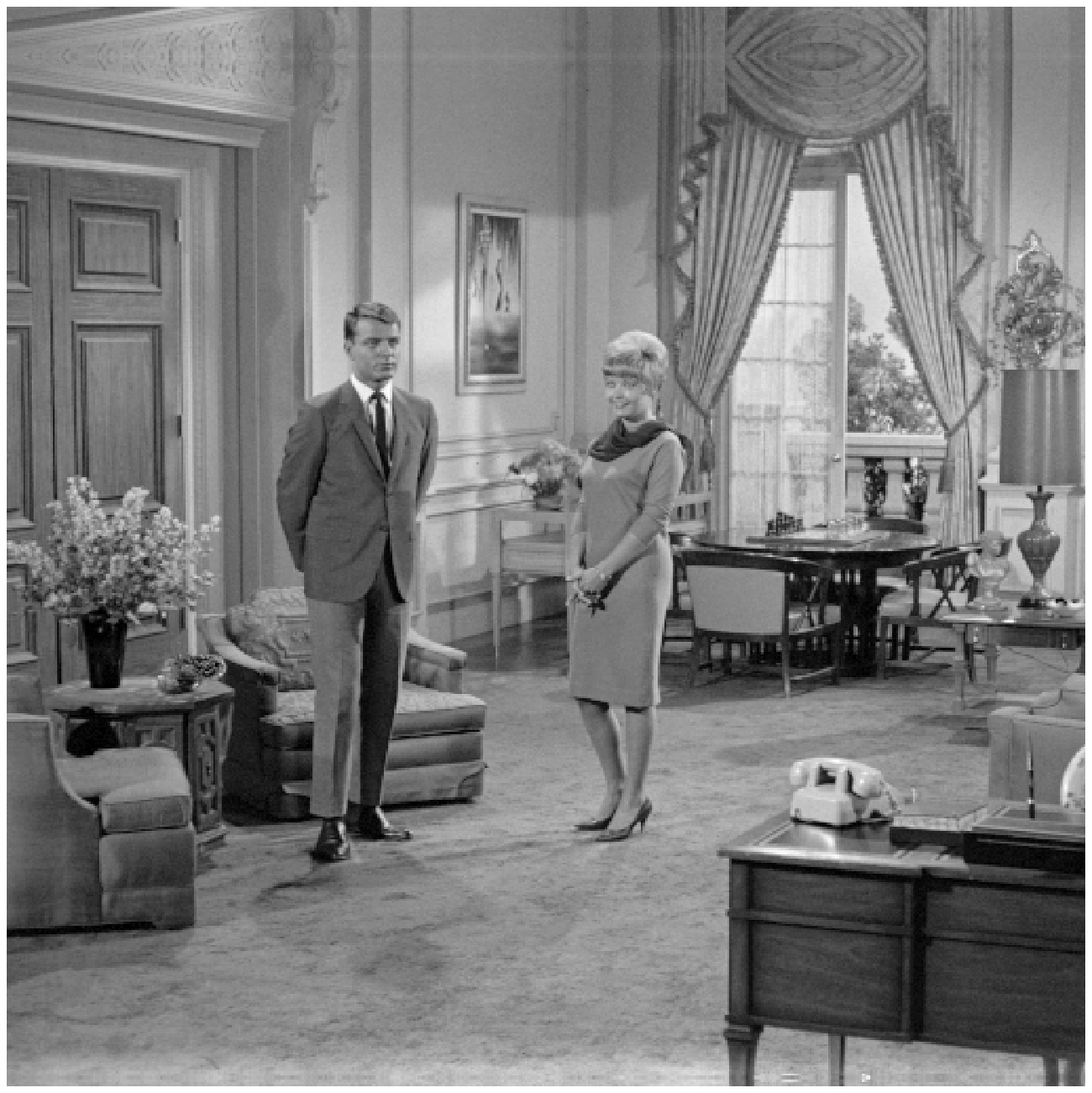}}
 \subfigure[Housepole]{\includegraphics[width=0.24\textwidth,clip]{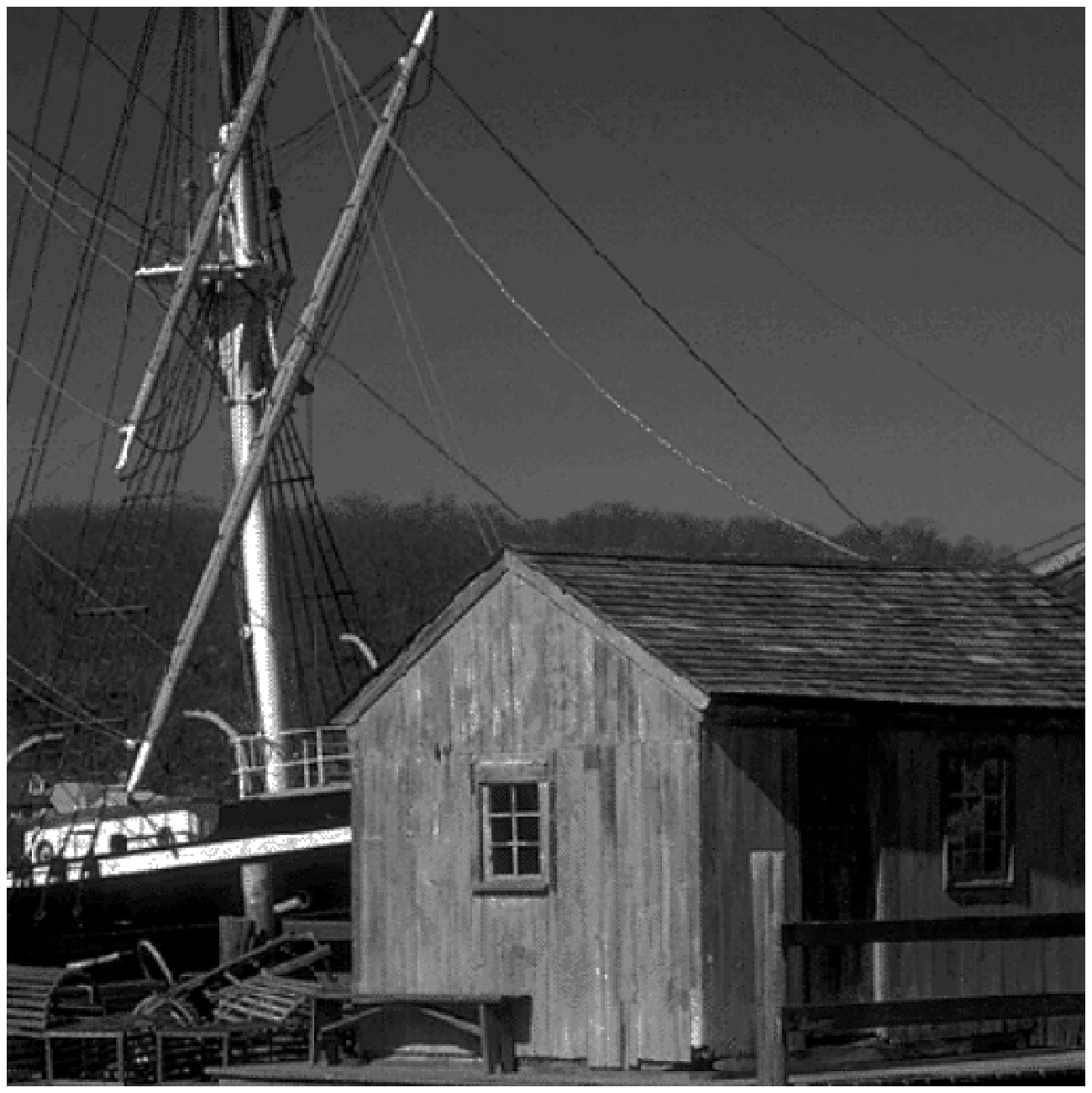}}
 \caption{\text{Original images}.}
  \label{originalimages}
 \end{figure}

In this section, we compare all our algorithms with other methods. All the test images are shown in Fig.~\ref{originalimages}, one 1024-by-1024 image as (a) Man, and three 512-by-512 images as: (b) Car, (c) Parlor, (d) Housepole.

The  quality  of  the  restoration  results  is  measured  quantitatively
by  using  the  peak signal-to-noise ratio (PSNR) in decibel (dB) and the  relative  error  (ReE):
\begin{equation*}
\textrm{PSNR} = 10\log_{10} \frac{n^2 {\rm Max}_I^2}{\|f-\bar{f}\|_2^2},\quad
\textrm{ReE} = \frac{\|f -\bar{f}\|_2}{\|\bar{f}\|_2},
\end{equation*}
where $\bar{f}$ and $f$ denote the original and restored images
respectively, and ${\rm Max}_I$ represents the maximum possible pixel value of the image.
In our experiments,  ${\rm Max}_I= 1$.

The stopping criterion used in our work is set to be as other methods
\begin{equation}
 \frac{|\mathcal{F}^{k+1}-\mathcal{F}^k|}{|\mathcal{F}^k|}< 10^{-5},
\end{equation}
where $\mathcal{F}^k$ is the objective function value of the respective model in the $k$th iteration.

we compare our methods with some other methods, such as Chan's TV method proposed in \cite{CTY2013} , Liu's method proposed in \cite{LHS2013}, and Liu's method proposed in \cite{LHLL2013}. Both the latter two methods \cite{LHS2013} and \cite{LHLL2013} are with inner iterations MM methods for OGS TV problems, where the number of the inner iterations is set 5 by them. Particularly, for a fair comparison, we set weighted matrix $W_g\in\mathbb{R}^{3\times 3}((W_g)_{i,j}\equiv 1$) in all the experiments of our methods as in \cite{LHS2013,LHLL2013}.

\subsubsection{Experiments for the constrained TV OGS $L_2$ model}
In this section, we compare our methods (CATVOGSL2 and CITVOGSL2) with some other methods, such as Chan's method proposed in \cite{CTY2013} (Algorithm 1 in \cite{CTY2013} for the constrained TV-L2 model) and Liu's method proposed in \cite{LHS2013}.

Regarding the penalty parameters $\beta$＊s in all our algorithms, theoretically any positive
values of  $\beta$＊s ensure the convergence. Therefore, we set the penalty parameters $\beta_1=35$, $\beta_2=20$, for the ATV case, $\beta_1=100$, $\beta_2=20$ for the ITV case
and relax parameter $\gamma=1.618$. The blur kernels are generated by Matlab built-in function (i)\texttt{fspecial('average',9)} for $9\times 9$  average blur. We
generate all blurring effects using the Matlab built-in function
\texttt{imfilter(I,psf,} \texttt{'circular','conv')} under periodic boundary conditions with ``\texttt{I}'' the original image and ``\texttt{psf}'' the blur kernel. We first generated the blurred images operating on images (a)-(c) by the former Gaussian blurs
and further corrupted by zero
mean Gaussian noise with BSNR = 40. The BSNR is given by
\begin{equation*}
\textrm{BSNR} = 20\log_{10}\frac{\|g\|_2}{\|\eta\|_2},
\end{equation*}
where $g$ and $\eta$ are the observed image and the noise, respectively.
\begin{table}[!hbp]
\setlength{\abovecaptionskip}{0pt}
\setlength{\belowcaptionskip}{10pt} \centering \caption{Numerical comparison of Chan. \cite{CTY2013}, Liu. \cite{LHS2013}, CATVOGSL2 and CITVOGSL2 for images (a)--(c) in Figure~\ref{originalimages}. PSNR: dB, Time: s.} \centering\tiny{\small{
\begin{tabular}{|@{~}c@{~}|@{~}c@{~}|@{~}r@{/}l@{/}r@{/}l@{~}|@{~}r@{/}l@{/}r@{/}l@{~}|@{~}r@{/}l@{/}r@{/}l@{~}|}
\hline 
 \multicolumn{2}{|c}{Images}&\multicolumn{4}{c}{(a) Man} & \multicolumn{4}{c}{(b) Car}&  \multicolumn{4}{c|}{(c) Parlor} \\
 \hline
 \multirow{1}{*}{Method} &\multirow{1}{*}{$\mu$ ($\times 10^5$)}  & Itrs&PSNR&Time&ReE & Itrs&PSNR&Time&ReE& Itrs&PSNR&Time&ReE \\
 \hline
  \multirow{1}{*}{Chan.\cite{CTY2013}} & 0.5 &15&30.34&6.02&0.0730     &11&31.13&\text{2.00}&0.0426       &12&31.70&\text{2.01}&0.0511     \\
   \hline
  \multirow{1}{*}{Liu.\cite{LHS2013}} & 1&12&30.60&9.06&0.0708  &  12&31.98&\text{3.93}&0.0386    &13&32.46&\text{4.21}&0.0464    \\
   \hline
  \multirow{1}{*}{CATVOGSL2} &1&7&30.62&\textbf{3.37}&0.0711      &9&31.68&\textbf{1.76}&0.0399       &7&32.40&\textbf{1.34}&0.0472   \\
   \hline
  \multirow{1}{*}{CITVOGSL2} &1&8&30.59&\textbf{3.70}&0.0716     &11&31.60&\textbf{2.04}&0.0403        &8&32.03&\textbf{1.61}&0.0492        \\
         \hline
\end{tabular}}}
\label{cameranmantol2}\end{table}

\begin{figure}
  \centering
 \subfigure{\includegraphics[width=0.3\textwidth,clip]{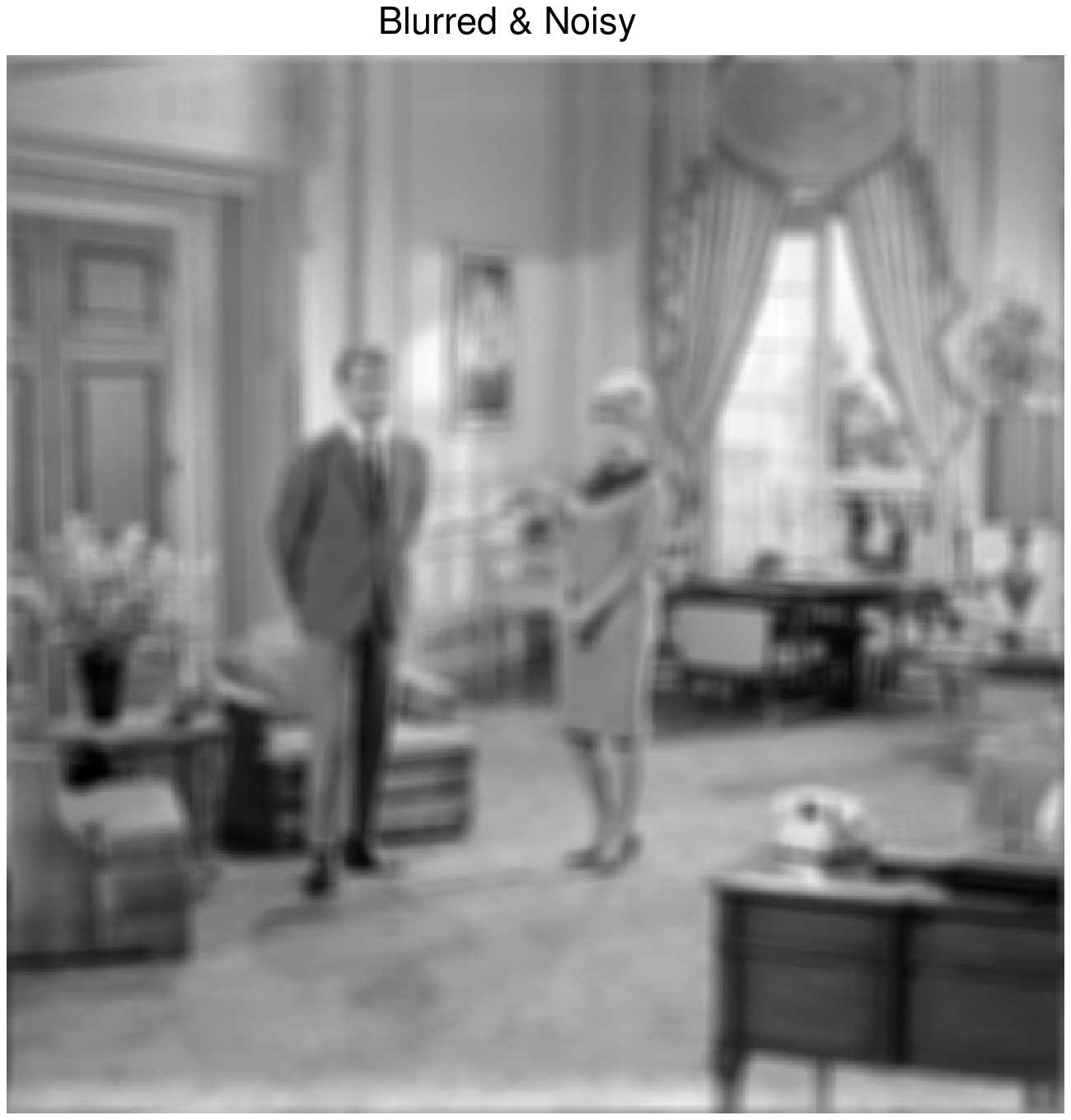}}
 \subfigure{\includegraphics[width=0.3\textwidth,clip]{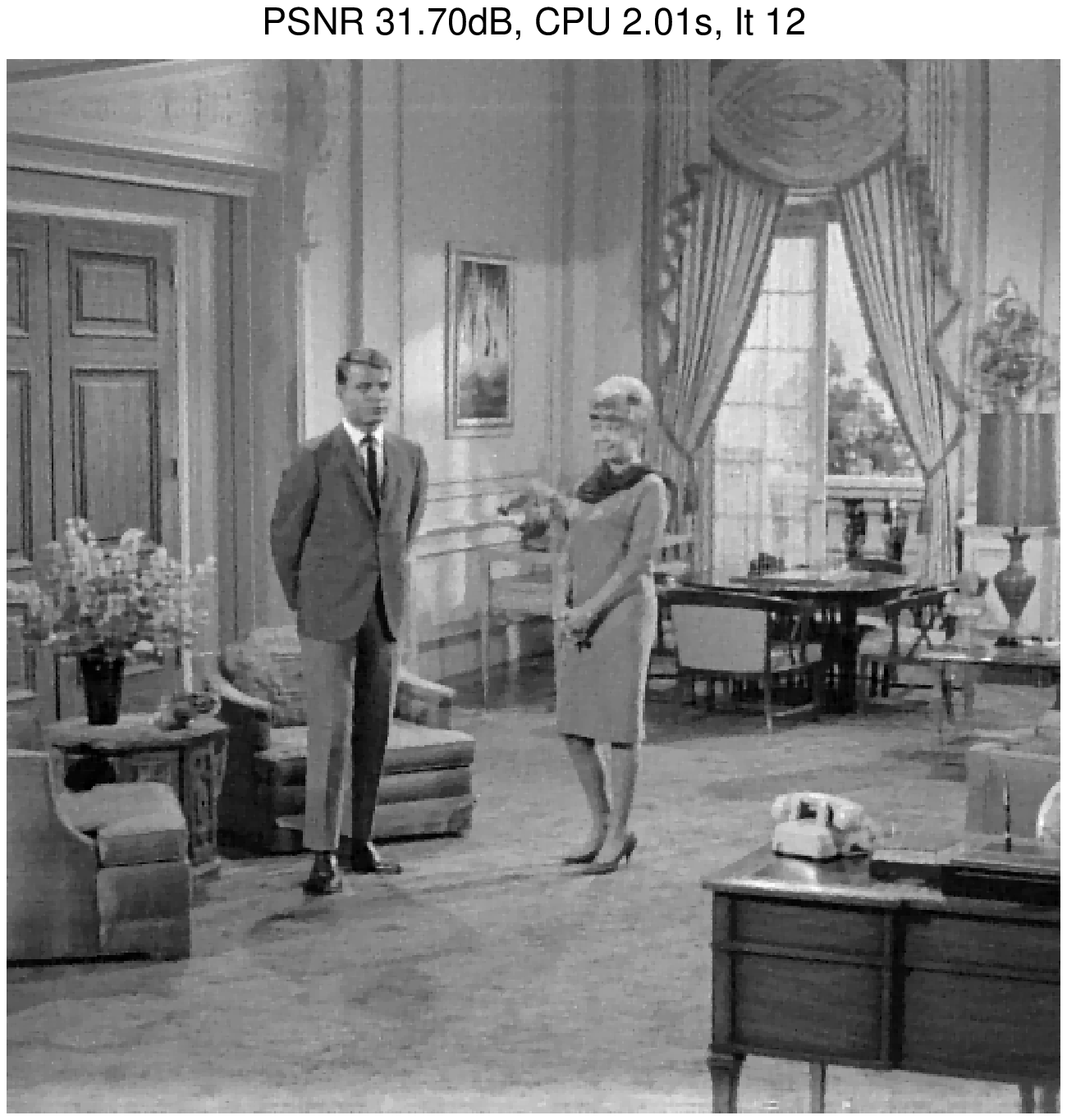}}
   \subfigure{\includegraphics[width=0.3\textwidth,clip]{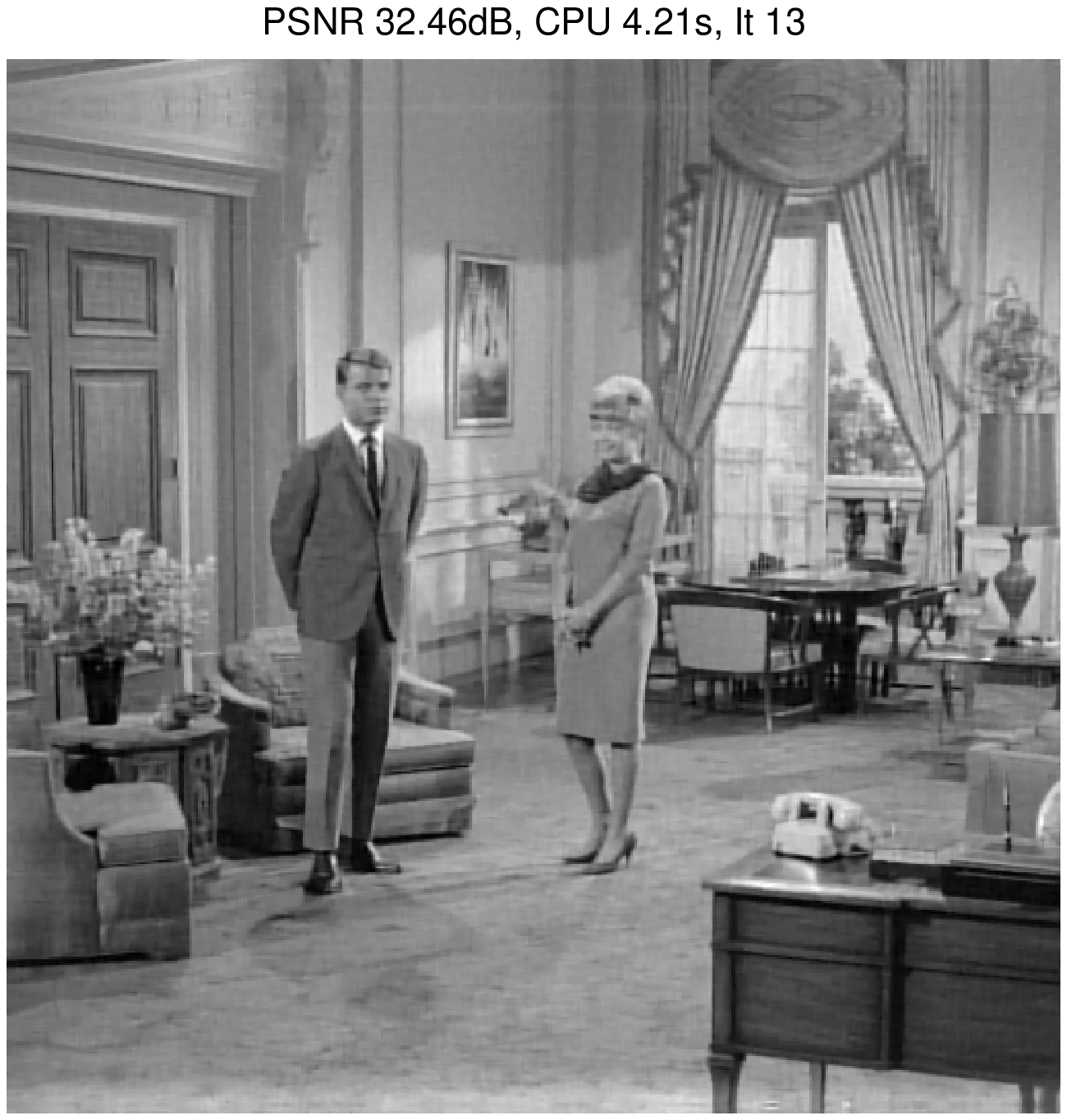}}
 \subfigure{\includegraphics[width=0.3\textwidth,clip]{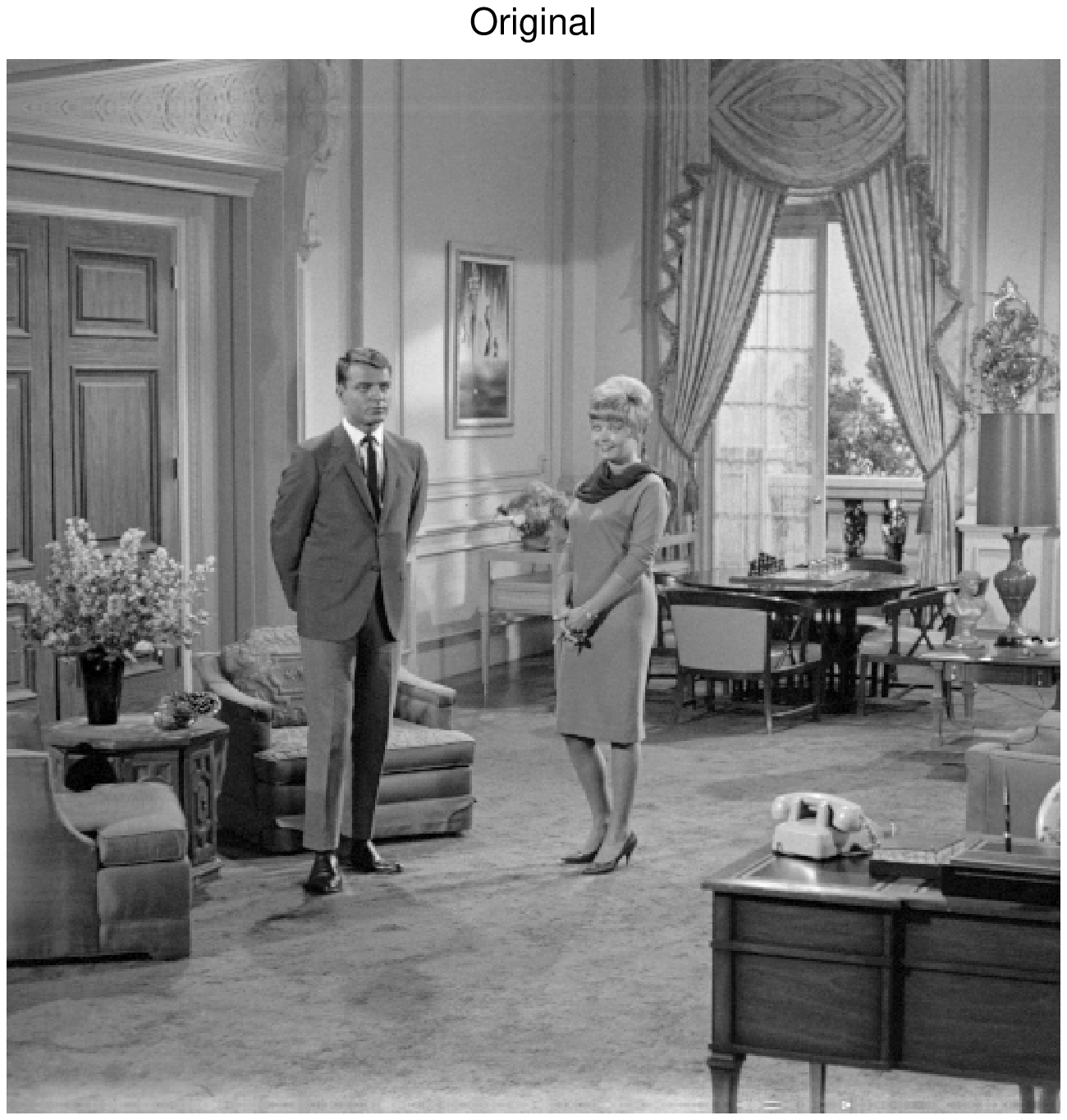}}
 \subfigure{\includegraphics[width=0.3\textwidth,clip]{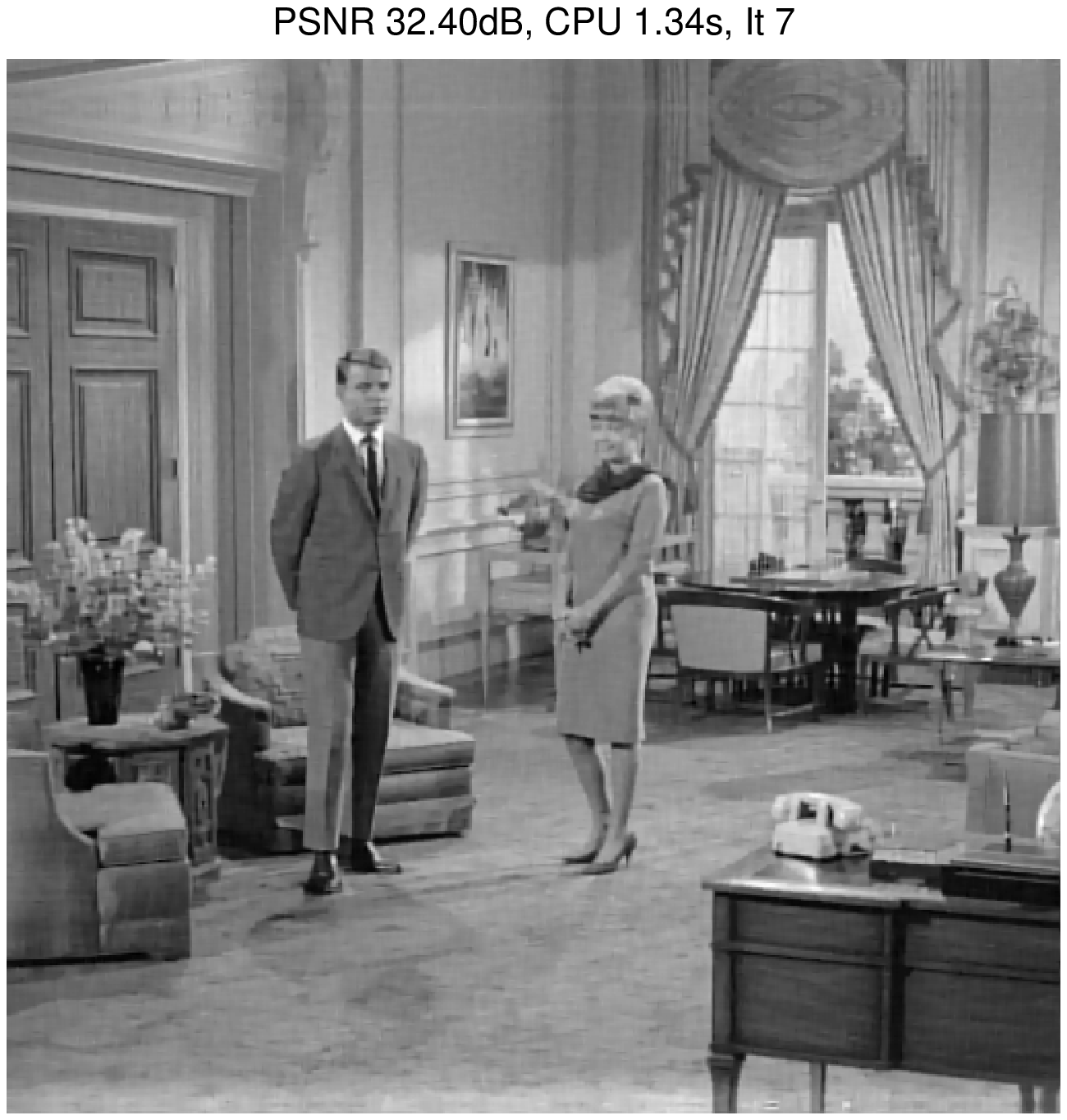}}
 \subfigure{\includegraphics[width=0.3\textwidth,clip]{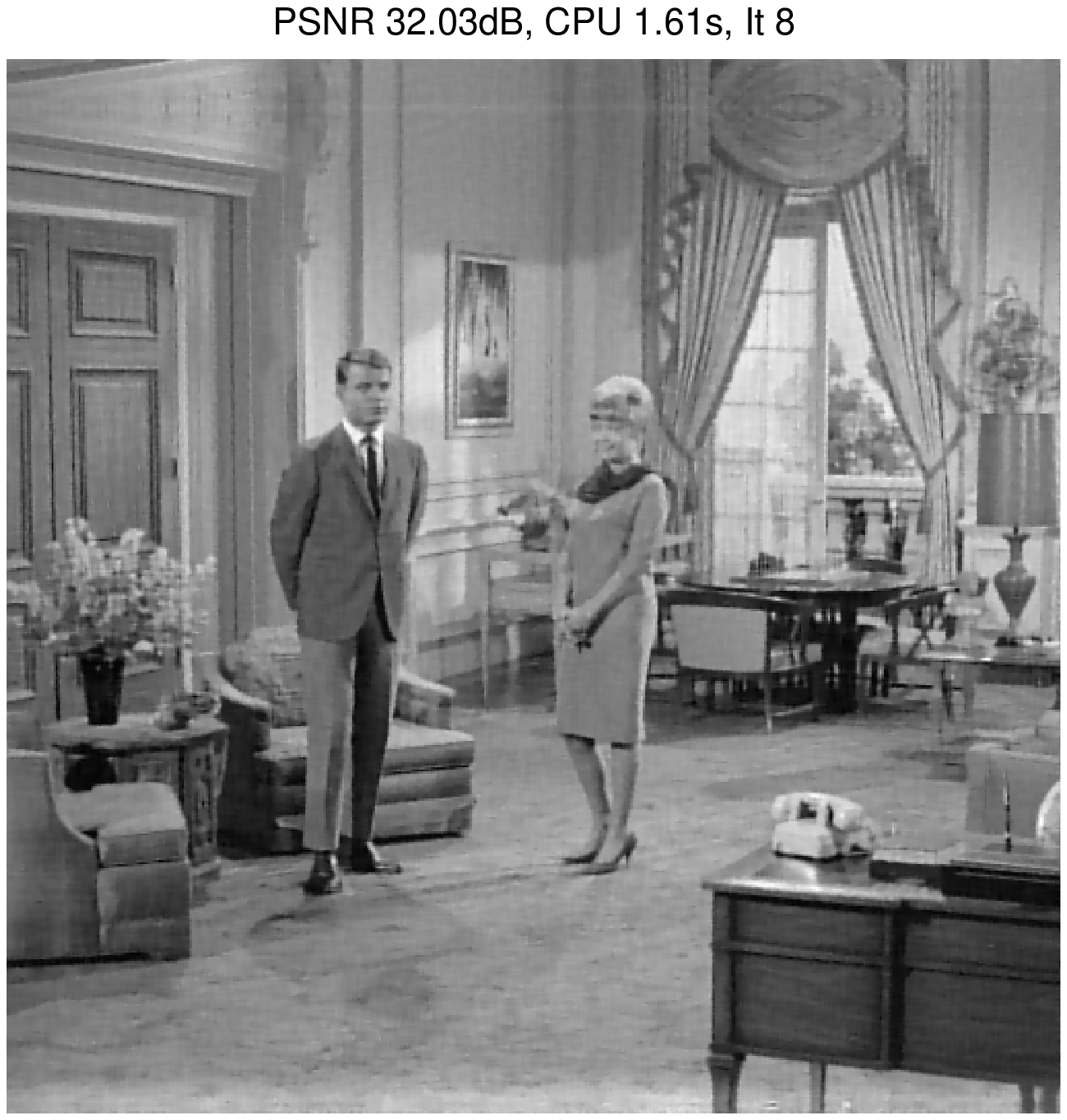}}
 \caption{Top row: blurred and noisy image (left), restoration images of Chan. \cite{CTY2013} (middle), Liu. \cite{LHS2013}. Bottom row: original image (left), restoration images of CATVOGSL2, CITVOGSL2.}
  \label{conparisonofoffice}
 \end{figure}

The numerical results of the three methods are shown in Table~\ref{cameranmantol2}. We have tuned the parameters for all the methods as in Table~\ref{cameranmantol2}. From Table~\ref{cameranmantol2}, we can see that the PSNR and ReE of our methods (both ATV and ITV cases) are almost same as Liu. \cite{LHS2013}, which used MM inner iterations to solve the subproblems (\ref{COTVSUB1}) and (\ref{COTVSUB11}) (only for the ATV case). However, each outer iteration of our methods is nearly twice faster than Liu. \cite{LHS2013} from the experiments. The time of each outer iteration of our methods is almost the same as the tradition TV method in Chan. \cite{CTY2013}. In Figure~\ref{conparisonofoffice} we display the restored ``Parlor'' images from different algorithms. We can see that OGS TV regularizers can get clearer edges on the desk of the image than TV regularizer.

Now we compute the complexity of each step of our methods and Liu's method \cite{LHS2013}. Firstly,
we know that the complexity of all the methods is $512\times 512\times 18 \times 4$ (4 times $n\log_2 n$) except the OGS subproblems.
Then, for the OGS subproblems, the MM method in \cite{LHS2013} with 5 steps inner iteration, the complexity is $512\times 512\times 90 \times 2$ (2 subproblems).
The complexity of our methods is $512\times 512\times 18 \times 2$ in CATVOGSL2 and $512\times 512\times 18$ in CITVOGSL2.
Therefore, the total complexity of \cite{LHS2013} is $512^2 \times 252$, the total complexity of CATVOGSL2 is $512^2 \times 108$ and the total complexity of CITVOGSL2 is $512^2 \times 90$.
That means each step of our methods is more than double faster than the inner iteration method \cite{LHS2013}. In the next section,
the common computation parts of our methods and the inner iteration method are much more, and then our methods are only nearly double faster.\\
\textbf{Remark 5}. Here we do not list the results of Image (d) because they are almost the same as Image (b) and (c).
Moreover, when $\beta_1<30$ which is not sufficiently large, the numerical results are also good although we do not have the convergence
 theorem as Theorem 1 and Theorem 2. This shows that the approximate part of our shrinkage formula is also good,
 and that when the inner step in \cite{LHLL2013,LHS2013} is chose to be 5 the numerical experiments are convergent although they did not find a convergence control sequence.

\subsubsection{Experiments for the constrained TV OGS $L_1$ model}
In this section, we compare our methods (CATVOGSL1 and CITVOGSL1) with some other methods, such as Chan's method proposed in \cite{CTY2013} (Algorithm 2 in \cite{CTY2013} for the constrained TV-L1 model) and Liu's method proposed in \cite{LHLL2013}.

\begin{table}[!hbp]
\setlength{\abovecaptionskip}{0pt}
\setlength{\belowcaptionskip}{10pt} \centering \caption{Numerical comparison of Chan. \cite{CTY2013}, Liu. \cite{LHLL2013}, CATVOGSL1 and CITVOGSL1 for images (b)--(d) in Figure~\ref{originalimages}. PSNR: dB, Time: s, ``Is'' is short for Images.} \centering{\footnotesize{
\begin{tabular}{|@{~}c@{~}|@{~}c@{~}|@{~}r@{/}r@{/}l@{/}r@{/}l@{~}|@{~}r@{/}l@{/}r@{/}l@{~}|@{~}r@{/}l@{/}r@{/}l@{~}|@{~}r@{/}l@{/}r@{/}l@{~}|}
\hline 
  \multirow{2}{*}{Is} & Noise & \multicolumn{5}{@{}@{}c@{}@{}}{Chan. \cite{CTY2013}}& \multicolumn{4}{@{}c@{}}{Liu. \cite{LHLL2013}}&\multicolumn{4}{@{}c@{}}{CATVOGSL1} &\multicolumn{4}{c|}{CITVOGSL1}\\
 \cline{3-19}
     &                 level              & $\mu$&Itrs&PSNR&Time&ReE & Itrs&PSNR&Time&ReE& Itrs&PSNR&Time&ReE &Itrs&PSNR&Time&ReE\\
 \hline
  \multirow{3}{*}{(a)}&  30\% & 25&130&29.59&\textsf{61.57}&0.0796   &35&30.92&\textsf{30.23}&0.0683        &32&31.17&\textbf{17.97}&0.0663      &41&31.34&\textbf{21.56}&0.0651 \\
    & 40\%&18&99&28.85&\textsf{46.41}&0.0866    &37&30.11&\textsf{31.54}&0.0749         &32&30.10&\textbf{18.03}&0.0751      &37&30.06&\textbf{19.95}&0.0754 \\
        & 50\%&15&83&28.03&39.53&0.0953    &48&29.02&41.50&0.0850         &45&28.64&\textbf{24.38}&0.0888     &56&28.60&\textbf{29.69}&0.0892 \\
  \hline
  \multirow{3}{*}{(b)} &  30\% &25&128&29.71&\textsf{22.90}&0.0501     &35&31.81&\textsf{12.93}&0.0393       &26&31.77&\textbf{5.97}&0.0395     &29&31.75&\textbf{6.26}&0.0396 \\
       & 40\%&20&98&28.59&\textsf{17.80}&0.0570     & 40&30.70&\textsf{14.99}&0.0447       &26&30.47&\textbf{6.12}&0.0459     &29&30.22&\textbf{6.38}&0.0472  \\
        & 50\%&15&74&27.18&\textsf{13.81}&0.0670      &49&28.92&\textsf{17.63}&0.0549        &39&28.59&\textbf{9.27}&0.0570    &44&28.16&\textbf{9.64}&0.0599  \\
  \hline
  \multirow{3}{*}{(c)}&  30\% & 26&138&30.21&\textsf{25.65}&0.0607   &  35&32.15&\textsf{12.84}&0.0485     &25&32.30&\textbf{5.82}&0.0477    &30&32.08&\textbf{6.75}&0.0489 \\
    & 40\%&22&105&29.10&\textsf{19.27}&0.0689     &37&31.04&\textsf{13.67}&0.0551        &26&31.01&\textbf{5.76}&0.0553     &27&30.50&\textbf{6.19}&0.0587  \\
        & 50\%&15&76&27.80&\textsf{14.37}&0.0801      &43&29.28&\textsf{15.43}&0.0675        &35&28.85&\textbf{7.85}&0.0710     &42&28.40&\textbf{8.78}&0.0748 \\
         \hline
   \multirow{3}{*}{(d)} &  30\% &26&127&30.41&\textsf{23.17}&0.0975     &36&32.47&\textsf{13.15}&0.0769       &25&32.54&\textbf{5.68}&0.0764     &29&32.51&\textbf{6.75}&0.0766 \\
       & 40\%&20&95&29.44&\textsf{17.43}&0.1091      &39&31.43&\textsf{14.27}&0.0867        &28&31.35&\textbf{6.40}&0.0875     &32&31.06&\textbf{7.33}&0.0905  \\
        & 50\%&16&77&28.33&\textsf{14.29}&0.1239      &47&29.81&\textsf{17.05}&0.1045        &41&29.44&\textbf{8.83}&0.1116     &45&29.29&\textbf{9.98}&0.1148 \\
         \hline
\end{tabular}}}
\label{cameranmantol1CA}\end{table}

Similarly as the last section, we set the penalty parameters $\beta_1=80$, $\beta_2=2000$, $\beta_3=1$,for the ATV case, $\beta_1=80$, $\beta_2=2000$, $\beta_3=1$, for the ITV case
and relax parameter $\gamma=1.618$. The blur kernel is generated by Matlab built-in function \texttt{fspecial('gaussian',7,5)} for $7\times 7$ Gaussian blur with standard deviation 5. We first generated the blurred images operating on images (b)-(d) by the former Gaussian blur and further corrupted them by salt-and-pepper noise from 30\% to 50\%. We generate all noise effects by Matlab built-in function \texttt{imnoise(B,'salt \& pepper',level)} with ``\texttt{B}'' the blurred image and fix the same random matrix for different methods.

The numerical results by the three methods are shown in Table~\ref{cameranmantol1CA}. We have tuned the parameters manually to give the best PSNR improvement for Chan. \cite{CTY2013} as in Table~\ref{cameranmantol1CA} for different images. And for Liu. \cite{LHLL2013} we choose the given parameters $\mu$ default as 100,80,60 for 30\% to 50\% respectively. For our method CATVOGSL1, we set $\mu$ as 180,140,100 for 30\% to 50\% respectively. For our method CITVOGSL1, we set $\mu$ as 140,100,80 for 30\% to 50\% respectively. In experiments, we find that the parameters of our methods are robust and have wide rage to choose. Therefore, we set the same $\mu$ for different images.

\begin{figure}
  \centering
 \subfigure{\includegraphics[width=0.3\textwidth,clip]{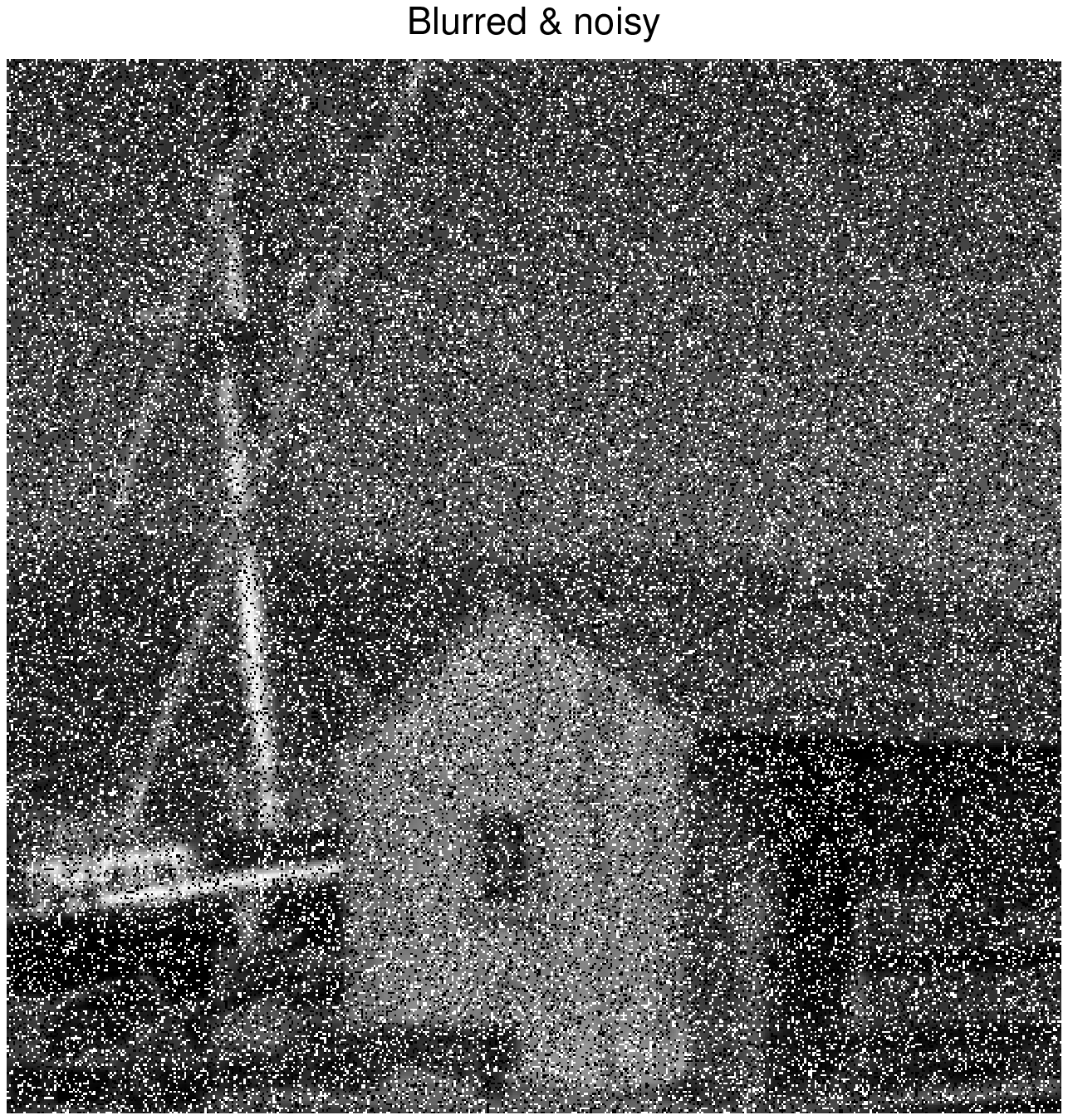}}
 \subfigure{\includegraphics[width=0.3\textwidth,clip]{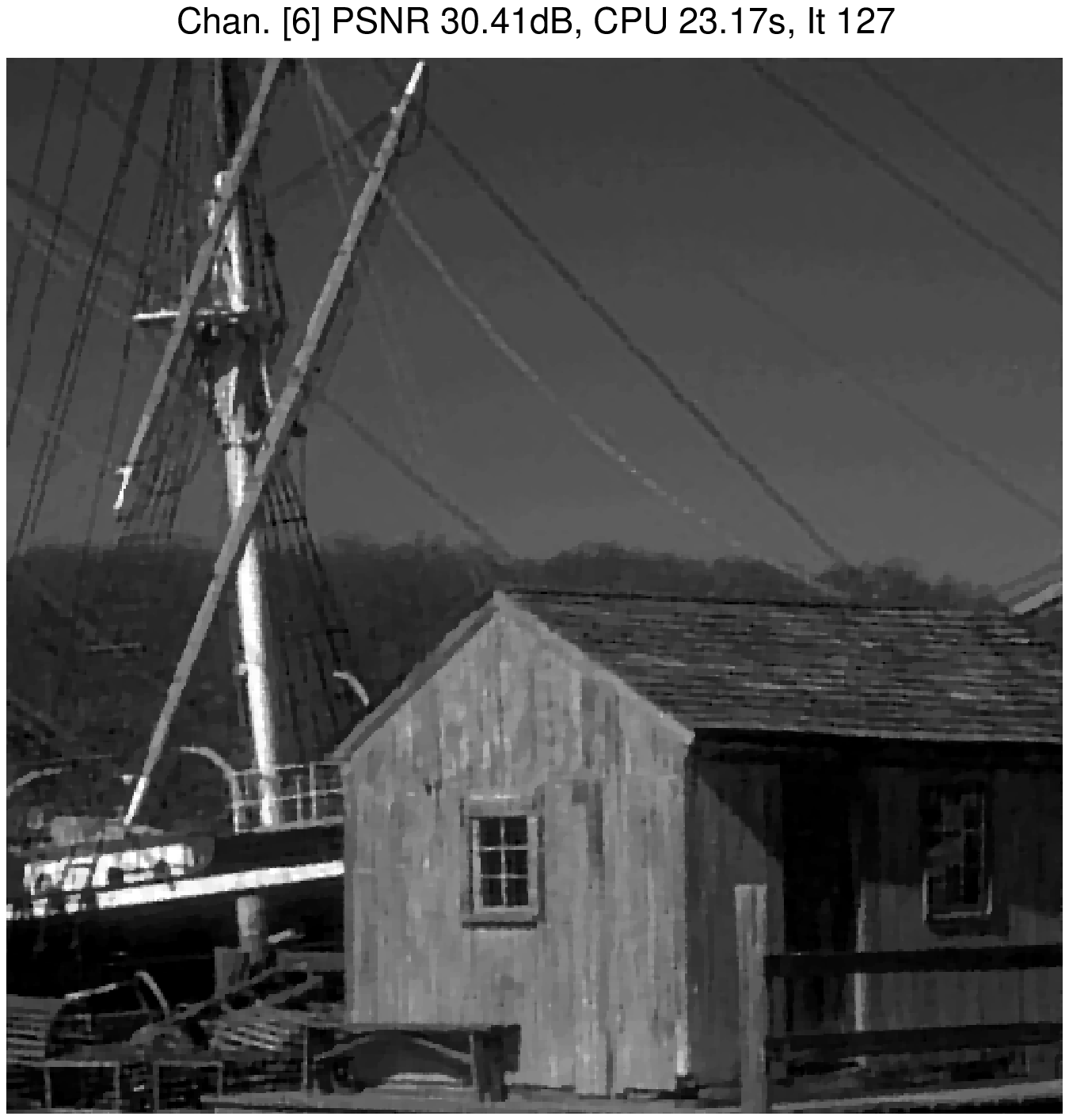}}
   \subfigure{\includegraphics[width=0.3\textwidth,clip]{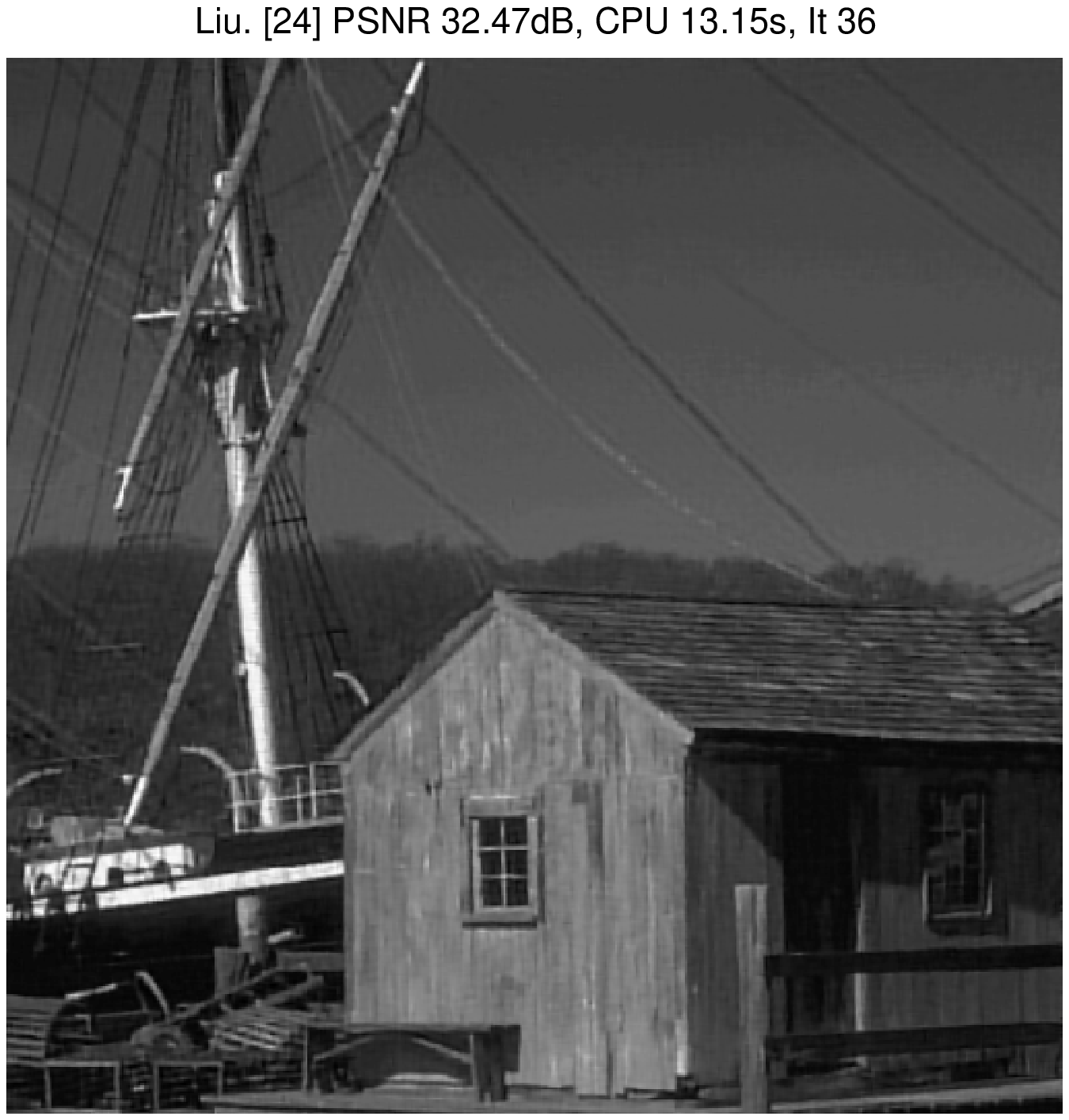}}
 \subfigure{\includegraphics[width=0.3\textwidth,clip]{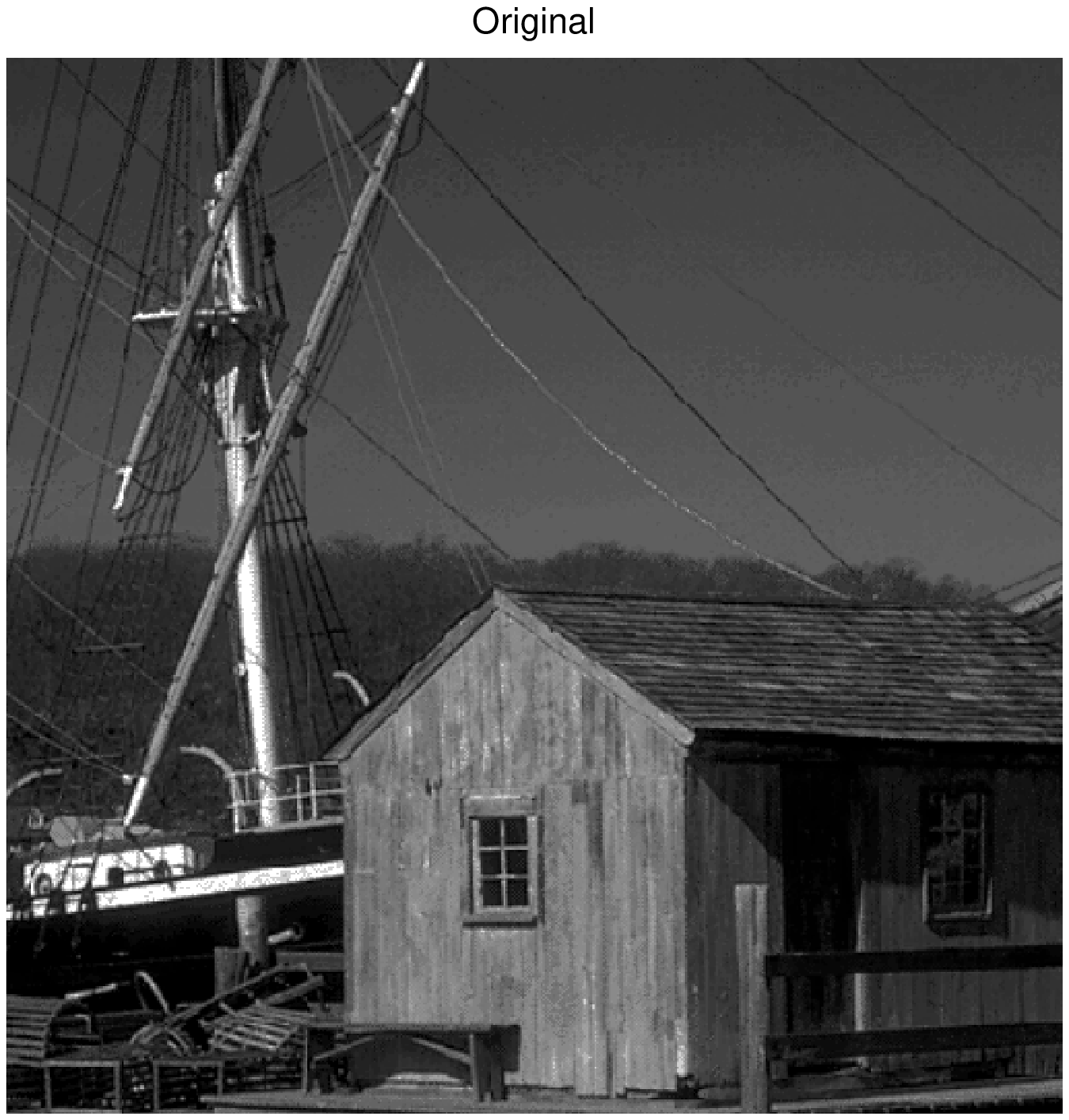}}
 \subfigure{\includegraphics[width=0.3\textwidth,clip]{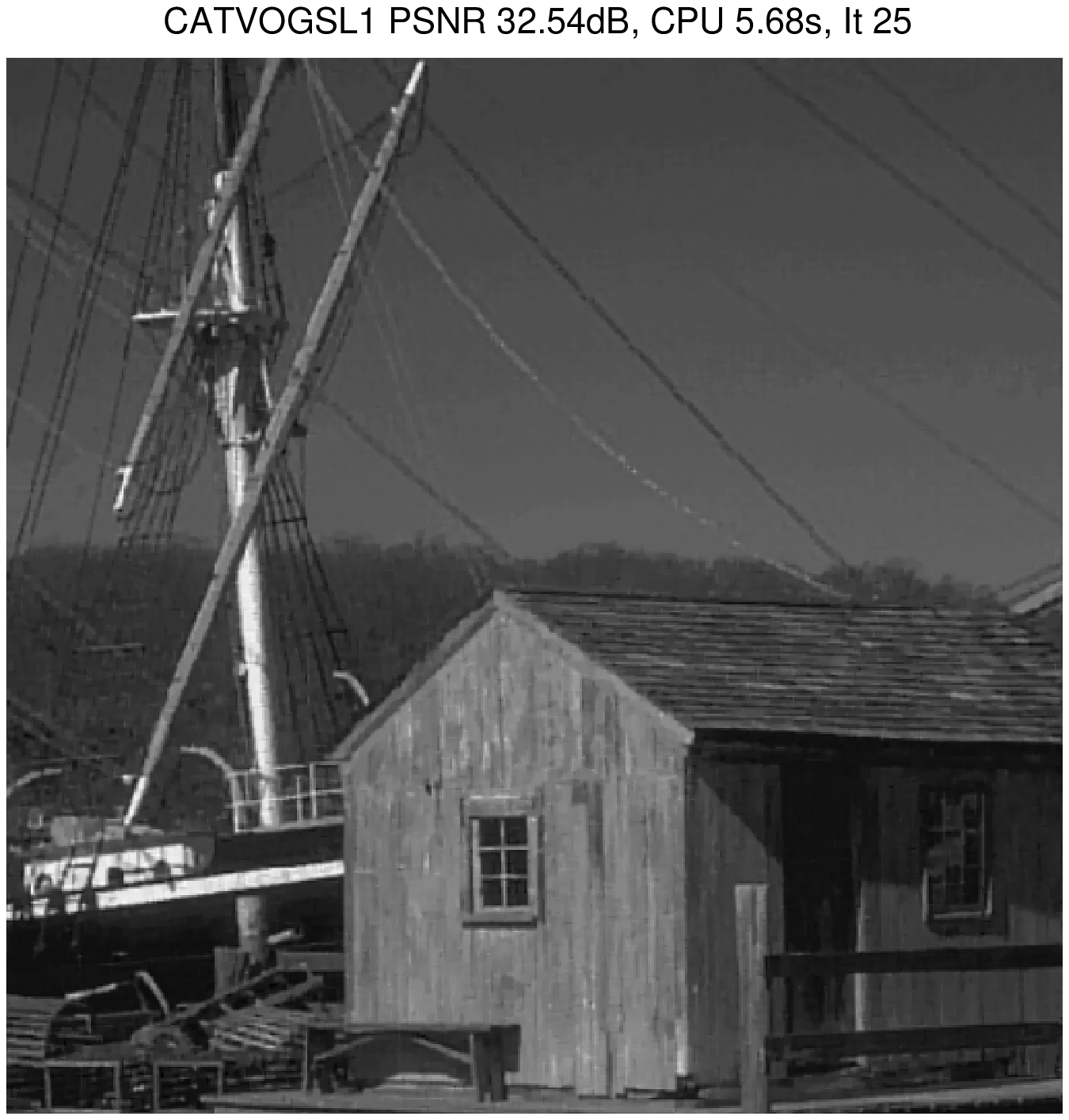}}
 \subfigure{\includegraphics[width=0.3\textwidth,clip]{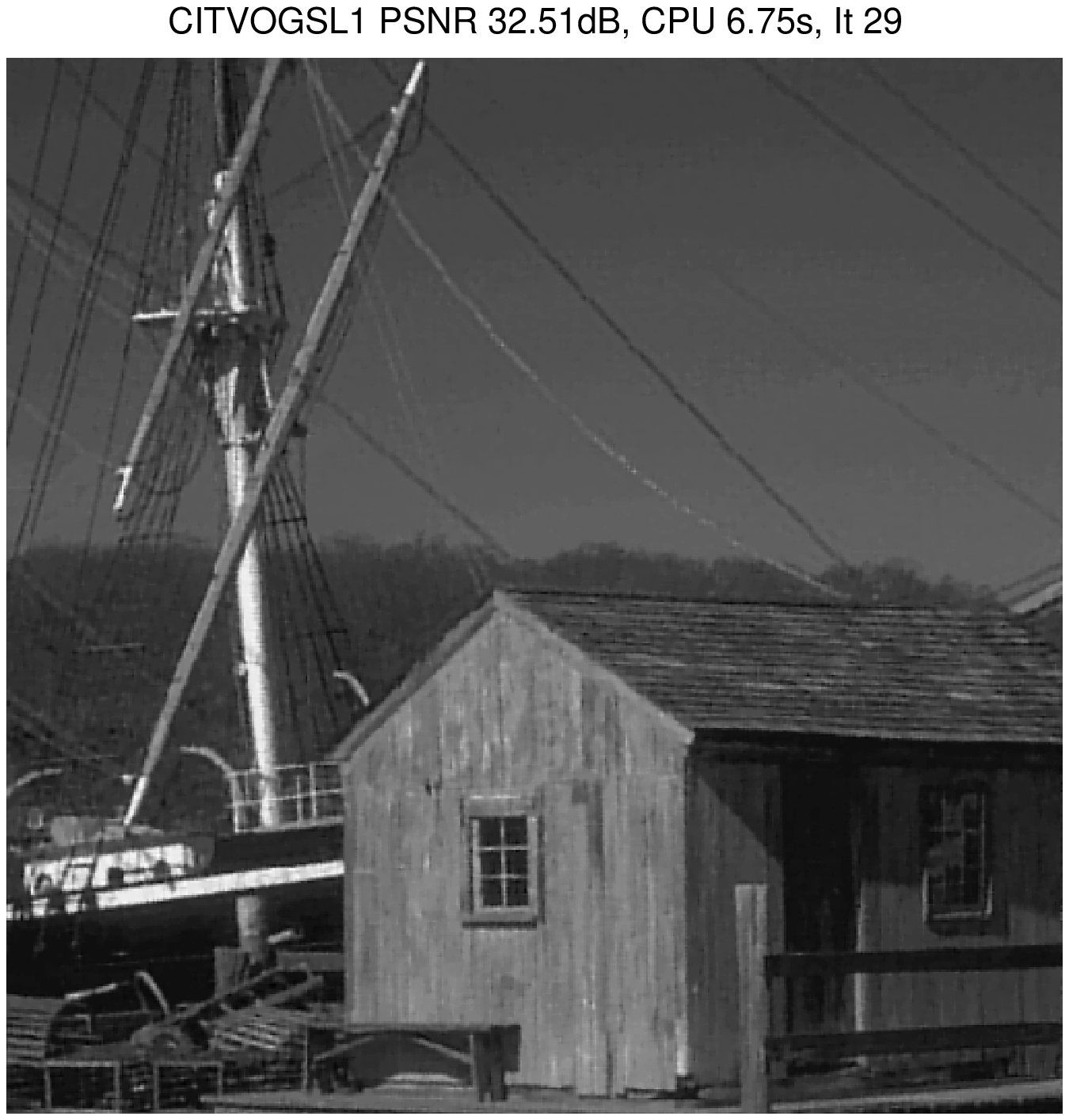}}
 \caption{Top row: blurred and noisy image (left), restoration images of Chan. \cite{CTY2013} (middle), Liu. \cite{LHLL2013}. Bottom row: original image (left), restoration images of CATVOGSL1, CITVOGSL1.}
  \label{conparisonofhousepole}
 \end{figure}

From Table~\ref{cameranmantol1CA}, we can also see that the PSNR and ReE of our methods (both ATV and ITV cases) are almost same as Liu. \cite{LHLL2013}, which used MM inner iterations to solve the subproblems (\ref{COTVSUB1}) and (\ref{COTVSUB11}) (only for the ATV case). However, each outer iteration of our methods is nearly twice faster than Liu. \cite{LHLL2013} from the experiments. The time of each outer iteration of our methods is almost the same as the tradition TV method in Chan. \cite{CTY2013}. Moreover, we can also see that sometimes the ATV is better than ITV and sometimes on the contrary for OGS TV. Finally, in Figure~\ref{conparisonofhousepole}, we display the degraded image, the original image and the restored images for 30\%
level of noise on Image (d) by the four methods. From the figure, we can see that both our methods and Liu. \cite{LHLL2013} can get better edges (handrail and window) than Chan. \cite{CTY2013}.

\section{Conclusion}

In this paper, we propose the explicit shrinkage formulas for one class of OGS regularization problems, which are with translation invariant overlapping groups. These formulas can be extended to several other regularization problems for instance nonconvex regularizers with overlapping group sparsity. We apply our results in OGS TV OGS regularization problems---deblurring and denoising problems, and get theorem convergence results and good experiments results. Furthermore, we also extend the image deblurring problems with OGS ATV in \cite{LHLL2013,LHS2013} to both ATV and ITV cases. For both of them, we have theorem convergence results by using ADMM both for $L_1$ and $L_2$ model. Since the theorem results and formulas are very simple, these results can be easily extended to many other application such as multichannel deconvolution and compress sensing, which we will consider in future. In addition, in this work we only choose all the entries of the weight matrix $W_g$ equal to 1. We will test for other weights in future on more experiments in order to choose the better or best weights for some applications.

\end{document}